\documentclass[12pt]{amsart}%
\usepackage{amsmath}
\usepackage{amsfonts}
\usepackage{amssymb}
\usepackage{amstext}
\usepackage{graphicx}%
\providecommand{\U}[1]{\protect\rule{.1in}{.1in}}
\providecommand{\U}[1]{\protect\rule{.1in}{.1in}}
\providecommand{\U}[1]{\protect\rule{.1in}{.1in}}
\newtheorem{theorem}{Theorem}[section]
\newtheorem{lemma}{Lemma}[section]
\newtheorem{proposition}{Proposition}[section]
\newtheorem{corollary}{Corollary}[section]

\newtheorem{definition}{Definition}[section]

\newtheorem{example}{Example}[section]
\numberwithin{equation}{section}

\theoremstyle{remark}
\newtheorem{remark}{Remark}[section]
\numberwithin{equation}{section}
\begin{document}
\title[Liouville properties for $p$-harmonic maps]{Liouville properties for $p$-harmonic maps with finite $q$-energy}
\author{Shu-Cheng Chang}
\address{Department of Mathematics, National Taiwan University, Taipei 10617, Taiwan, R.O.C.}
\email{scchang@math.ntu.edu.tw}
\author{Jui-Tang Chen}
\address{Department of Mathematics, National Taiwan Normal University, Taipei 11677,
Taiwan, R.O.C. }
\email{jtchen@ntnu.edu.tw }
\author{Shihshu Walter Wei}
\address{Department of Mathematics, University of Oklahoma, Norman, Oklahoma
73019-0315, U.S.A.}
\email{wwei@ou.edu}
\thanks{S.C. Chang and J.T. Chen were partially supported by NSC, and S.W. Wei was
partially supported by NSF(DMS-1447008) and the OU Arts and Sciences Travel
Assistance Program Fund.}
\date{}
\dedicatory{ }
\begin{abstract}
We introduce and study an approximate solution of the $p$-Laplace equation,
and a linearlization $\mathcal{L}_{\epsilon}$ of a perturbed $p$-Laplace
operator. By deriving an $\mathcal{L}_{\epsilon}$-type Bochner's formula and
Kato type inequalities, we prove a Liouville type theorem for weakly
$p$-harmonic functions with finite $p$-energy on a complete noncompact
manifold $M$ which supports a weighted Poincar\'{e} inequality and satisfies a
curvature assumption. This nonexistence result, when combined with an
existence theorem, yields in turn some information on topology, i.e. such an
$M$ has at most one $p$-hyperbolic end. Moreover, we prove a Liouville type
theorem for strongly $p$-harmonic functions with finite $q$-energy on
Riemannian manifolds. As an application, we extend this theorem to some
$p$-harmonic maps such as $p$-harmonic morphisms and conformal maps between
Riemannian manifolds. In particular, we obtain a Picard-type Theorem for
$p$-harmonic morphisms.

\end{abstract}
\maketitle

\section{Introduction}

The study of $p$-harmonic maps and in particular $p$-harmonic functions is
central to $p$-harmonic geometry and related problems.

A real-valued $C^{3}$ function on a Riemannian $m$-manifold $M$ with a
Riemannian metric $\langle\,\,,\,\rangle\,$ is said to be \textit{strongly
}$p$-\textit{harmonic} if $u$ is a (strong) solution of the $p$-Laplace
equation (\ref{1.0}), $p>1,$
\begin{equation}%
\begin{array}
[c]{l}%
\Delta_{p}u:=\text{\textrm{div}}\left(  |\nabla u|^{p-2}\nabla u\right)  =0.
\end{array}
\label{1.0}%
\end{equation}
where $\nabla u$ is the gradient vector field of $u$ on $M\,,$ and $|\nabla
u|=\langle\nabla u,\nabla u\rangle^{\frac{1}{2}}\,.$

A function $u\in W_{loc}^{1,p}\left(  M\right)  $ is said to be \textit{weakly
}$p$-\textit{harmonic} if $u$ is a (Sobolev) weak solution of the $p$-Laplace
equation (\ref{1.0}), i.e.
\[%
\begin{array}
[c]{l}%
\int_{M}\left\vert \nabla u\right\vert ^{p-2}\left\langle \nabla u,\nabla
\phi\right\rangle dv=0
\end{array}
\]
holds for every $\phi\in C_{0}^{\infty}\left(  M\right)  ,$ where $dv$ is the
volume element of $M\,.$

The $p$-Laplace equation (\ref{1.0}) arises as the Euler-Lagrange equation of
the $p$-energy $E_{p}$ functional given by $E_{p}(u) = \int_{M} |\nabla
u|^{p}\, dv\, . $ Ural'tseva \cite{U}, Evans \cite{E} and Uhlenbeck \cite{Ur}
proved that weak solutions of the equation (\ref{1.0}) have H\"{o}lder
continuous derivatives for $p\geq2$. Tolksdorff \cite{To}, Lewis \cite{Le} and
DiBenedetto \cite{D} extended the result to $p>1.$ In fact, weak solutions of
(\ref{1.0}), in general do not have any regularity better than $C_{loc}%
^{1,\alpha}.$

When $p=2,$ $p$-harmonic functions are simply harmonic functions. Liouville
type properties or topological end properties have been studied by a long list
of authors. We refer the reader to, for example \cite{L}, \cite{LT},
\cite{LT1}, \cite{LT2}, \cite{LT3}, \cite{LW1}, \cite{LW2}, \cite{LW3},
\cite{PRS}, \cite{SY} for further references. In particular, P. Li and J. Wang
showed Liouville type properties and splitting type properties on complete
noncompact manifolds with positive spectrum $\lambda$ when the Ricci curvature
has a lower bound depending on $\lambda.$ They also extended their work to a
complete noncompact manifold with weighted Poincar\'{e} inequality $(P_{\rho
}).$

For $p>1,$ We refer the works, for example \cite{CHS}, \cite{DW}, \cite{H},
\cite{H1}, \cite{H2}, \cite{H3}, \cite{HK}, \cite{HPV}, \cite{HV}, \cite{KN},
\cite{P}, \cite{PRS2}, \cite{W}, \cite{WLW}, to the reader. In particular, I.
Holopainen \cite{H2} proved a sharp $L^{q}$-Liouville properties for
$p$-harmonic functions, i.e. if $u\in L^{q}\left(  M\right)  $ is $p$-harmonic
(or more generally, $\mathcal{A}$-harmonic) in $M$ with $q>p-1,$ then $u$ is
constant. For $q=p-1$ and $m\geq2,$ there exist a complete Riemannian
$m$-manifold $M$ and a nonconstant positive $p$-harmonic function $f$ with
$\left\Vert f\right\Vert _{L^{p-1}\left(  M\right)  }<\infty.$ In
\cite{WLW}\cite{WLW2}, S.W.Wei, J.F. Li and L. Wu proved sharp Liouville
Theorems for $\mathcal{A}$-harmonic function $u$ with $p$-balanced growth
(e.g. $u \in L^{q}(M)\, ,$ for $q > p-1\, ,$ cf. \cite{W} 6.3). In \cite{HPV},
I. Holopainen, S. Pigola and G. Veronelli showed that if $u,v\in W_{loc}%
^{1,p}\left(  M\right)  \cap C^{0}\left(  M\right)  $ satisfy $\Delta_{p}%
u\geq\Delta_{p}v$ weakly and $\left\vert \nabla u\right\vert ,$ $\left\vert
\nabla v\right\vert \in L^{p}\left(  M\right)  ,$ for $p>1,$ then $u-v$ is
constant provided $M$ is connected, possibly incomplete, $p$-parabolic
Riemannian manifold. They also discussed $L^{q}$ comparison principles in the
non-parabolic setting. In \cite{PRS2}, S. Pigola, M. Rigoli and A.G. Setti
showed the constancy of $p$-harmonic map homotopic to a constant and with
finite $p$-energy from $p$-parabolic manifolds to manifolds with non-positive
sectional curvature. Moreover, if manifold $M$ has Poincar\'{e}-Sobolev
inequality, and $Ric_{M}\geq-k\left(  x\right)  $ with $k\left(  x\right)
\geq0$ and the integral type of $k\left(  x\right)  $ has upper bound
depending on Poincar\'{e}-Sobolev constant, $p\geq2$ and $p\geq q,$ then they
obtained constancy properties of $p$-harmonic map with some finite energy
types from $M$ to manifolds with non-positive sectional curvature. In
\cite{DW}, by a conservation law originated from E. Noether and comparison
theorems in Riemannian Geometry, Y.X. Dong and S.W. Wei obtained some
vanishing theorems for vector bundle valued differential forms. In particular,
they prove some Liouville type Theorems for $p$-harmonic maps with finite
$p$-energy under various curvature conditions.

In \cite{KN}, B. Kotschwar and L. Ni use a Bochner's formula on a neighborhood
of the maximum point (i.e. the $p$-Laplace operator is neither degenerate nor
singular elliptic on this neighborhood) to prove a gradient estimate for
positive $p$-harmonic functions. This also implies Liouville type properties
of positive $p$-harmonic functions on complete noncompact manifolds with
nonnegative Ricci curvature, and sectional curvature bounded below.

However, the approach of Kotschwar-Ni's gradient estimate for positive
$p$-harmonic functions, does not seem to work in this paper, since we need a
Bochner's formula which is unambiguously defined at every point in the manifold.

To overcome the difficulty, in this paper, we introduce and study an
approximate solution $u_{\epsilon}$ of the weakly $p$-harmonic function $u.$
This $u_{\epsilon}$ is the Euler-Lagrange equation of the $\left(
p,\epsilon\right)  $-energy%
\[%
\begin{array}
[c]{l}%
E_{p,\epsilon}=\int_{\Omega}\left(  \left\vert \nabla u_{\epsilon}\right\vert
^{2}+\epsilon\right)  ^{\frac{p}{2}}dv
\end{array}
\]
with $u-u_{\epsilon}\in W_{0}^{1,p}\left(  \Omega\right)  ,$ where $\Omega$ is
a domain on $M.$ That is, $u_{\epsilon}$ is the weak solution of a perturbed
$p$-Laplace equation
\begin{equation}%
\begin{array}
[c]{l}%
\Delta_{p,\epsilon}u_{\epsilon}=\text{\textrm{div}}\left(  \left(  \left\vert
\nabla u_{\epsilon}\right\vert ^{2}+\epsilon\right)  ^{\frac{p-2}{2}}\nabla
u_{\epsilon}\right)  =0.
\end{array}
\label{1.2}%
\end{equation}
Moreover, we consider a linearization $\mathcal{L}_{\epsilon}$ of the
perturbed operator $\Delta_{p,\epsilon}\,,$ given by
\begin{equation}%
\begin{array}
[c]{l}%
\mathcal{L}_{\epsilon}\left(  \Psi\right)  =\text{\textrm{div}}\left(
f_{\epsilon}^{p-2}A_{\epsilon}\left(  \nabla\Psi\right)  \right)  ,
\end{array}
\label{1.3}%
\end{equation}
for $\Psi\in C^{2}\left(  \Omega\right)  ,$ where $p>1,$ $f_{\epsilon}%
=\sqrt{\left\vert \nabla u_{\epsilon}\right\vert ^{2}+\epsilon}$ and%
\[%
\begin{array}
[c]{l}%
A_{\epsilon}:=\mathrm{id}+\left(  p-2\right)  \frac{\nabla u_{\epsilon}%
\otimes\nabla u_{\epsilon}}{f_{\epsilon}^{2}}.
\end{array}
\]

We observe that since $\Delta_{p,\epsilon}$ is no longer degenerate, by the
existence and $\epsilon$-Regularization results (Proposition \ref{ex} and
Proposition \ref{3.1}), $u_{\epsilon}$ exists and is infinitely
differentiable. Then we can derive an $\mathcal{L}_{\epsilon}$-type Bochner's
formula and a Kato type inequality, and apply them to $u_{\epsilon}.$ Hence,
using the convergence of the approximate solutions $u_{\epsilon}$ in $W^{1,p}$
on every domain in $M$, as $\epsilon\rightarrow0,$ we prove a Liouville type
property of weakly $p$-harmonic functions with finite $p$-energy. This
nonexistence result, when combined with the result of Proposition \ref{2 E},
yields in turn the topological information that such manifold has at most one
$p$-hyperbolic end.

We also note that, the perturbation method we employed in studying the
$p$-Laplace equation is in contrast to the methods in \cite{SU} for harmonic
maps on surfaces, in \cite{ES} for the level-set formulation of the mean
curvature flow, in \cite{HI} for the inverse mean curvature flow, and in
\cite{KN} for certain parabolic equations associated to the $p$-Laplacian.

\begin{theorem}
\label{T1}Let $M$ be a complete noncompact Riemannian $m$-manifold, $m\geq2$
supporting a weighted Poincar\'{e} inequality $\left(  P_{\rho}\right)  \,,$
with Ricci curvature
\begin{equation}%
\begin{array}
[c]{lll}%
Ric_{M}(x) & \geq & -\tau\rho\left(  x\right)
\end{array}
\label{1.4}%
\end{equation}
for all $x\in M,$ where $\tau$ is a constant such that%
\[%
\begin{array}
[c]{l}%
\tau<\frac{4\left(  p-1+\kappa\right)  }{p^{2}},\text{ }%
\end{array}
\]
in which $p>1,$ and%
\[%
\begin{array}
[c]{l}%
\kappa=\left\{
\begin{array}
[c]{ll}%
\max\left\{  \frac{1}{m-1},\min\left\{  \frac{\left(  p-1\right)  ^{2}}%
{m},1\right\}  \right\}  & \text{if }p>2,\\
\frac{\left(  p-1\right)  ^{2}}{m-1} & \text{if }1<p\leq2.
\end{array}
\right.
\end{array}
\]

Then every weakly $p$-harmonic function $u$ with finite $p$-energy $E_{p}$ is
constant. Moreover, $M$ has at most one $p$-hyperbolic end.
\end{theorem}

In Theorem \ref{T1}, we say that $M$ supports a weighted Poincar\'{e}
inequality $\left(  P_{\rho}\right)  $, if there exists a positive function
$\rho(x)$ a.e. on $M$ such that, for every $\Psi\in W_{0}^{1.2}\left(
M\right)  \,,$%

\begin{equation}%
\begin{array}
[c]{lll}%
\int_{M}\rho\left(  x\right)  \Psi^{2}\left(  x\right)  dv & \leq & \int
_{M}\left\vert \nabla\Psi\left(  x\right)  \right\vert ^{2}dv.
\end{array}
\label{WP}%
\end{equation}

If $\rho(x)$ is no less than a positive constant $\lambda\,,$ then $M$ has
positive spectrum. For example, the hyperbolic space $H^{m}$ has positive
spectrum, and $\rho\left(  x\right)  =\frac{\left(  m-1\right)  ^{2}}{4}.$ In
$\mathbb{R}^{m},$ if we select $\rho\left(  x\right)  =\frac{\left(
m-2\right)  ^{2}}{4|x|^{2}}\left(  x\right)  ,$ then (\ref{WP}) is Hardy's
inequality. For more examples, see \cite{CLW}\cite{LW3}\cite{WL}. \bigskip

\bigskip

If $u$ is a $C^{3}$ strongly $p$-harmonic function with finite $q$-energy,
then we have a Liouville type property as follows. \bigskip

\begin{theorem}
\label{T2}Let $M$ be a complete noncompact Riemannian $m$-manifold, $m\geq2, $
satisfying $\left(  P_{\rho}\right)  \,,$ with Ricci curvature
\begin{equation}%
\begin{array}
[c]{lll}%
Ric_{M}(x) & \geq & -\tau\rho(x)
\end{array}
\label{Rs}%
\end{equation}
for all $x\in M,$ where $\tau$ is a constant such that
\begin{equation}%
\begin{array}
[c]{l}%
\tau<\frac{4\left(  q-1+\kappa+b\right)  }{q^{2}},\text{ }%
\end{array}
\label{1.7}%
\end{equation}
in which
\[%
\begin{array}
[c]{l}%
\kappa=\min\{\frac{\left(  p-1\right)  ^{2}}{m-1},1\}\text{ and}\text{ }%
b=\min\{0,(p-2)(q-p)\},\text{ where }p>1.
\end{array}
\]
Let $u\in C^{3}\left(  M\right)  $ be a strongly $p$-harmonic function with
finite $q$-energy $E_{q}\left(  u\right)  <\infty.$ \newline(I). Then $u$ is
constant under each one of the following conditions: \newline(1) $p=2$ and
$q>\frac{m-2}{m-1},$ \newline(2) $p=4,$ $q>\max\left\{  1,1-\kappa-b\right\}
,$ \newline(3) $p>2,$ $p\neq4,$ and either%
\[%
\begin{array}
[c]{l}%
\max\left\{  1,p-1-\frac{\kappa}{p-1}\right\}  <q\leq\min\left\{
2,p-\frac{\left(  p-4\right)  ^{2}m}{4\left(  p-2\right)  }\right\}
\end{array}
\]
or%
\[%
\begin{array}
[c]{l}%
\max\left\{  2,1-\kappa-b\right\}  <q,
\end{array}
\]
\newline(II) $u$ does not exist for $1<p<2$ and $q>2.$
\end{theorem}

We remark that the curvature assumption (\ref{Rs}) or the assumption
(\ref{1.7}) on the constant $\tau$ in (\ref{Rs}) cannot be dropped, due to the
nontrivial $p$-harmonic functions with finite $q$-energy that are constructed
in Sect. 6.3.\medskip

As an application, we also extend Theorem \ref{T2} to $p$-harmonic morphisms
and conformal maps in Sections \ref{morphism} and \ref{Maps} respectively. In
particular, we obtain a Picard-type Theorem for $p$-harmonic morphisms. Some
applications to such Picard-type Theorems on stable minimal hypersurfaces in
Riemannian manifolds can be found in \cite{CLW}. \bigskip

The paper is organized as follows. In section $2$, we recall some facts about
$p$-hyperbolic and $p$-parabolic ends from \cite{LT} and \cite{H1}, and prove
an existence theorem on manifolds with two $p$-hyperbolic ends. In section
$3$, we introduce the linearization $\mathcal{L}_{\epsilon}$ (\ref{1.3}) of
the perturbed operator$\Delta_{p,\epsilon}\,,$ and derive the $\mathcal{L}%
_{\epsilon}$-type Bochner's formula (\ref{b}) and Kato type inequalities
(\ref{k'})(\ref{k''}) for the solution $u_{\epsilon}$ of the perturbed
equation (\ref{1.2}). In section $4$, by applying Bochner's formula and Kato's
inequality, we show a Liouville type theorem and one $p$-hyperbolic end
property for a weakly $p$-harmonic function with finite $p$-energy in a
complete noncompact manifold which supports a weighted Poincar\'{e} inequality
and satisfies a curvature assumption. In section $5$, we show Liouville type
theorems for strongly $p$-harmonic functions with finite $q$-energy, and we
also extend our results to some $p$-harmonic maps such as $p$-harmonic
morphisms and conformal maps between Riemannian manifolds. In section $6$ of
the Appendix, we prove the existence of the approximate solution $u_{\epsilon
},$ Proposition \ref{3.1}, and volume estimate of complete noncompact
manifolds with $p$-Poincar\'{e} inequality. We also construct an example of
non-trivial $p$-harmonic function with finite $q$-energy on manifolds with
weighted Poincar\'{e} inequality.

\section{$p$-Hyperbolicity}

We recall some basic facts about capacities from \cite{H}, \cite{H1} and
\cite{Tr}.

Let $M$ be a Riemannian manifold, $G\subset M$ a connected open set in $M.$ If
$D$ and $\Omega$ are nonempty, disjoint, and closed sets contained in the
closure of $G.$ A triple $\left(  \Omega,D;G\right)  $ is called a condenser.
The $p$-capacity of $\left(  \Omega,D;G\right)  $ is defined by%
\[%
\begin{array}
[c]{lll}%
\emph{Cap}_{p}\left(  \Omega,D;G\right)  & = & \displaystyle\inf\limits_{u}
\int_{G}\left\vert \nabla u\right\vert ^{p}dv,
\end{array}
\]
for $1\leq p<\infty\, ,$ where the infimum is taken over all $u\in
W^{1,p}\left(  G\right)  \cap C^{0}(G) $ with $u=1$ in $\Omega$ and $u=0$ in
$D.$

Above and in what follows, $W^{1,p}\left(  M\right)  $ is the Sobolev space of
all function \thinspace$u\in L^{p}\left(  M\right)  $ and whose distributional
gradient $\nabla u$ also belongs to $L^{p}\left(  M\right)  ,$ with respect to
the Sobolev norm
\[%
\begin{array}
[c]{lll}%
\left\Vert u\right\Vert _{1,p} & = & \left\Vert u\right\Vert _{L_{p}%
}+\left\Vert \nabla u\right\Vert _{L_{p}}.
\end{array}
\]
The space $W_{0}^{1,p}\left(  M\right)  $ is the closure of $C_{0}^{\infty
}\left(  M\right)  $ in $W^{1,p}\left(  M\right)  \, ,$ with respect to the
$\left\Vert \quad\right\Vert _{1,p}$ norm.

The following properties of capacities are well known (see e.g. \cite{Tr}).

\begin{itemize}
\item $\Omega_{2}\subset\Omega_{1}\Longrightarrow\emph{Cap}_{p}\left(
\Omega_{2},D;G\right)  \leq\emph{Cap}_{p}\left(  \Omega_{1},D;G\right)  ;$

\item $D_{2}\subset D_{1}\Longrightarrow\emph{Cap}_{p}\left(  \Omega
,D_{2};G\right)  \leq\emph{Cap}_{p}\left(  \Omega,D_{1};G\right)  ;$

\item If $\Omega_{1}\supset\Omega_{2}\cdots\supset\cap_{i}\Omega_{i}=\Omega$
and $D_{1}\supset D_{2}\cdots\supset\cap_{i}D_{i}=D,$ then
\[
\emph{Cap}_{p}\left(  \Omega,D;G\right)  =\lim_{i\rightarrow\infty}%
\emph{Cap}_{p}\left(  \Omega_{i},D_{i};G\right)  .
\]

\item If $\overline{G\backslash\left(  \Omega\cup D\right)  }$ is compact,
then there exists a unique weak solution $u:\overline{G\backslash\left(
\Omega\cup D\right)  }\rightarrow\mathbb{R} $ to the Dirichlet problem%
\[%
\begin{array}
[c]{l}%
\left\{
\begin{array}
[c]{ll}%
\Delta_{p}u=0\text{ \ \ \ \ \ } & \text{on }G\backslash\left(  \Omega\cup
D\right)  ,\text{ }\\
u=1 & \text{on }\Omega,\\
u=0 & \text{on }D,
\end{array}
\right.
\end{array}
\]
with $\emph{Cap}_{p}\left(  \Omega,D;G\right)  =\int_{G}\left\vert
du\right\vert ^{p}dv.$
\end{itemize}

Given a compact set $\Omega$ in $M$, an \textit{end} $E_{\Omega}$ with respect
to $\Omega$ is an unbounded connected component of $M\backslash\Omega\, .$ By
a compactness argument, it is readily seen that the number of ends with
respect to $\Omega$ is finite, it is also clear that if $\Omega\subset
\Omega^{\prime}$ , then every end $E_{\Omega^{\prime}}$ is contained in
$E_{\Omega}$, so that the number of ends increases as the compact $\Omega$
enlarges. Let $x_{0}\in\Omega.$ We denote $E_{\Omega}\left(  R\right)
=B_{x_{0}}\left(  R\right)  \cap E_{\Omega},$ $\partial E_{\Omega}\left(
R\right)  =\partial B_{x_{0}}\left(  R\right)  \cap E_{\Omega}$ and $\partial
E_{\Omega}=\partial\Omega\cap E_{\Omega}.$

In \cite{LT} (or see e.g. \cite{L}, \cite{LT1}-\cite{LT3}, \cite{PRS}),
$2$-parabolic and $2$-nonparabolic manifolds and ends are introduced. In
\cite{H1}, I. Holopainen defined the $p$-parabolic end as follows:

\begin{definition}
Let $E$ be an end of $M$ with respect to $\Omega.$ $E$ is $p$-parabolic, or,
equivalently, has zero $p$-capacity at infinity if,%
\[%
\begin{array}
[c]{lllll}%
\emph{Cap}_{p}\left(  \Omega,\infty;E\right)  & := & \lim_{i\rightarrow\infty
}\emph{Cap}_{p}\left(  \Omega,\overline{E}\backslash\Omega_{i};E\right)  & = &
0,
\end{array}
\]
where $\left\{  \Omega_{i}\right\}  _{i=1}^{\infty}$ is an exhaustion of $M$
by relatively compact open domains with smooth boundary and $\Omega_{i}
\subset\subset\Omega_{i+1}, $ for every integer $i\, .$
\end{definition}

This definition also implies: if $E$ is an end with respect to $\Omega,$ there
are sequence of weakly $p$-harmonic functions $\left\{  u_{i}\right\}  ,$
$u_{i}\in W^{1,p},$ defined on $E,$ satisfying%
\begin{equation}%
\begin{array}
[c]{lllll}%
\Delta_{p}u_{i} & = & 0 & \text{on} & E\left(  r_{i}\right)
\end{array}
\label{p-1}%
\end{equation}
with boundary conditions%
\begin{equation}%
\begin{array}
[c]{lll}%
u_{i} & = & \left\{
\begin{array}
[c]{ll}%
1\text{ \ \ \ \ \ \ \ \ } & \text{on }\Omega,\\
0 & \text{on }\overline{E\backslash\Omega_{i}},
\end{array}
\right.
\end{array}
\label{p-2}%
\end{equation}
then $\left\{  u_{i}\right\}  $ converges (converges uniformly on each compact
set of $E$) to the constant function $u=1$ on $E$ as $i\rightarrow\infty.$

\begin{definition}
An end $E$ is $p$-hyperbolic (or $p$-nonparabolic) if $E$ is not $p$-parabolic.
\end{definition}

If $h_{i}$ is a weakly $p$-harmonic function satisfying (\ref{p-1}) and
(\ref{p-2}), then $E$ is $p$-hyperbolic if and only if $\left\{
h_{i}\right\}  $ converges to a weakly $p$-harmonic function $h$ with $h=1$ on
$\partial E$, $\inf_{E}h=0$ and finite $p$-energy.

\begin{definition}
A manifold $M$ is $p$-parabolic, or, equivalently, has zero $p$-capacity at
infinity if, for each compact set $\Omega\subset M,$
\[%
\begin{array}
[c]{lllll}%
\emph{Cap}_{p}\left(  \Omega,\infty;M\right)  & := & \lim_{i\rightarrow\infty
}\emph{Cap}_{p}\left(  \Omega,M\backslash\Omega_{i};M\right)  & = & 0,
\end{array}
\]
where $\left\{  \Omega_{i}\right\}  _{i=1}^{\infty}$ is an exhaustion of $M$
by domains with smooth boundary and $\Omega_{i} \subset\subset\Omega_{i+1}, $
for every integer $i\, .$
\end{definition}

\begin{definition}
A manifold $M$ is $p$-hyperbolic (or $p$-nonparabolic) if $M$ is not $p$-parabolic.
\end{definition}

This definition also implies that a manifold $M$ is $p$-parabolic if each end
of $M$ is $p$-parabolic, $M$ is $p$-hyperbolic if $M$ has at least one $p
$-hyperbolic end.

Now we focus on manifold $M$ with two $p$-hyperbolic ends (cf. \cite{H}).

\begin{proposition}
\label{2 E}Let $M$ be a complete noncompact manifold, and assume $M$ has two
$p$-hyperbolic ends $E_{1}$ and $E_{2}.$ Then there exists a weakly
$p$-harmonic function $h:M\rightarrow\mathbb{R}$ with finite $p$-energy such
that $0<h<1,$ $\sup_{E_{1}}h=1$ and $\inf_{E_{2}}h=0.$ Moreover, $h$ is
$C^{1,\alpha}$.
\end{proposition}

\proof

Given $\Omega\subset M,$ we fix an exhaustion $\left\{  \Omega_{i}\right\}  $
of $M$ by domains with smooth boundary and $\Omega_{i}\subset\subset
\Omega_{i+1}$ for every integer $i\, .$

Denote by $E_{A}$ the $p$-hyperbolic ends of $M$ with respect to $A\,.$ For
every $A$, let $u_{i}^{E_{A}}$ be the $p$-harmonic function satisfying%
\[
\left\{
\begin{array}
[c]{ll}%
\Delta_{p}u_{i}^{E_{A}}=0 & \text{in }E_{A}\cap\Omega_{i},\\
u_{i}^{E_{A}}=1 & \text{on }\partial E_{A},\\
u_{i}^{E_{A}}=0 & \text{on }\partial E_{A,\Omega_{i}}=\partial\left(
E_{A}\cap\Omega_{i}\right)  \backslash\partial E_{A}.
\end{array}
\right.
\]
By the monotone property, $u_{i}^{E_{A}}$ converges uniformly to $u^{E_{A}}$
on every compact subset of $E_{A}.$

For every $i$, let $h_{i}$ be the weak solution of the boundary value problem
\[
\left\{
\begin{array}
[c]{ll}%
\Delta_{p}h_{i}=0 & \text{in }\Omega_{i},\\
h_{i}=1 & \text{on }\partial\Omega_{i}\cap E_{1},\\
h_{i}=0 & \text{on }\partial\Omega_{i}\cap\left(  M\backslash E_{1}\right)  .
\end{array}
\right.
\]
Then, $0\leq h_{i}\leq1$, and by gradient estimate (\cite{KN}), there are
subsequence, say $\left\{  h_{i}\right\}  ,$ converges, locally uniformly, to
a weakly $p$-harmonic function $h$ on $M$, satisfying $0\leq h\leq1.$

On $E_{1}$, the maximum principle implies $1-u_{i}^{E_{1}}\leq h_{i}<1.$ Hence
$1-u^{E_{1}}\leq h<1$ on $E_{1}$, so that $\sup_{E_{1}}\left(  1-u^{E_{1}%
}\right)  \leq\sup_{E_{A}}h=1$ gives $\sup_{E_{1}}h=1$ since $\inf_{E_{1}%
}u^{E_{1}}=0.$

On $E_{2},$ the maximum principle implies $0<h_{i}\leq u_{i}^{E_{2}}.$ Hence
we have $0<h\leq u^{E_{2}}$ on $E_{2}$, so that $0\leq\inf_{E_{2}}h\leq
\inf_{E_{2}}u^{E_{2}}=0$.

Now we have $\sup_{E_{1}}h=1$ and $\inf_{E_{2}}h=0,$ so $h$ is a nonconstant
$p$-harmonic function on $M.$

Finally, $h$ has finite $p$-energy by
\[%
\begin{array}
[c]{lll}%
\emph{Cap}_{p}\left(  E_{1}\backslash\Omega_{i},M\backslash\left(  \Omega
_{i}\cup E_{1}\right)  ;M\right)  & = & \int_{M}\left\vert \nabla
h_{i}\right\vert ^{p}dv\neq0,
\end{array}
\]
and the monotonic properties of capacities. \endproof

\bigskip

\section{Bochner's formula and Kato's inequality}

First of all, we define $N=M\times%
\mathbb{R}
$ with metric $g_{N}=g_{M}+dt^{2},$ and let
\begin{equation}
v_{\epsilon}\left(  x,t\right)  =u_{\epsilon}\left(  x\right)  +\sqrt
{\epsilon}t \label{1.2.1}%
\end{equation}
for $x\in\Omega\subset M,$ $t\in%
\mathbb{R}
,$ and $\epsilon>0,$ where $u_{\epsilon}$ is the solution of the perturbed
$p$-Laplace equation (\ref{1.2}). Then $v_{\epsilon}\in C^{\infty}$ is a
strongly $p$-harmonic function on $\Omega_{N}=\Omega\times%
\mathbb{R}
,$ i.e. if $\Delta_{p}^{N}$ is the $p$-Laplace operator on $\left(  \Omega
_{N},g_{N}\right)  ,$ we have $\Delta_{p}^{N}v_{\epsilon}=0$ with $\left\vert
\nabla^{N}v_{\epsilon}\right\vert ^{2}\ge\epsilon>0$ and $Ric_{N}\left(
\nabla^{N}v_{\epsilon},\nabla^{N}v_{\epsilon}\right)  =Ric\left(  \nabla
u_{\epsilon},\nabla u_{\epsilon}\right)  .$ Moreover, if $f=\left\vert \nabla
u_{\epsilon}\right\vert ,$ then $f_{\epsilon}=\left\vert \nabla^{N}%
v_{\epsilon}\right\vert =\sqrt{f^{2}+\epsilon}$ which is independent of $t.$
Hence, we have $\nabla^{N}f_{\epsilon}=\nabla f_{\epsilon}$ and $\Delta
^{N}f_{\epsilon}=\Delta f_{\epsilon}.$

According to the argument of Kotschwar-Ni \cite{KN}, we define the linearized
operator $\mathcal{L}_{0}^{N}$ of the $p$-Laplace operator $\Delta_{p}^{N}$ on
$\left(  \Omega_{N},g_{N}\right)  $ as follows:%
\[
\mathcal{L}_{0}^{N}\left(  \Psi\right)  =\text{\textrm{div}}^{N}\left(
f_{\epsilon}^{p-2}A_{0}\left(  \nabla^{N}\Psi\right)  \right)  ,
\]
for $\Psi\in C^{2}\left(  \Omega_{N}\right)  ,$ where $\operatorname{div}^{N}$
is the divergence on $\left(  \Omega_{N},g_{N}\right)  $ and%
\[%
\begin{array}
[c]{l}%
A_{0}:=\mathrm{id}+\left(  p-2\right)  \frac{\nabla^{N}v_{\epsilon}%
\otimes\nabla^{N}v_{\epsilon}}{f_{\varepsilon}^{2}}.
\end{array}
\]

Now we show Bochner's formula as the following:

\begin{lemma}
\label{Bo0} Let $v_{\epsilon}$ be the $p$-harmonic function on $\left(
\Omega_{N},g_{N}\right)  ,$ and $\left(  \nabla d\right)  ^{N} v_{\epsilon}$
be the Hessian of $v_{\epsilon}$ on $\left(  \Omega_{N},g_{N}\right)  \,.$
Then for every $p>1,$%
\begin{equation}%
\begin{array}
[c]{lll}%
\frac{1}{2}\mathcal{L}_{0}^{N}\left(  f_{\epsilon}^{2}\right)  & = &
\frac{p-2}{4}f_{\epsilon}^{p-4}\left\vert \nabla^{N}f_{\epsilon}%
^{2}\right\vert ^{2}+f_{\epsilon}^{p-2}\left(  \left\vert \left(  \nabla
d\right)  ^{N}v_{\epsilon}\right\vert ^{2}+Ric_{N}\left(  \nabla
^{N}v_{\epsilon},\nabla^{N}v_{\epsilon}\right)  \right)  .
\end{array}
\label{b00}%
\end{equation}

\end{lemma}

\proof Since $f_{\epsilon}>0\,,$ for every $p>1\,,$ the $p$-harmonic equation
$\Delta_{p}^{N} v_{\epsilon} = 0$ is equivalent to
\begin{equation}%
\begin{array}
[c]{l}%
\frac{p-2}{2}\left\langle \nabla^{N}f_{\epsilon}^{2},\nabla^{N}v_{\epsilon
}\right\rangle =-f_{\epsilon}^{2}\Delta^{N}v_{\epsilon}%
\end{array}
\label{b0}%
\end{equation}

which implies
\begin{equation}%
\begin{array}
[c]{lll}%
\frac{p-2}{2}f_{\epsilon}^{p-6}\left\langle \nabla^{N}v_{\epsilon},\nabla
^{N}f_{\epsilon}^{2}\right\rangle ^{2} & = & -f_{\epsilon}^{p-4}\left\langle
\nabla^{N}v_{\epsilon},\nabla^{N}f_{\epsilon}^{2}\right\rangle \Delta
^{N}v_{\epsilon}.
\end{array}
\label{b4}%
\end{equation}

On the other hand, taking the gradient of both sides of (\ref{b0}), and then
taking the inner product with $\nabla^{N}v_{\epsilon}\,,$ we have
\begin{equation}%
\begin{array}
[c]{lll}%
\frac{p-2}{2}\left\langle \nabla^{N}\left\langle \nabla^{N}f_{\epsilon}%
^{2},\nabla^{N}v_{\epsilon}\right\rangle ,\nabla^{N}v_{\epsilon}\right\rangle
& = & -\left\langle \nabla^{N}f_{\epsilon}^{2},\nabla^{N}v_{\epsilon
}\right\rangle \Delta^{N}v_{\epsilon}\\
&  & -f_{\epsilon}^{2}\left\langle \nabla^{N}\left(  \Delta^{N}v_{\epsilon
}\right)  ,\nabla^{N}v_{\epsilon}\right\rangle .
\end{array}
\label{b1}%
\end{equation}

Now we compute
\begin{equation}%
\begin{array}
[c]{ll}%
\frac{1}{2}\mathcal{L}_{0}^{N}\left(  f_{\epsilon}^{2}\right)  & = \frac{1}%
{2}\text{\textrm{div}}^{N}\left(  f_{\epsilon}^{p-2}\nabla^{N}f_{\epsilon}%
^{2}+\left(  p-2\right)  f_{\epsilon}^{p-4}\left\langle \nabla^{N}v_{\epsilon
},\nabla^{N}f_{\epsilon}^{2}\right\rangle \nabla^{N}v_{\epsilon}\right) \\
& = \frac{p-2}{4}f_{\epsilon}^{p-4}\left\vert \nabla^{N}f_{\epsilon}%
^{2}\right\vert ^{2}+\frac{1}{2}f_{\epsilon}^{p-2}\Delta^{N}f_{\epsilon}%
^{2}+\frac{\left(  p-2\right)  \left(  p-4\right)  }{4}f_{\epsilon}%
^{p-6}\left\langle \nabla^{N}v_{\epsilon},\nabla^{N}f_{\epsilon}%
^{2}\right\rangle ^{2}\\
& \quad+\frac{p-2}{2}f_{\epsilon}^{p-4}\left\langle \nabla^{N}\left\langle
\nabla^{N}v_{\epsilon},\nabla^{N}f_{\epsilon}^{2}\right\rangle ,\nabla
^{N}v_{\epsilon}\right\rangle \\
& \quad+\frac{p-2}{2}f_{\epsilon}^{p-4}\left\langle \nabla^{N}v_{\epsilon
},\nabla^{N}f_{\epsilon}^{2}\right\rangle \Delta^{N}v_{\epsilon}.
\end{array}
\label{b2}%
\end{equation}
Substituting (\ref{b1}) into (\ref{b2}), one gets%
\begin{equation}%
\begin{array}
[c]{ll}%
\frac{1}{2}\mathcal{L}_{0}^{N}\left(  f_{\epsilon}^{2}\right)  = & \frac
{p-2}{4}f_{\epsilon}^{p-4}\left\vert \nabla^{N}f_{\epsilon}^{2}\right\vert
^{2}+\frac{p-4}{2}f_{\epsilon}^{p-4}\left\langle \nabla^{N}v_{\epsilon}%
,\nabla^{N}f_{\epsilon}^{2}\right\rangle \Delta^{N}v_{\epsilon}\\
& +\frac{1}{2}f_{\epsilon}^{p-2}\Delta^{N}f_{\epsilon}^{2} - f_{\epsilon
}^{p-2}\left\langle \nabla^{N}\left(  \Delta^{N}v_{\epsilon}\right)
,\nabla^{N}v_{\epsilon}\right\rangle \\
& +\frac{\left(  p-2\right)  \left(  p-4\right)  }{4}f_{\epsilon}%
^{p-6}\left\langle \nabla^{N}v_{\epsilon},\nabla^{N}f_{\epsilon}%
^{2}\right\rangle ^{2}.
\end{array}
\label{b3}%
\end{equation}

Applying Bochner's formula%
\[%
\begin{array}
[c]{l}%
\frac{1}{2}\Delta^{N}f_{\epsilon}^{2}=\left\vert \left(  \nabla d\right)
^{N}v_{\epsilon}\right\vert ^{2}+\left\langle \nabla^{N}v_{\epsilon}%
,\nabla^{N}\left(  \Delta^{N}v_{\epsilon}\right)  \right\rangle +Ric_{N}%
\left(  \nabla^{N}v_{\epsilon},\nabla^{N}v_{\epsilon}\right)
\end{array}
\]
and the equation (\ref{b4}) to the third term and the last term of right hand
side of (\ref{b3}) respectively, one obtains the desired formula (\ref{b00}).
\endproof

If $\Psi$ is independent of $t,$ then%
\[%
\begin{array}
[c]{lll}%
\mathcal{L}_{0}^{N}\left(  \Psi\right)  & = & \text{\textrm{div}}^{N}\left(
f_{\epsilon}^{p-2}\nabla^{N}\Psi+\left(  p-2\right)  \left\langle \nabla
^{N}v_{\epsilon},\nabla^{N}\Psi\right\rangle \nabla^{N}v_{\epsilon}\right) \\
& = & \text{\textrm{div}}\left(  f_{\epsilon}^{p-2}\nabla\Psi+\left(
p-2\right)  f_{\epsilon}^{p-4}\left\langle \nabla u_{\epsilon},\nabla
\Psi\right\rangle \nabla u_{\epsilon}\right) \\
& = & \mathcal{L}_{\epsilon}\left(  \Psi\right)
\end{array}
\]
where $\mathcal{L}_{\epsilon}$ is defined by (\ref{1.3}). Hence Lemma
\ref{Bo0} implies the following Lemma.

\begin{lemma}
\label{Bo} Let $u_{\epsilon}$ be a solution of (\ref{1.2}) on $\Omega\subset
M\,,$ $f_{\epsilon}=\sqrt{\left\vert \nabla u_{\epsilon}\right\vert
^{2}+\epsilon}\,,$ and $\nabla du_{\epsilon}$ be the Hessian of $u_{\epsilon
}\, $ on $M\, .$ Then for every $p>1,$%
\begin{equation}%
\begin{array}
[c]{lll}%
\frac{1}{2}\mathcal{L}_{\epsilon}\left(  f_{\epsilon}^{2}\right)  & = &
\frac{p-2}{4}f_{\epsilon}^{p-4}\left\vert \nabla f_{\epsilon}^{2}\right\vert
^{2}+f_{\epsilon}^{p-2}\left(  \left\vert \nabla du_{\epsilon}\right\vert
^{2}+Ric\left(  \nabla u_{\epsilon},\nabla u_{\epsilon}\right)  \right)  .
\end{array}
\label{b}%
\end{equation}

\end{lemma}

\bigskip

Next, we derive the following Kato type inequalities for the approximate
solution $u_{\epsilon}:$

\begin{lemma}
\label{ka1}Let $u_{\epsilon}$ be a solution of (\ref{1.2}) on $\Omega\subset
M^{m},$ $p>1.$ Then the Hessian of $u_{\epsilon}$ satisfies
\begin{equation}%
\begin{array}
[c]{lll}%
\left\vert du_{\epsilon}\right\vert ^{2}\left\vert \nabla du_{\epsilon
}\right\vert ^{2} & \geq & \frac{1+\kappa_{1}}{4}\left\vert \nabla\left\vert
du_{\epsilon}\right\vert ^{2}\right\vert ^{2}%
\end{array}
\label{k'}%
\end{equation}
at $x\in\Omega,$ where%
\[%
\begin{array}
[c]{lll}%
\kappa_{1} & = & \left\{
\begin{array}
[c]{ll}%
\frac{1}{m-1} & \text{if }p\geq2,\\
\frac{\left(  p-1\right)  ^{2}}{m-1} & \text{if }1<p<2.
\end{array}
\right.
\end{array}
\]

\end{lemma}

\proof

Fix $x\in\Omega\subset M$ with $du_{\epsilon}\neq0,$ we select a local
orthonormal frame field $\left\{  e_{1},e_{2},\ldots e_{m}\right\}  $ such
that at $x,$ $\nabla_{e_{i}}e_{j}=0,$ $\nabla u_{\epsilon}=\left\vert \nabla
u_{\epsilon}\right\vert e_{1},$ $u_{\epsilon, 1} = f\,,$ and $u_{\epsilon
,\alpha}=0$ for all $i,j=1,\ldots,m,$ $\alpha=2,\ldots\,, m\, $ where
$u_{\epsilon,\alpha}=\left\langle \nabla u_{\epsilon},e_{\alpha}\right\rangle
.$

Let $f=\left\vert \nabla u_{\epsilon}\right\vert ,$ $f_{\epsilon}%
=\sqrt{\left\vert \nabla u_{\epsilon}\right\vert ^{2}+\epsilon}$ and the
directional derivative $f_{\epsilon,i}=\left\langle \nabla f_{\epsilon}%
,e_{i}\right\rangle \,.$ Denote the directional derivative $\left\langle
\nabla u_{\epsilon,i},e_{j}\right\rangle $ by $u_{\epsilon,ij}\,.$ Then
(\ref{b0}) implies
\[%
\begin{array}
[c]{lllll}%
\Delta u_{\epsilon} & = & -\frac{p-2}{2f_{\epsilon}^{2}}\left\langle \nabla
f_{\epsilon}^{2},\nabla u_{\epsilon}\right\rangle  & = & -\frac{p-2}%
{2f_{\epsilon}^{2}}\sum_{i=1}^{m}\left(  f_{\epsilon}^{2}\right)
_{,i}u_{\epsilon,i}\\
&  &  & = & -\frac{p-2}{2f_{\epsilon}^{2}}\left(  f_{\epsilon}^{2}\right)
_{,1}u_{\epsilon,1}\\
&  &  & = & -\frac{p-2}{2f_{\epsilon}^{2}}\left(  f_{\epsilon}^{2}\right)
_{,1}f.
\end{array}
\]
Moreover, by using the following property%
\[%
\begin{array}
[c]{lllllll}%
\left(  f_{\epsilon}^{2}\right)  _{,j} & = & \left(  f^{2}\right)  \,_{,j} &
= & \sum_{i=1}^{m}\left(  u_{\epsilon,i}^{2}\right)  \,_{,j} & = & 2\sum
_{i=1}^{m}u_{\epsilon,i}u_{\epsilon,ij}\\
&  &  &  &  & = & 2u_{\epsilon,1}u_{\epsilon,1j}\\
&  &  &  &  & = & 2fu_{\epsilon,1j}.
\end{array}
\]
We have%
\begin{equation}%
\begin{array}
[c]{l}%
\Delta u_{\epsilon}=\frac{-\left(  p-2\right)  f^{2}}{f_{\epsilon}^{2}%
}u_{\epsilon,11},
\end{array}
\label{k3}%
\end{equation}
and
\begin{equation}%
\begin{array}
[c]{l}%
u_{\epsilon,1j}=f_{,j}.
\end{array}
\label{k4}%
\end{equation}
On the other hand,%

\begin{equation}%
\begin{array}
[c]{lll}%
\sum_{i,j=1}^{m}\left(  u_{\epsilon,ij}\right)  ^{2} & \geq & \left(
u_{\epsilon,11}\right)  ^{2}+2\sum_{\alpha=2}^{m}\left(  u_{\epsilon,1\alpha
}\right)  ^{2}+\sum_{\alpha=2}^{m}\left(  u_{\epsilon,\alpha\alpha}\right)
^{2}\\
& \geq & \left(  u_{\epsilon,11}\right)  ^{2}+2\sum_{\alpha=2}^{m}\left(
u_{\epsilon,1\alpha}\right)  ^{2}+\frac{\left(  \sum_{\alpha=2}^{m}%
u_{\epsilon,\alpha\alpha}\right)  ^{2}}{m-1}\\
& = & \left(  u_{\epsilon,11}\right)  ^{2}+2\sum_{\alpha=2}^{m}\left(
u_{\epsilon,1\alpha}\right)  ^{2}+\frac{\left(  \Delta u_{\epsilon
}-u_{\epsilon,11}\right)  ^{2}}{m-1}.
\end{array}
\label{k1}%
\end{equation}
Therefore, by using (\ref{k3}) and (\ref{k4}), the inequality (\ref{k1})
implies%
\[%
\begin{array}
[c]{lll}%
\sum_{i,j=1}^{m}\left(  u_{\epsilon,ij}\right)  ^{2} & \geq & \left(
u_{\epsilon,11}\right)  ^{2}+2\sum_{\alpha=2}^{m}\left(  u_{\epsilon,1\alpha
}\right)  ^{2}+\frac{\left(  \left(  \frac{\left(  p-2\right)  f^{2}%
}{f_{\epsilon}^{2}}+1\right)  u_{\epsilon,11}\right)  ^{2}}{m-1}\\
& = & \left(  1+\frac{\left(  \left(  p-1\right)  f^{2}+\epsilon\right)  ^{2}%
}{\left(  m-1\right)  f_{\epsilon}^{4}}\right)  \left(  u_{\epsilon
,11}\right)  ^{2}+2\sum_{\alpha=2}^{m}\left(  u_{\epsilon,1\alpha}\right)
^{2}\\
& \geq & \left(  1+\kappa\right)  \left\vert \nabla f\right\vert ^{2}%
\end{array}
\]
which can be written as
\[%
\begin{array}
[c]{lll}%
\left\vert du_{\epsilon}\right\vert ^{2}\left\vert \nabla du_{\epsilon
}\right\vert ^{2} & \geq & \frac{1+\kappa_{1}}{4}\left\vert \nabla\left\vert
du_{\epsilon}\right\vert ^{2}\right\vert ^{2}%
\end{array}
\]
for all $x\in\Omega.$ This completes the proof. \endproof

\bigskip

\begin{lemma}
\label{ka2}Let $u_{\epsilon}$ be a solution of (\ref{1.2}) on $\Omega\subset
M^{m},$ $p>1.$ Then the Hessian of $u_{\epsilon}$ satisfies
\begin{equation}%
\begin{array}
[c]{lll}%
\left(  \left\vert du_{\epsilon}\right\vert ^{2}+\epsilon\right)  \left\vert
\nabla du_{\epsilon}\right\vert ^{2} & \geq & \frac{1+\kappa_{2}}{4}\left\vert
\nabla\left\vert du_{\epsilon}\right\vert ^{2}\right\vert ^{2}%
\end{array}
\label{k''}%
\end{equation}
at $x\in\Omega,$ where $\kappa_{2}=\min\left\{  \frac{\left(  p-1\right)
^{2}}{m},1\right\}  .$
\end{lemma}

\proof Since $v_{\epsilon}\in C^{\infty}\left(  N\right)  $ is the strongly
$p$-harmonic function on $\left(  \Omega_{N},g_{N}\right)  ,$ then Kato's
inequality for strongly $p$-harmonic function on $\left(  \Omega_{N}%
,g_{N}\right)  $ (see Lemma \ref{KS}) implies%
\begin{equation}%
\begin{array}
[c]{l}%
\left\vert \left(  \nabla d\right)  ^{N}v_{\epsilon}\right\vert ^{2}%
\geq\left(  1+\kappa_{2}\right)  \left\vert \nabla^{N}\left\vert
d^{N}v_{\epsilon}\right\vert \right\vert ^{2}%
\end{array}
\label{k5}%
\end{equation}
where $\left(  \nabla d\right)  ^{N}$ is the Hessian on $\left(  \Omega
_{N},g_{N}\right)  ,$ and $\kappa_{2}=\min\left\{  \frac{\left(  p-1\right)
^{2}}{m},1\right\}  .$ Moreover, (\ref{k5}) can be rewritten as%
\[%
\begin{array}
[c]{lll}%
\left(  \left\vert du_{\epsilon}\right\vert ^{2}+\epsilon\right)  \left\vert
\nabla du_{\epsilon}\right\vert ^{2} & \geq & \left(  1+\kappa_{2}\right)
\left(  \left\vert du_{\epsilon}\right\vert ^{2}+\epsilon\right)  \left\vert
\nabla\sqrt{\left\vert du_{\epsilon}\right\vert ^{2}+\epsilon}\right\vert
^{2}\\
& = & \frac{1+\kappa_{2}}{4}\left\vert \nabla\left\vert du_{\epsilon
}\right\vert ^{2}\right\vert ^{2}.
\end{array}
\]
for all $x\in\Omega.$ \ \endproof

\section{The Proof of Theorem \ref{T1}}

Now we use Lemmas \ref{Bo} - \ref{ka2} and weighted Poincar\'{e} inequality
(\ref{WP}) to obtain the following inequality (\ref{v3}):

\begin{lemma}
\label{m1}Let $M$ be a manifold satisfying the hypothesis of Theorem \ref{T1}.
Let $u_{\epsilon}$ be a solution of (\ref{1.2}) on $B\left(  2R\right)
\subset M.$ Then we have
\begin{equation}%
\begin{array}
[c]{lll}%
C\int_{B\left(  R\right)  }\rho\left\vert \nabla u_{\epsilon}\right\vert
^{p}dv & \leq & \frac{100\cdot B}{R^{2}}\int_{B\left(  2R\right)  \backslash
B\left(  R\right)  }\left(  \left\vert \nabla u_{\epsilon}\right\vert
^{2}+\epsilon\right)  ^{\frac{p}{2}}dv,
\end{array}
\label{v3}%
\end{equation}
where $C\left(  p,m,\kappa,\tau,\epsilon_{1},\epsilon_{2}\right)  >0$ and
$B\left(  p,m,\kappa,\epsilon_{1},\epsilon_{2}\right)  >0$ are positive
constants for sufficiently small constants $\epsilon_{1}, \epsilon_{2} > 0\,
.$
\end{lemma}

\proof

Let $\Omega=B\left(  2R\right)  $ be a geodesic ball of radius $2R$ centered
at a fixed point.

Let $f=\left\vert \nabla u_{\epsilon}\right\vert $ and $f_{\epsilon}%
=\sqrt{f^{2}+\epsilon}.$ In view of Lemma \ref{ka1} and Lemma \ref{ka2},%

\[%
\begin{array}
[c]{ccc}%
f_{\epsilon}^{2}\left\vert \nabla du_{\epsilon}\right\vert ^{2} & \geq &
\frac{1+\kappa}{4}\left\vert \nabla f^{2}\right\vert ^{2}%
\end{array}
\]
holds for all on $M,$ where $\kappa=\max\left\{  \kappa_{1},\kappa
_{2}\right\}  .$ Then by Lemma \ref{Bo}, we rewrite Bochner's formula as%
\begin{equation}%
\begin{array}
[c]{lll}%
\frac{1}{2}\mathcal{L}_{\epsilon}\left(  f_{\epsilon}^{2}\right)  & \geq &
\left(  p-1+\kappa\right)  f_{\epsilon}^{p-2}\left\vert \nabla f_{\epsilon
}\right\vert ^{2}+f_{\epsilon}^{p-2}Ric\left(  \nabla u_{\epsilon},\nabla
u_{\epsilon}\right)  ,
\end{array}
\label{v0}%
\end{equation}
here we use $\nabla f_{\epsilon}^{2}=\nabla f^{2}\,.$

We multiply both sides of (\ref{v0}) by $\eta^{2}$ and integrate over $M,$
\begin{equation}%
\begin{array}
[c]{lll}%
\frac{1}{2}\int_{M}\eta^{2}\mathcal{L}_{\epsilon}\left(  f_{\epsilon}%
^{2}\right)  dv & \geq & \left(  p-1+\kappa\right)  \int_{M}\eta
^{2}f_{\epsilon}^{p-2}\left\vert \nabla f_{\epsilon}\right\vert ^{2}dv\\
&  & +\int_{M}\eta^{2}f_{\epsilon}^{p-2}Ric\left(  \nabla u_{\epsilon},\nabla
u_{\epsilon}\right)  dv
\end{array}
\label{v1}%
\end{equation}
where $\eta\in C_{0}^{\infty}(M)$ is a cut-off function with $0\leq\eta\left(
x\right)  \leq1$ on $M$ satisfying%
\[
\left\{
\begin{array}
[c]{ll}%
\eta\left(  x\right)  =1 & \text{if }x\in\overline{B\left(  R\right)  },\\
\left\vert \nabla\eta\left(  x\right)  \right\vert \leq\frac{10}%
{R}\text{\ \ \ \ \ } & \text{if }x\in B\left(  2R\right)  \backslash
\overline{B\left(  R\right)  },\\
\eta\left(  x\right)  =0 & \text{if }x\in M\backslash B\left(  2R\right)  .
\end{array}
\right.
\]

On the other hand, applying integration by parts and Cauchy-Schwarz inequality
one has%
\[%
\begin{array}
[c]{lll}%
\frac{1}{2}\int_{M}\eta^{2}\mathcal{L}_{\epsilon}\left(  f_{\epsilon}%
^{2}\right)  dv & = & \frac{-1}{2}\int_{M}\left\langle \nabla\eta
^{2},f_{\epsilon}^{p-2}\nabla f_{\epsilon}^{2}+\left(  p-2\right)
f_{\epsilon}^{p-4}\left\langle \nabla u_{\epsilon},\nabla f_{\epsilon}%
^{2}\right\rangle \nabla u_{\epsilon}\right\rangle dv\\
& \leq & 2\int_{M}\eta\left\vert \nabla\eta\right\vert \left(  f_{\epsilon
}^{p-1}\left\vert \nabla f_{\epsilon}\right\vert +\left\vert p-2\right\vert
f_{\epsilon}^{p-3}f^{2}\left\vert \nabla f_{\epsilon}\right\vert \right)  dv\\
& \leq & 2\left(  1+\left\vert p-2\right\vert \right)  \int_{M}\eta\left\vert
\nabla\eta\right\vert f_{\epsilon}^{p-1}\left\vert \nabla f_{\epsilon
}\right\vert dv\\
& \leq & \epsilon_{1}\int_{M}\eta^{2}f_{\epsilon}^{p-2}\left\vert \nabla
f_{\epsilon}\right\vert ^{2}dv+\frac{\left(  1+\left\vert p-2\right\vert
\right)  ^{2}}{\epsilon_{1}}\int_{M}\left\vert \nabla\eta\right\vert
^{2}f_{\epsilon}^{p}dv,
\end{array}
\]
where $\epsilon_{1}$ is a positive constant satisfying%
\[%
\begin{array}
[c]{l}%
p-1-\epsilon_{1}>0.
\end{array}
\]

Then (\ref{v1}) implies%
\begin{equation}%
\begin{array}
[c]{lll}%
\frac{\left(  1+\left\vert p-2\right\vert \right)  ^{2}}{\epsilon_{1}}\int
_{M}\left\vert \nabla\eta\right\vert ^{2}f_{\epsilon}^{p}dv & \geq & \int
_{M}\left(  p-1+\kappa-\epsilon_{1}\right)  \eta^{2}f_{\epsilon}%
^{p-2}\left\vert \nabla f_{\epsilon}\right\vert ^{2}dv\\
&  & +\int_{M}\eta^{2}f_{\epsilon}^{p-2}Ric\left(  \nabla u_{\epsilon},\nabla
u_{\epsilon}\right)  dv.
\end{array}
\label{v1-2}%
\end{equation}

Besides, we may rewrite the first term in the right hand side of (\ref{v1-2})
by
\[%
\begin{array}
[c]{lll}
&  & \left(  p-1+\kappa-\epsilon_{1}\right)  \int_{M}\eta^{2}f_{\epsilon
}^{p-2}\left\vert \nabla f_{\epsilon}\right\vert ^{2}dv\\
& = & \frac{4\left(  p-1+\kappa-\epsilon_{1}\right)  }{p^{2}}\int_{M}\eta
^{2}\left\vert \nabla f_{\epsilon}^{\frac{p}{2}}\right\vert ^{2}dv\\
& = & \frac{4\left(  p-1+\kappa-\epsilon_{1}\right)  }{p^{2}}\int
_{M}\left\vert \nabla\left(  \eta f_{\epsilon}^{\frac{p}{2}}\right)  -\left(
\nabla\eta\right)  f_{\epsilon}^{\frac{p}{2}}\right\vert ^{2}dv\\
& = & \frac{4\left(  p-1+\kappa-\epsilon_{1}\right)  }{p^{2}}\int_{M}\left\{
\left\vert \nabla\left(  \eta f_{\epsilon}^{\frac{p}{2}}\right)  \right\vert
^{2}-2\left\langle \nabla\left(  \eta f_{\epsilon}^{\frac{p}{2}}\right)
,f_{\epsilon}^{\frac{p}{2}}\nabla\eta\right\rangle +\left\vert \nabla
\eta\right\vert ^{2}f_{\epsilon}^{p}\right\}  dv\\
& \geq & \frac{4\left(  1-\epsilon_{2}\right)  \left(  p-1+\kappa-\epsilon
_{1}\right)  }{p^{2}}\int_{M}\left\vert \nabla\left(  \eta f_{\epsilon}%
^{\frac{p}{2}}\right)  \right\vert ^{2}+\frac{4\left(  1-\frac{1}{\epsilon
_{2}}\right)  \left(  p-1+\kappa-\epsilon_{1}\right)  }{p^{2}}\int
_{M}\left\vert \nabla\eta\right\vert ^{2}f_{\epsilon}^{p}dv.
\end{array}
\]
where $\epsilon_{2}$ is a positive constant satisfying $\epsilon_{2}<1.$ Thus,
we have%
\begin{equation}%
\begin{array}
[c]{lll}
&  & \frac{4\left(  1-\epsilon_{2}\right)  \left(  p-1+\kappa-\epsilon
_{1}\right)  }{p^{2}}\int_{M}\left\vert \nabla\left(  \eta f_{\epsilon}%
^{\frac{p}{2}}\right)  \right\vert ^{2}dv+\int_{M}\eta^{2}f_{\epsilon}%
^{p-2}Ric\left(  \nabla u_{\epsilon},\nabla u_{\epsilon}\right)  \,dv\\
& \leq & \left(  \frac{\left(  1+\left\vert p-2\right\vert \right)  ^{2}%
}{\epsilon_{1}}+\frac{4\left(  \frac{1}{\epsilon_{2}}-1\right)  \left(
p-1+\kappa-\epsilon_{1}\right)  }{p^{2}}\right)  \int_{M}\left\vert \nabla
\eta\right\vert ^{2}f_{\epsilon}^{p}dv.
\end{array}
\label{v2}%
\end{equation}

According to the weighted Poincar\'{e} inequality (\ref{WP})
\[%
\begin{array}
[c]{lll}%
\int_{M}\rho\Psi^{2}dv & \leq & \int_{M}\left\vert \nabla\Psi\right\vert
^{2}dv
\end{array}
\]
with $\Psi=\eta f_{\epsilon}^{\frac{p}{2}},$ then (\ref{v2}) implies%
\begin{equation}%
\begin{array}
[c]{lll}%
\int_{B\left(  R\right)  }Af_{\epsilon}^{p-2}dv & \leq & \frac{100\cdot
B}{R^{2}}\int_{B\left(  2R\right)  \backslash B\left(  R\right)  }f_{\epsilon
}^{p}dv,
\end{array}
\label{v3-1}%
\end{equation}
for all fixed $R>0,$ where
\[%
\begin{array}
[c]{lll}%
A & = & \frac{4\left(  1-\epsilon_{2}\right)  \left(  p-1+\kappa-\epsilon
_{1}\right)  }{p^{2}}\rho f_{\epsilon}^{2} + Ric\left(  \nabla u_{\epsilon
},\nabla u_{\epsilon}\right)  \,,
\end{array}
\]
and%
\[%
\begin{array}
[c]{lll}%
B & = & \left(  \frac{\left(  1+\left\vert p-2\right\vert \right)  ^{2}%
}{\epsilon_{1}}+\frac{4\left(  \frac{1}{\epsilon_{2}}-1\right)  \left(
p-1+\kappa-\epsilon_{1}\right)  }{p^{2}}\right)  .
\end{array}
\]

Since the curvature condition (\ref{1.4}) means that there exists a constant
$0<\tau<\frac{4\left(  p-1+\kappa\right)  }{p^{2}}$ such that
\[%
\begin{array}
[c]{lll}%
Ric_{M} & \geq & -\tau\rho,
\end{array}
\]
Then%
\[%
\begin{array}
[c]{lll}%
A & \geq & C\left(  p,m,\kappa,\tau,\epsilon_{1},\epsilon_{2}\right)  \rho
f^{2}%
\end{array}
\]
with $C>0$ whenever we select $\epsilon_{1}$ and $\epsilon_{2}$ small enough.

Hence, (\ref{v3-1}) gives
\[%
\begin{array}
[c]{lll}%
C\int_{B\left(  R\right)  }\rho\left\vert \nabla u_{\epsilon}\right\vert
^{p}dv & \leq & \frac{100\cdot B}{R^{2}}\int_{B\left(  2R\right)  \backslash
B\left(  R\right)  }\left(  \left\vert \nabla u_{\epsilon}\right\vert
^{2}+\epsilon\right)  ^{\frac p2}dv,
\end{array}
\]
where $C\left(  p,m,\kappa,\tau,\epsilon_{1},\epsilon_{2}\right)  >0,$ and
$B\left(  p,m,\kappa,\epsilon_{1},\epsilon_{2}\right)  >0.$

\endproof

\bigskip

\textbf{Proof of Theorem \ref{T1}.}

\proof

Given $B\left(  R_{0}\right)  \subset M,$ for every $a>0,$ we let $\Omega
_{a}=\left\{  x\in B\left(  R_{0}\right)  :\rho\left(  x\right)  >1/a\right\}
.$ It is clear the measure of $\Omega_{a}$ tends to zero as $a\rightarrow
0^{+}.$ If we are able to show $\int_{B\left(  R_{0}\right)  \backslash
\Omega_{a}}\rho\left\vert \nabla u\right\vert ^{p}dv<\delta$ for any
$\delta>0\,,$ then it implies $\nabla u=0$ on $B\left(  R_{0}\right)  $ almost
everywhere. This also infers $\nabla u=0$ on $B\left(  R_{0}\right)  $ by the
fact $u\in C_{loc}^{1,\alpha}\left(  M\right)  . $ Moreover, since $B\left(
R_{0}\right)  $ is arbitrary, $u$ must be constant on $M.$

Moreover, if we assume $M$ has at least two $p$-hyperbolic ends. By
Proposition \ref{2 E}, one may construct a nontrivial bounded $p$-harmonic
function with finite $p$-energy on $M$, this gives a contradiction to our
conclusion, hence $M$ has only one $p$-hyperbolic end.

Now we prove the claim. By using the finite $p$-energy of $u,$ we may select
$0<<R<\infty$ large enough such that $B\left(  R_{0}\right)  \subset B\left(
R\right)  $ and
\[%
\begin{array}
[c]{lll}%
\frac{100B}{R^{2}C}\int_{B\left(  2R\right)  \backslash B\left(  R\right)
}\left\vert \nabla u\right\vert ^{p}dv & < & \delta
\end{array}
\]
where $B$ and $C$ are defined as (\ref{v3}).

Now we construct $u_{\epsilon}\in C^{\infty}\left(  B\left(  2R\right)
\right)  $ such that $u_{\epsilon}=u$ on $\partial B\left(  2R\right)  $ and
$u_{\epsilon}$ satisfies (\ref{1.2}). Then (\ref{v3}) implies
\[%
\begin{array}
[c]{lll}%
C\int_{B\left(  R_{0}\right)  \backslash\Omega_{a}}\rho\left\vert \nabla
u\right\vert ^{p}dv & \leq & \frac{100\cdot B}{R^{2}}\int_{B\left(  2R\right)
\backslash B\left(  R\right)  }\left\vert \nabla u\right\vert ^{p}dv,
\end{array}
\]
as $\epsilon\rightarrow0,$ we may therefore conclude that
\[%
\begin{array}
[c]{lll}%
\int_{B\left(  R_{0}\right)  \backslash\Omega_{a}}\rho\left\vert \nabla
u\right\vert ^{p}dv & < & \delta.
\end{array}
\]

\endproof

\bigskip

\bigskip

If $M$ has positive spectrum $\lambda>0,$ then $M$ has $p$-Poincar\'{e}
inequality%
\[%
\begin{array}
[c]{l}%
\lambda_{p}\int_{M}\Psi^{p}\leq\int_{M}\left\vert \nabla\Psi\right\vert
^{p},\text{ }\lambda_{p}>0
\end{array}
\]
for all $\Psi\in W_{0}^{1,p}\left(  M\right)  $ and $p\geq2$ (cf. \cite{HKM}
Theorem 1.8). Since $p$-Poincar\'{e} inequality and Caccioppoli type estimate
imply decay estimate (see Lemma \ref{decay estimate on E} which is similar to
the work of \cite{LW1} Lemma 1.1 and Lemma 1.2), then $p$-Poincar\'{e}
inequality infers that $M$ must be a $p$-hyperbolic manifold (see Theorem
\ref{hy}). So we have the following:

\begin{corollary}
Let $M^{m},m\geq2,$ be a complete noncompact Riemannian manifold with positive
spectrum $\lambda>0$ and
\[%
\begin{array}
[c]{lll}%
Ric_{M} & \geq & -\tau\lambda
\end{array}
\]
where $p\geq2,$ and constant $\tau$ is the same as in Theorem \ref{T1}. Then
every weakly $p$-harmonic function $u$ with finite $p$-energy is constant.
Moreover, $M$ has only one $p$-hyperbolic end.
\end{corollary}

\begin{remark}
Similarly, if $M$ has $p$-Poincar\'{e} inequality, $1<p<2,$ then $M$ has
positive spectrum $\lambda>0.$ Hence, if $M$ is a complete noncompact
Riemannian manifold with $p$-Poincar\'{e} inequality, $1<p<2,$ and
$Ric_{M}\geq-\tau\lambda$ where $\tau<\frac{4\left(  p-1\right)  \left(
p+m-2\right)  }{p^{2}\left(  m-1\right)  }.$ Then $M$ has only one
$p$-hyperbolic end.
\end{remark}

\bigskip

\section{Strongly $p$-harmonic functions with applications}

\subsection{Bochner's formula}

Let $u$ be a $C^{3}$ strongly $p$-harmonic function for $p>1$ on $M\, .$ Then
$\left\vert \nabla u\right\vert ^{p-2}\nabla u$ must be $C^{1}$ on $M,$ and
hence $u$ is a solution of (\ref{pl}) as follows:

\begin{lemma}
\label{PL}If $u \in C^{3}(M)$ is a strongly $p$-harmonic function for $p>1$,
then $u$ is a solution of
\begin{equation}%
\begin{array}
[c]{lll}%
f^{2}\Delta u+\frac{p-2}{2}\left\langle \nabla f^{2},\nabla u\right\rangle  &
= & 0,
\end{array}
\label{pl}%
\end{equation}
on $M\, ,$ where $f=\left\vert \nabla u\right\vert .$
\end{lemma}

\proof

First, we multiply both side of (\ref{1.0}) by $f^{4},$ because of $f^{4}\in
C^{2}\left(  M\right)  ,$ then%
\[%
\begin{array}
[c]{l}%
f^{4}\text{\textrm{div}}(f^{p-2}\nabla u)=0
\end{array}
\]
implies%
\[%
\begin{array}
[c]{lll}%
0 & = & \text{\textrm{div}}(f^{p+2}\nabla u)-2f^{p}\left\langle \nabla
f^{2},\nabla u\right\rangle \\
& = & f^{p+2}\Delta u+\left\langle \nabla f^{p+2},\nabla u\right\rangle
-2f^{p}\left\langle \nabla f^{2},\nabla u\right\rangle .
\end{array}
\]
Since $p>1$ and%
\[%
\begin{array}
[c]{l}%
\nabla f^{p+2}=\nabla\left(  \left(  f^{2}\right)  ^{\frac{p+2}{2}}\right)
=\frac{p+2}{2}f^{p}\nabla f^{2},
\end{array}
\]
so we have%
\[%
\begin{array}
[c]{l}%
f^{p+2}\Delta u+\frac{p-2}{2}f^{p}\left\langle \nabla f^{2},\nabla
u\right\rangle =0.
\end{array}
\]
which implies%
\[%
\begin{array}
[c]{l}%
f^{2}\Delta u+\frac{p-2}{2}\left\langle \nabla f^{2},\nabla u\right\rangle =0
\end{array}
\]
on all of $M.$ \ 

\endproof

\begin{remark}
(1). If $u$ is a solution of (\ref{pl}), $u$ may be not a strongly
$p$-harmonic function. For example, any constant function is a solution of
(\ref{pl}), but it is not a strongly $p$-harmonic function for $1<p<2.$
\newline(2). For $p\geq4, $ $u\in C^{2}\left(  M\right)  $ is a solution of
(\ref{pl}) if and only if $u$ is a strongly $p$-harmonic function.
\end{remark}

Now we define an operator $\mathcal{L}_{s,\varepsilon}$ by%
\[%
\begin{array}
[c]{l}%
\mathcal{L}_{s,\varepsilon}\left(  \Psi\right)  =\text{\textrm{div}}\left(
f_{\varepsilon}^{s}A_{\varepsilon}\left(  \nabla\Psi\right)  \right)  ,
\end{array}
\]
for $\Psi\in C^{2}\left(  M\right)  ,$ where $s\in\mathbb{R},$ $p>1,$
$\varepsilon>0,$ $f_{\varepsilon}=\sqrt{f^{2}+\varepsilon}$ and%
\[%
\begin{array}
[c]{l}%
A_{\varepsilon}:=\mathrm{id}+\left(  p-2\right)  \frac{\nabla u\otimes\nabla
u}{f_{\varepsilon}^{2}}.
\end{array}
\]
Note that $\mathcal{L}_{s,\varepsilon}$ is a linearized operator of the
nonlinear equation (\ref{1.0}), and $\mathcal{L}_{s,\varepsilon}\left(
f_{\varepsilon}^{2}\right)  \left(  x\right)  $ is well define for all $x\in
M$ since $f_{\varepsilon}>0$ and $f_{\varepsilon}^{2}\in C^{2}\left(
M\right)  .$ \bigskip

Next we use the operator $\mathcal{L}_{s,\varepsilon}$ to derive the Bochner's
formula for the solution of (\ref{pl}).

\begin{lemma}
[Bochner's formula]\label{BS}If $u\in C^{3}\left(  M\right)  $ is a strongly
$p$-harmonic function. Let $f=\left\vert \nabla u\right\vert $ and
$f_{\varepsilon}=\sqrt{f^{2}+\varepsilon},$ then for all $p>1$ and
$s\in\mathbb{R},$ the formula
\[%
\begin{array}
[c]{lll}%
\frac{1}{2}\mathcal{L}_{s,\varepsilon}\left(  f_{\varepsilon}^{2}\right)  &
= & \frac{s}{4}f_{\varepsilon}^{s-2}\left\vert \nabla f_{\varepsilon}%
^{2}\right\vert ^{2}+f_{\varepsilon}^{s}\sum_{i,j=1}^{m}\left(  u_{ij}%
^{2}+R_{ij}u_{i}u_{j}\right) \\
&  & +\frac{\left(  p-2\right)  \left(  s-p+2\right)  }{4}f_{\varepsilon
}^{s-4}\left\langle \nabla u,\nabla f_{\varepsilon}^{2}\right\rangle ^{2}\\
&  & +\varepsilon\left(  f_{\varepsilon}^{s-2}\left\langle \nabla
u,\nabla\Delta u\right\rangle +\frac{p-4}{2}f_{\varepsilon}^{s-4}\left\langle
\nabla u,\nabla f_{\varepsilon}^{2}\right\rangle \Delta u\right)
\end{array}
\]
holds on all of $M\, ,$ where $R_{ij} = \sum_{k=1}^{m} \langle R (e_{i},
e_{k}) e_{k}, e_{j} \rangle$ is the Ricci curvature tensor of $M\, .$ In
particular, if $p=2,$ then%
\[%
\begin{array}
[c]{l}%
\frac{1}{2}\mathcal{L}_{s,\varepsilon}\left(  f_{\varepsilon}^{2}\right)
=\frac{s}{4}f_{\varepsilon}^{s-2}\left\vert \nabla f_{\varepsilon}%
^{2}\right\vert ^{2}+f_{\varepsilon}^{s}\sum_{i,j=1}^{m}\left(  u_{ij}%
^{2}+R_{ij}u_{i}u_{j}\right)
\end{array}
\]
holds on all of $M$ and for all $s\in\mathbb{R}.$
\end{lemma}

\proof

By Lemma \ref{PL}, $u$ must be a solution of (\ref{pl}). Taking the gradient
of both sides of (\ref{pl}), and then taking the inner product with $\nabla
u\,,$ we have
\begin{equation}%
\begin{array}
[c]{lll}%
0 & = & \frac{p-2}{2}\left\langle \nabla\left\langle \nabla f^{2},\nabla
u\right\rangle ,\nabla u\right\rangle +\left\langle \nabla f^{2},\nabla
u\right\rangle \Delta u\\
&  & +f^{2}\left\langle \nabla\left(  \Delta u\right)  ,\nabla u\right\rangle
.
\end{array}
\label{w1}%
\end{equation}

Now we rewrite $\mathcal{L}_{s,\varepsilon}\left(  f_{\varepsilon}^{2}\right)
$ as the following formula,
\begin{equation}%
\begin{array}
[c]{lll}%
\frac{1}{2}\mathcal{L}_{s,\varepsilon}\left(  f_{\varepsilon}^{2}\right)  &
= & \frac{1}{2}\text{\textrm{div}}\left(  f_{\varepsilon}^{s}\nabla
f_{\varepsilon}^{2}+\left(  p-2\right)  f_{\varepsilon}^{s-2}\left\langle
\nabla u,\nabla f_{\varepsilon}^{2}\right\rangle \nabla u\right) \\
& = & \frac{s}{4}f_{\varepsilon}^{s-2}\left\vert \nabla f_{\varepsilon}%
^{2}\right\vert ^{2}+\frac{1}{2}f_{\varepsilon}^{s}\Delta f_{\varepsilon}%
^{2}\\
&  & +\frac{\left(  p-2\right)  \left(  s-2\right)  }{4}f_{\varepsilon}%
^{s-4}\left\langle \nabla u,\nabla f_{\varepsilon}^{2}\right\rangle ^{2}\\
&  & +\frac{p-2}{2}f_{\varepsilon}^{s-2}\left\langle \nabla\left\langle \nabla
u,\nabla f_{\varepsilon}^{2}\right\rangle ,\nabla u\right\rangle \\
&  & +\frac{p-2}{2}f_{\varepsilon}^{s-2}\left\langle \nabla u,\nabla
f_{\varepsilon}^{2}\right\rangle \Delta u.
\end{array}
\label{w2}%
\end{equation}
Combining (\ref{w1}), one has%
\begin{equation}%
\begin{array}
[c]{lll}%
\frac{1}{2}\mathcal{L}_{s,\varepsilon}\left(  f_{\varepsilon}^{2}\right)  &
= & \frac{s}{4}f_{\varepsilon}^{s-2}\left\vert \nabla f_{\varepsilon}%
^{2}\right\vert ^{2}+\frac{1}{2}f_{\varepsilon}^{s}\Delta f_{\varepsilon}%
^{2}\\
&  & +\frac{\left(  p-2\right)  \left(  s-2\right)  }{4}f_{\varepsilon}%
^{s-4}\left\langle \nabla u,\nabla f_{\varepsilon}^{2}\right\rangle ^{2}\\
&  & -f_{\varepsilon}^{s-2}f^{2}\left\langle \nabla\left(  \Delta u\right)
,\nabla u\right\rangle \\
&  & +\frac{p-4}{2}f_{\varepsilon}^{s-2}\left\langle \nabla u,\nabla
f_{\varepsilon}^{2}\right\rangle \Delta u,
\end{array}
\label{w3}%
\end{equation}
here we use the fact $\nabla f_{\varepsilon}^{2}=\nabla f^{2}.$

According to (\ref{pl}), the last term of right hand side can be rewritten as%
\[%
\begin{array}
[c]{lll}%
f_{\varepsilon}^{s-2}\left\langle \nabla u,\nabla f_{\varepsilon}%
^{2}\right\rangle \Delta u & = & f_{\varepsilon}^{s-4}\left(  f^{2}%
+\varepsilon\right)  \left\langle \nabla u,\nabla f_{\varepsilon}%
^{2}\right\rangle \Delta u\\
& = & f_{\varepsilon}^{s-4}f^{2}\left\langle \nabla u,\nabla f_{\varepsilon
}^{2}\right\rangle \Delta u+\varepsilon f_{\varepsilon}^{s-4}\left\langle
\nabla u,\nabla f_{\varepsilon}^{2}\right\rangle \Delta u\\
& = & -\frac{p-2}{2}f_{\varepsilon}^{s-4}\left\langle \nabla u,\nabla
f_{\varepsilon}^{2}\right\rangle ^{2}+\varepsilon f_{\varepsilon}%
^{s-4}\left\langle \nabla u,\nabla f_{\varepsilon}^{2}\right\rangle \Delta u.
\end{array}
\]
Using Bochner's formula%
\[%
\begin{array}
[c]{l}%
\frac{1}{2}\Delta f^{2}=\sum_{i,j=1}^{m}u_{ij}^{2}+\left\langle \nabla
u,\nabla\Delta u\right\rangle +\sum_{i,j=1}^{m}R_{ij}u_{i}u_{j}%
\end{array}
\]
and the equality $\Delta f^{2}=\Delta f_{\varepsilon}^{2},$ then (\ref{w3})
gives the desired%
\[%
\begin{array}
[c]{lll}%
\frac{1}{2}\mathcal{L}_{s,\varepsilon}\left(  f_{\varepsilon}^{2}\right)  &
= & \frac{s}{4}f_{\varepsilon}^{s-2}\left\vert \nabla f_{\varepsilon}%
^{2}\right\vert ^{2}+f_{\varepsilon}^{s}\sum_{i,j=1}^{m}\left(  u_{ij}%
^{2}+R_{ij}u_{i}u_{j}\right) \\
&  & +\frac{\left(  p-2\right)  \left(  s-p+2\right)  }{4}f_{\varepsilon
}^{s-4}\left\langle \nabla u,\nabla f_{\varepsilon}^{2}\right\rangle ^{2}\\
&  & +\varepsilon\left(  f_{\varepsilon}^{s-2}\left\langle \nabla
u,\nabla\Delta u\right\rangle +\frac{p-4}{2}f_{\varepsilon}^{s-4}\left\langle
\nabla u,\nabla f_{\varepsilon}^{2}\right\rangle \Delta u\right)  .
\end{array}
\]
\endproof

\begin{lemma}
[Refined Kato's inequality]\label{KS}Let $u\in C^{2}\left(  M\right)  $ be
$p$-harmonic function on a complete manifold $M^{m},$ $p>1$ and $\kappa
=\min\left\{  \frac{\left(  p-1\right)  ^{2}}{m-1},1\right\}  .$ Then at any
$x\in M$ with $du\left(  x\right)  \neq0,$
\begin{equation}%
\begin{array}
[c]{l}%
\left\vert \nabla\left(  du\right)  \right\vert ^{2}\geq\left(  1+\kappa
\right)  \left\vert \nabla\left\vert du\right\vert \right\vert ^{2},
\end{array}
\label{ks}%
\end{equation}
and "$=$" holds if and only if
\[
\left\{
\begin{array}
[c]{ll}%
u_{\alpha\beta}=0\text{ and }u_{11}=-\frac{m-1}{p-1}u_{\alpha\alpha}, &
\text{for }\left(  p-1\right)  ^{2}=m-1,\\
u_{\alpha\beta}=0,\text{ }u_{1\alpha}=0\text{ and }u_{11}=-\frac{m-1}%
{p-1}u_{\alpha\alpha}, & \text{for }\left(  p-1\right)  ^{2}<m-1,\\
u_{\alpha\beta}=0\text{ and }u_{ii}=0,\text{ } & \text{for }\left(
p-1\right)  ^{2}>m-1,
\end{array}
\right.
\]
for all $\alpha,\beta=2,\ldots,m,$ $\alpha\neq\beta$ and $i=1,\ldots,m.$
\end{lemma}

\proof

Fix a point $x\in M.$ If $du\neq0$ at $x,$ we are able to select a local
orthonormal frame field $\left\{  e_{1},e_{2},\ldots e_{m}\right\}  $ such
that, at $x,$ $\nabla_{e_{i}}e_{j}=0,$ $\nabla u=\left\vert \nabla
u\right\vert e_{1},$ and $u_{\alpha}=0$ for all $i,j=1,\ldots,m,$
$\alpha=2,\ldots,m.$ Here we use the convenient notation $u_{i}=\left\langle
\nabla u,e_{i}\right\rangle .$

Observing that
\begin{equation}%
\begin{array}
[c]{lll}%
\sum_{i,j=1}^{m}\left(  u_{ij}\right)  ^{2} & \geq & \left(  u_{11}\right)
^{2}+2\sum_{\alpha=2}^{m}\left(  u_{1\alpha}\right)  ^{2}+\sum_{\alpha=2}%
^{m}\left(  u_{\alpha\alpha}\right)  ^{2}\\
& \geq & \left(  u_{11}\right)  ^{2}+2\sum_{\alpha=2}^{m}\left(  u_{1\alpha
}\right)  ^{2}+\frac{\left(  \sum_{\alpha=2}^{m}u_{\alpha\alpha}\right)  ^{2}%
}{m-1}\\
& = & \left(  u_{11}\right)  ^{2}+2\sum_{\alpha=2}^{m}\left(  u_{1\alpha
}\right)  ^{2}+\frac{\left(  \Delta u-u_{11}\right)  ^{2}}{m-1}.
\end{array}
\label{ks1}%
\end{equation}
However, letting $f=\left\vert \nabla u\right\vert \, $ and using $f = u_{1}\,
,$ $f_{1} = \langle\nabla f, e_{1} \rangle\, , $%
\[%
\begin{array}
[c]{lll}%
0 & = & div\left(  f^{p-2}\nabla u\right) \\
& = & f^{p-2}\Delta u+\left(  p-2\right)  f^{p-3}\left\langle \nabla f,\nabla
u\right\rangle \\
& = & u_{1}^{p-2}\Delta u+\left(  p-2\right)  u_{1}^{p-2}f_{1},
\end{array}
\]
and%
\begin{equation}%
\begin{array}
[c]{lll}%
f_{j}=\frac{\left(  f^{2}\right)  _{,j}}{2f}=\frac{\left(  \sum_{i=1}^{m}%
u_{i}^{2}\right)  _{,j}}{2f}=\frac{\sum_{i=1}^{m}u_{i}u_{ij}}{f} & = &
\frac{u_{1}u_{1j}}{f}\\
& = & u_{1j},
\end{array}
\label{ks2}%
\end{equation}
we then obtain%
\begin{equation}%
\begin{array}
[c]{l}%
\Delta u=-\left(  p-2\right)  u_{11}.
\end{array}
\label{ks3}%
\end{equation}
Therefore the inequality (\ref{ks1}) can be written as%
\begin{equation}%
\begin{array}
[c]{lll}%
\sum_{i,j=1}^{m}\left(  u_{ij}\right)  ^{2} & \geq & \left(  u_{11}\right)
^{2}+2\sum_{\alpha=2}^{m}\left(  u_{1\alpha}\right)  ^{2}+\frac{\left(
p-1\right)  ^{2}}{m-1}\left(  u_{11}\right)  ^{2}\\
& = & \left(  1+\frac{\left(  p-1\right)  ^{2}}{m-1}\right)  \left(
u_{11}\right)  ^{2}+2\sum_{\alpha=2}^{m}\left(  u_{1\alpha}\right)  ^{2}\\
& \geq & \left(  1+\kappa\right)  \sum_{j=1}^{m}\left(  u_{1j}\right)  ^{2}\\
& = & \left(  1+\kappa\right)  \left\vert \nabla f\right\vert ^{2}.
\end{array}
\label{ks4}%
\end{equation}
Then (\ref{ks}) follows.

When "$=$" holds in the inequality (\ref{ks}), then by (\ref{ks1}), we have
\[%
\begin{array}
[c]{ll}%
u_{\alpha\beta}=0 & \text{for all }\alpha\neq\beta, \text{ where }\alpha
,\beta=2,\ldots m
\end{array}
\]
and%
\begin{equation}%
\begin{array}
[c]{ll}%
u_{\alpha\alpha}=u_{\beta\beta}\text{ \ \ } & \text{for all }\alpha
,\beta=2,\ldots m.
\end{array}
\label{ks5}%
\end{equation}
Using (\ref{ks3}), (\ref{ks5}) then gives%
\[%
\begin{array}
[c]{ll}%
u_{11}=-\frac{m-1}{p-1}u_{\alpha\alpha}\text{ \ \ } & \text{for all }%
\alpha=2,\ldots m.
\end{array}
\]
Moreover, by (\ref{ks4}),

\begin{itemize}
\item If $\left(  p-1\right)  ^{2}<m-1,$ then $u_{1\alpha}=0$ for all
$\alpha=2,\ldots m.$

\item If $\left(  p-1\right)  ^{2}>m-1,$ then $u_{11}=0,$ i.e. $u_{ii}=0$ for
all $i=1,\ldots m.$
\end{itemize}

Hence we complete the proof. \ 

\endproof

\bigskip

Next, we show two examples to verify Lemma \ref{KS} is sharp.

\begin{example}
If $u\left(  x\right)  =\log\left\vert x\right\vert $ in $\mathbb{R}
^{m}\backslash\left\{  0\right\}  ,$ then it is easy to check that $\Delta
_{m}u=0$ for all $m\geq2.$ Since%
\[%
\begin{array}
[c]{lll}%
\left\vert \nabla du\right\vert ^{2}=\sum_{i,j=1}^{m}\left(  \frac{\delta
_{ij}}{\left\vert x\right\vert ^{2}}-\frac{2x_{i}x_{j}}{\left\vert
x\right\vert ^{4}}\right)  ^{2} & \text{and} & \left\vert \nabla\left\vert
\nabla u\right\vert \right\vert ^{2}=\frac{1}{\left\vert x\right\vert ^{4}},
\end{array}
\text{ }%
\]
we obtain%
\[%
\begin{array}
[c]{l}%
\left\vert \nabla du\right\vert ^{2}=m\left\vert \nabla\left\vert \nabla
u\right\vert \right\vert ^{2}%
\end{array}
\]
for $m\geq2.$ This example implies Lemma \ref{KS} is sharp in the case of
$p=m=2.$
\end{example}

\begin{example}
Let $u\left(  x\right)  =\left\vert x\right\vert ^{\frac{p-m}{p-1}}$ in
$\mathbb{R} ^{m}\backslash\left\{  0\right\}  ,$ $p\neq m,$ then $u$ is a
$p$-harmonic function. Since%
\[%
\begin{array}
[c]{l}%
\left\vert \nabla du\right\vert ^{2}=\left(  \frac{p-m}{p-1}\right)
^{2}\left\vert x\right\vert ^{\frac{2\left(  1-m\right)  }{p-1}-2}\sum
_{i,j=1}^{m}\left\{  \delta_{ij}+\left(  \frac{1-m}{p-1}-1\right)  \frac
{x_{i}x_{j}}{\left\vert x\right\vert ^{2}}\right\}  ^{2}%
\end{array}
\]
and%
\[%
\begin{array}
[c]{l}%
\left\vert \nabla\left\vert \nabla u\right\vert \right\vert ^{2}=\left(
\frac{\left(  p-m\right)  \left(  1-m\right)  }{\left(  p-1\right)  ^{2}%
}\right)  ^{2}\left\vert x\right\vert ^{\frac{2\left(  1-m\right)  }{p-1}-2},
\end{array}
\]
we have%
\[%
\begin{array}
[c]{l}%
\left\vert \nabla du\right\vert ^{2}=\left(  1+\frac{\left(  p-1\right)  ^{2}%
}{m-1}\right)  \left\vert \nabla\left\vert \nabla u\right\vert \right\vert
^{2}.
\end{array}
\]
This example implies Lemma \ref{KS} is sharp in the case of $\left(
p-1\right)  ^{2}\leq m-1.$
\end{example}

\bigskip

\subsection{The Proof of Theorem \ref{T2}}

\smallskip

We need several Lemmas:

\begin{lemma}
\label{S-1}Suppose $M^{m}$ is a complete noncompact Riemannian manifold
satisfying $\left(  P_{\rho}\right)  $ and (\ref{Rs}). Let $u\in C^{3}\left(
M^{m}\right)  $ be a strongly $p$-harmonic function, $p>1,$ $p\neq2.$ Then,
for every $0<\varepsilon<1,$
\begin{equation}%
\begin{array}
[c]{l}%
\int_{B\left(  R\right)  }A_{2}f_{\varepsilon}^{q-2}dv+\varepsilon B_{0}%
\leq\frac{100\cdot B_{1}}{R^{2}}\int_{B\left(  2R\right)  \backslash B\left(
R\right)  }f_{\varepsilon}^{q}dv,
\end{array}
\label{vs3-1}%
\end{equation}
where $f=\left\vert \nabla u\right\vert \, , $ $f_{\varepsilon}=\sqrt
{f^{2}+\varepsilon},$ $q=s+2,$ $q-1+\kappa+b>\varepsilon_{1},$ $b=\min\left\{
0,\left(  p-2\right)  \left(  q-p\right)  \right\}  ,$
\begin{equation}%
\begin{array}
[c]{lll}%
B_{0} & = & \int_{M}\eta^{2}\bigg(f_{\varepsilon}^{s-2}\sum_{i,j=1}^{m}%
u_{ij}^{2}+f_{\varepsilon}^{s-2}\left\langle \nabla u,\nabla\Delta
u\right\rangle \\
&  & +\frac{p-4}{2}f_{\varepsilon}^{s-4}\left\langle \nabla u,\nabla
f_{\varepsilon}^{2}\right\rangle \Delta u-bf_{\varepsilon}^{s-2}\left\vert
\nabla f_{\varepsilon}\right\vert ^{2}\bigg)dv\, ,
\end{array}
\label{vs5}%
\end{equation}
\[%
\begin{array}
[c]{l}%
A_{2}=\frac{4\left(  1-\varepsilon_{2}\right)  \left(  q-1+\kappa
+b-\varepsilon_{1}\right)  }{q^{2}}\rho f_{\varepsilon}^{2} + \sum_{i,j=1}%
^{m}R_{ij}u_{i}u_{j}\, ,
\end{array}
\]
and%
\[%
\begin{array}
[c]{l}%
B_{1}=\frac{\left(  1+\left\vert p-2\right\vert \right)  ^{2}}{\varepsilon
_{1}}+\frac{4\left(  \frac{1}{\varepsilon_{2}}-1\right)  \left(
q-1+\kappa+b-\varepsilon_{1}\right)  }{q^{2}},
\end{array}
\]
for some $0<\varepsilon_{1},\varepsilon_{2}<1.$
\end{lemma}

\proof

Combining Lemma \ref{KS} and Lemma \ref{BS}, and using the formula%
\[%
\begin{array}
[c]{l}%
f^{2}\left\vert \nabla\left(  du\right)  \right\vert ^{2}\geq\frac{1+\kappa
}{4}\left\vert \nabla\left\vert du\right\vert ^{2}\right\vert ^{2}%
\end{array}
\]
holds on all of $M$, we have the following.%
\begin{equation}%
\begin{array}
[c]{lll}%
\frac{1}{2}\mathcal{L}_{s,\varepsilon}\left(  f_{\varepsilon}^{2}\right)  &
\geq & \left(  s+1+\kappa\right)  f_{\varepsilon}^{s}\left\vert \nabla
f_{\varepsilon}\right\vert ^{2}+f_{\varepsilon}^{s}\sum_{i,j=1}^{m}R_{ij}%
u_{i}u_{j}\\
&  & +\frac{\left(  p-2\right)  \left(  s-p+2\right)  }{4}f_{\varepsilon
}^{s-4}\left\langle \nabla u,\nabla f_{\varepsilon}^{2}\right\rangle ^{2}\\
&  & +\varepsilon\bigg(f_{\varepsilon}^{s-2}\sum_{i,j=1}^{m} u_{ij}%
^{2}+f_{\varepsilon}^{s-2}\left\langle \nabla u,\nabla\Delta u\right\rangle \\
&  & +\frac{p-4}{2}f_{\varepsilon}^{s-4}\left\langle \nabla u,\nabla
f_{\varepsilon}^{2}\right\rangle \Delta u\bigg).
\end{array}
\label{bs1}%
\end{equation}

We multiply both sides of (\ref{bs1}) by \ a cut off function $\eta^{2}\in
C_{0}^{\infty}(M)$ and integrate over $M,$
\begin{equation}%
\begin{array}
[c]{lll}
&  & \frac{1}{2}\int_{M}\eta^{2}\mathcal{L}_{s,\varepsilon}\left(
f_{\varepsilon}^{2}\right)  dv\\
& \geq & \int_{M}\eta^{2}f_{\varepsilon}^{s}\left(  ( s+1+\kappa) \left\vert
\nabla f_{\varepsilon}\right\vert ^{2}+\sum_{i,j=1}^{m}R_{ij}u_{i}%
u_{j}\right)  dv\\
&  & +\frac{\left(  p-2\right)  \left(  s-p+2\right)  }{4}\int_{M}\eta
^{2}f_{\varepsilon}^{s-4}\left\langle \nabla u,\nabla f_{\varepsilon}%
^{2}\right\rangle ^{2}dv\\
&  & +\varepsilon\int_{M}\eta^{2}\bigg(f_{\varepsilon}^{s-2}\sum_{i,j=1}%
^{m}u_{ij}^{2}+f_{\varepsilon}^{s-2}\left\langle \nabla u,\nabla\Delta
u\right\rangle \\
&  & +\frac{p-4}{2}f_{\varepsilon}^{s-4}\left\langle \nabla u,\nabla
f_{\varepsilon}^{2}\right\rangle \Delta u\bigg)dv
\end{array}
\label{vs1}%
\end{equation}
where $\eta$ is a cut-off function on $M$ satisfying%
\[
\left\{
\begin{array}
[c]{ll}%
\eta\left(  x\right)  =1 & \text{if }x\in\overline{B\left(  R\right)  },\\
0<\eta\left(  x\right)  <1\text{\ \ \ \ \ } & \text{if }x\in B\left(
2R\right)  \backslash\overline{B\left(  R\right)  },\\
\eta\left(  x\right)  =0 & \text{if }x\in M\backslash B\left(  2R\right)  ,
\end{array}
\right.
\]
and%
\[
\left\{
\begin{array}
[c]{ll}%
\left\vert \nabla\eta\left(  x\right)  \right\vert =0 & \text{if }x\in
B\left(  R\right)  \text{ or }x\in M\backslash B\left(  2R\right)  ,\\
\left\vert \nabla\eta\left(  x\right)  \right\vert \leq\frac{10}%
{R}\text{\ \ \ \ \ } & \text{if }x\in B\left(  2R\right)  \backslash
\overline{B\left(  R\right)  },
\end{array}
\right.
\]

Since integration by parts and Cauchy-Schwarz inequality assert that%
\[%
\begin{array}
[c]{lll}%
\frac{1}{2}\int_{M}\eta^{2}\mathcal{L}_{s,\varepsilon}\left(  f_{\varepsilon
}^{2}\right)  dv & = & \frac{-1}{2}\int_{M}\left\langle \nabla\eta
^{2},f_{\varepsilon}^{s}\mathcal{\nabla}f_{\varepsilon}^{2}+\left(
p-2\right)  f_{\varepsilon}^{s-2}\left\langle \nabla u,\mathcal{\nabla
}f_{\varepsilon}^{2}\right\rangle \nabla u\right\rangle dv\\
& \leq & 2\int_{M}\eta\left\vert \nabla\eta\right\vert \left(  f_{\varepsilon
}^{s+1}\left\vert \mathcal{\nabla}f_{\varepsilon}\right\vert +\left(
p-2\right)  f_{\varepsilon}^{s-1}f^{2}\left\vert \mathcal{\nabla
}f_{\varepsilon}\right\vert \right)  dv\\
& \leq & 2\left(  1+\left\vert p-2\right\vert \right)  \int_{M}\eta\left\vert
\nabla\eta\right\vert f_{\varepsilon}^{s+1}\left\vert \mathcal{\nabla
}f_{\varepsilon}\right\vert dv\\
& \leq & \varepsilon_{1}\int_{M}\eta^{2}f_{\varepsilon}^{s}\left\vert
\mathcal{\nabla}f_{\varepsilon}\right\vert ^{2}dv+\frac{\left(  1+\left\vert
p-2\right\vert \right)  ^{2}}{\varepsilon_{1}}\int_{M}\left\vert \nabla
\eta\right\vert ^{2}f_{\varepsilon}^{s+2}dv,
\end{array}
\]
where $0<\varepsilon_{1}<1$ is a positive constant such that $q-1+\kappa
+b>\varepsilon_{1}.$

On the other hand,%
\[%
\begin{array}
[c]{lll}
&  & \frac{\left(  p-2\right)  \left(  s-p+2\right)  }{4}\int_{M}\eta
^{2}f_{\varepsilon}^{s-4}\left\langle \nabla u,\nabla f_{\varepsilon}%
^{2}\right\rangle ^{2}dv\\
& \geq & \frac{b}{4}\int_{M}\eta^{2}f_{\varepsilon}^{s-4}\left\vert \nabla
u\right\vert ^{2}\left\vert \nabla f_{\varepsilon}^{2}\right\vert ^{2}dv\\
& = & b\int_{M}\eta^{2}f_{\varepsilon}^{s-2}f^{2}\left\vert \nabla
f_{\varepsilon}\right\vert ^{2}dv\\
& = & b\int_{M}\eta^{2}f_{\varepsilon}^{s}\left\vert \nabla f_{\varepsilon
}\right\vert ^{2}dv-b\varepsilon\int_{M}\eta^{2}f_{\varepsilon}^{s-2}%
\left\vert \nabla f_{\varepsilon}\right\vert ^{2}dv
\end{array}
\]
where%
\[%
\begin{array}
[c]{l}%
b=\min\left\{  0,\left(  p-2\right)  \left(  s-p+2\right)  \right\}  .
\end{array}
\]
Then (\ref{vs1}) implies%

\begin{equation}%
\begin{array}
[c]{lll}
&  & A_{1}\int_{M}\eta^{2}f_{\varepsilon}^{s}\left\vert \nabla f_{\varepsilon
}\right\vert ^{2}dv + \int_{M}\eta^{2}f_{\varepsilon}^{s}\sum_{i,j=1}%
^{m}R_{ij}u_{i}u_{j}dv+\varepsilon B_{0}\\
& \leq & \frac{\left(  1+\left\vert p-2\right\vert \right)  ^{2}}%
{\varepsilon_{1}}\int_{M}\left\vert \nabla\eta\right\vert ^{2}f_{\varepsilon
}^{s+2}dv,
\end{array}
\label{vs1-2}%
\end{equation}
where $A_{1}=s+1+\kappa+b-\varepsilon_{1}>0\, .$

Now we compute the first term in the left hand side of (\ref{vs1-2}). Since
$q=s+2\, ,$
\[%
\begin{array}
[c]{ll}
& \int_{M}\eta^{2}f_{\varepsilon}^{s}\left\vert \nabla f_{\varepsilon
}\right\vert ^{2}dv\\
= & \frac{4}{q^{2}}\int_{M}\eta^{2}\left\vert \nabla f_{\varepsilon}^{\frac
{q}{2}}\right\vert ^{2}dv\\
= & \frac{4}{q^{2}}\int_{M}\left\vert \nabla\left(  \eta f_{\varepsilon
}^{\frac{q}{2}}\right)  -\left(  \nabla\eta\right)  f_{\varepsilon}^{\frac
{q}{2}}\right\vert ^{2}dv\\
= & \frac{4}{q^{2}}\int_{M}\left\vert \nabla\left(  \eta f_{\varepsilon
}^{\frac{q}{2}}\right)  \right\vert ^{2}-2\left\langle \nabla\left(  \eta
f_{\varepsilon}^{\frac{q}{2}}\right)  ,f_{\varepsilon}^{\frac{q}{2}}\nabla
\eta\right\rangle +\left\vert \nabla\eta\right\vert ^{2}f_{\varepsilon}%
^{q}dv\\
\geq & \frac{4\left(  1-\varepsilon_{2}\right)  }{q^{2}}\int_{M}\left\vert
\nabla\left(  \eta f_{\varepsilon}^{\frac{q}{2}}\right)  \right\vert
^{2}+\frac{4\left(  1-\frac{1}{\varepsilon_{2}}\right)  }{q^{2}}\int_{M^{m}%
}\left\vert \nabla\eta\right\vert ^{2}f_{\varepsilon}^{q}dv.
\end{array}
\]
where $\varepsilon_{2}$ is a positive constant satisfying $0<\varepsilon
_{2}<1.$ Thus, we have%
\begin{equation}%
\begin{array}
[c]{lll}
&  & \frac{4\left(  1-\varepsilon_{2}\right)  A_{1}}{q^{2}}\int_{M}\left\vert
\nabla\left(  \eta f_{\varepsilon}^{\frac{q}{2}}\right)  \right\vert ^{2}dv+
\int_{M}\eta^{2}f_{\varepsilon}^{q-2}\sum_{i,j=1}^{m}R_{ij}u_{i}%
u_{j}\,dv+\varepsilon B_{0}\\
& \leq & \left(  \frac{\left(  1+\left\vert p-2\right\vert \right)  ^{2}%
}{\varepsilon_{1}}+\frac{4\left(  \frac{1}{\varepsilon_{2}}-1\right)  A_{1}%
}{q^{2}}\right)  \int_{M}\left\vert \nabla\eta\right\vert ^{2}f_{\varepsilon
}^{q}dv.
\end{array}
\label{vs2}%
\end{equation}

According to weighted Poincar\'{e} inequality
\[%
\begin{array}
[c]{l}%
\int_{M}\rho\Psi^{2}dv\leq\int_{M}\left\vert \nabla\Psi\right\vert ^{2}dv,
\end{array}
\]
if we select $\Psi=\eta f_{\epsilon}^{\frac{q}{2}},$ then (\ref{vs2}) implies%
\[%
\begin{array}
[c]{l}%
\int_{B\left(  R\right)  }A_{2}f_{\varepsilon}^{q-2}dv+\varepsilon B_{0}%
\leq\frac{100\cdot B_{1}}{R^{2}}\int_{B\left(  2R\right)  \backslash B\left(
R\right)  }f_{\varepsilon}^{q}dv,
\end{array}
\]
for all fixed $R>0.$\ \endproof

\bigskip

\begin{lemma}
\label{S}Let $B_{0}$ be as in $(\ref{vs5})$, $p>1,$ $p\neq2,$ $q=s+2$ and
\[%
\begin{array}
[c]{l}%
b=\min\left\{  0,\left(  p-2\right)  \left(  q-p\right)  \right\}  .
\end{array}
\]
Then \newline(i) if $q>2,$ then $\varepsilon B_{0}\rightarrow0$ as
$\varepsilon\rightarrow0,$ \newline(ii) if $1<q\leq2$ and $b\leq-\frac{\left(
p-4\right)  ^{2}m}{4},$ then $\varepsilon B_{0}\geq0$ as $\varepsilon
\rightarrow0.$ \newline
\end{lemma}

\proof

First of all, we derive some properties.

For $s\geq2,$ it is easy to check that
\begin{equation}%
\begin{array}
[c]{l}%
\varepsilon f_{\varepsilon}^{s-2}\rightarrow0\text{ as }\varepsilon
\rightarrow0.
\end{array}
\label{s00}%
\end{equation}

If $0<s<2,$ then we also have%
\begin{equation}%
\begin{array}
[c]{l}%
\varepsilon f_{\varepsilon}^{s-2}=\frac{\varepsilon}{f_{\varepsilon}^{2-s}%
}\leq\frac{\varepsilon}{\varepsilon^{1-s/2}}=\varepsilon^{s/2}\rightarrow
0\text{ as }\varepsilon\rightarrow0,
\end{array}
\label{s0}%
\end{equation}

By using the estimates%
\[%
\begin{array}
[c]{lll}%
\varepsilon f_{\varepsilon}^{s-4}\left\langle \nabla u,\nabla f_{\varepsilon
}^{2}\right\rangle  & = & 2\varepsilon f_{\varepsilon}^{s-4}\sum_{i,j=1}%
^{m}u_{ij}u_{i}u_{j}\\
& \leq & 2\varepsilon f_{\varepsilon}^{s-4}\sup_{i,j=1,\cdots,m}\left\vert
u_{ij}\right\vert \sum_{i,j=1}^{m}\left\vert u_{i}u_{j}\right\vert \\
& \leq & 2m\varepsilon f_{\varepsilon}^{s-4}f^{2}\sup_{i,j=1,\cdots
,m}\left\vert u_{ij}\right\vert \\
& \leq & 2m\varepsilon f_{\varepsilon}^{s-2}\sup_{i,j=1,\cdots,m}\left\vert
u_{ij}\right\vert \text{ \ \ }%
\end{array}
\]
and%
\[%
\begin{array}
[c]{lll}%
\varepsilon f_{\varepsilon}^{s-2}\left\vert \nabla f_{\varepsilon}\right\vert
^{2} & = & \frac{\varepsilon}{4}f_{\varepsilon}^{s-4}\left\vert \nabla
f_{\varepsilon}^{2}\right\vert ^{2}\\
& = & \varepsilon f_{\varepsilon}^{s-4}\sum_{i,j,k=1}^{m}u_{ik}u_{kj}%
u_{i}u_{j}\\
& \leq & 2m\varepsilon f_{\varepsilon}^{s-4}f^{2}\sup_{i,j,k=1,\cdots
,m}\left\vert u_{ik}\right\vert \left\vert u_{kj}\right\vert \\
& \leq & 2m\varepsilon f_{\varepsilon}^{s-2}\sup_{i,j,k=1,\cdots,m}\left\vert
u_{ik}\right\vert \left\vert u_{kj}\right\vert
\end{array}
\]
then (\ref{s00}) and (\ref{s0}) imply
\begin{equation}
\left\{
\begin{array}
[c]{l}%
\varepsilon f_{\varepsilon}^{s-4}\left\vert \left\langle \nabla u,\nabla
f_{\varepsilon}^{2}\right\rangle \right\vert \rightarrow0,\\
\varepsilon f_{\varepsilon}^{s-2}\left\vert \nabla f_{\varepsilon}\right\vert
^{2}\rightarrow0,
\end{array}
\right.  \label{s1}%
\end{equation}
as $\varepsilon\rightarrow0,$ for all $s>0.$

In the case $-1<s\leq0,$%
\begin{equation}%
\begin{array}
[c]{l}%
\varepsilon ff_{\varepsilon}^{s-2}\leq\frac{\varepsilon}{\left(
f^{2}+\varepsilon\right)  ^{1/2-s/2}}\leq\varepsilon^{1/2+s/2}\rightarrow
0\text{ as }\varepsilon\rightarrow0.
\end{array}
\label{s3}%
\end{equation}

Now we prove Lemma as follows.

For any fixed $s>0,$ by (\ref{s00}), (\ref{s0}) and (\ref{s1}), then we
obtain,
\[%
\begin{array}
[c]{lll}%
\left\vert \varepsilon B_{0} \right\vert  & = & \varepsilon\bigg{|} \int
_{M}\eta^{2}\bigg(f_{\varepsilon}^{s-2}\sum_{i,j=1}^{m}u_{ij}^{2}%
+f_{\varepsilon}^{s-2}\left\langle \nabla u,\nabla\Delta u\right\rangle \\
&  & +\frac{p-4}{2}f_{\varepsilon}^{s-4}\left\langle \nabla u,\nabla
f_{\varepsilon}^{2}\right\rangle \Delta u-bf_{\varepsilon}^{s-2}\left\vert
\nabla f_{\varepsilon}\right\vert ^{2}\bigg)dv \bigg{|}\\
& \le & \int_{M}\eta^{2}\bigg(\left(  \varepsilon f_{\varepsilon}%
^{s-2}\right)  \sum_{i,j=1}^{m}u_{ij}^{2}+\left(  \varepsilon f_{\varepsilon
}^{s-2}\right)  \left\vert \left\langle \nabla u,\nabla\Delta u\right\rangle
\right\vert \\
&  & +\frac{|p-4|}{2}\left(  \varepsilon f_{\varepsilon}^{s-4}\left\vert
\left\langle \nabla u,\nabla f_{\varepsilon}^{2}\right\rangle \right\vert
\right)  \left\vert \Delta u\right\vert -b\left(  \varepsilon f_{\varepsilon
}^{s-2}\left\vert \nabla f_{\varepsilon}\right\vert ^{2}\right)  \bigg)dv\\
& \rightarrow & 0,\text{ as }\varepsilon\rightarrow0.
\end{array}
\]

If $s>-1,$ since $b\leq-\frac{\left(  p-4\right)  ^{2}m}{4},$ then%
\[%
\begin{array}
[c]{lll}
&  & f_{\varepsilon}^{s-2}\sum_{i,j=1}^{m}u_{ij}^{2}-bf_{\varepsilon}%
^{s-2}\left\vert \nabla f_{\varepsilon}\right\vert ^{2}+\frac{p-4}%
{2}f_{\varepsilon}^{s-4}\left\langle \nabla u,\nabla f_{\varepsilon}%
^{2}\right\rangle \Delta u\\
& \geq & f_{\varepsilon}^{s-2}\sum_{i,j=1}^{m}u_{ij}^{2}-bf_{\varepsilon
}^{s-2}\left\vert \nabla f_{\varepsilon}\right\vert ^{2}-\left\vert
p-4\right\vert f_{\varepsilon}^{s-2}\left\vert \nabla f_{\varepsilon
}\right\vert \left\vert \Delta u\right\vert \\
& \geq & f_{\varepsilon}^{s-2}\sum_{i,j=1}^{m}u_{ij}^{2}-bf_{\varepsilon
}^{s-2}\left\vert \nabla f_{\varepsilon}\right\vert ^{2}-\frac{f_{\varepsilon
}^{s-2}\left(  \Delta u\right)  ^{2}}{m}-\frac{\left(  p-4\right)  ^{2}m}%
{4}f_{\varepsilon}^{s-2}\left\vert \nabla f_{\varepsilon}\right\vert ^{2}\\
& \geq & 0,
\end{array}
\]
here we use $\sum_{i,j=1}^{m}u_{ij}^{2}\geq\frac{\left(  \Delta u\right)
^{2}}{m}.$ Hence by (\ref{s00}), (\ref{s0}) and (\ref{s3}),%
\[%
\begin{array}
[c]{lll}%
\varepsilon B_{0} & = & \varepsilon\int_{M}\eta^{2}\bigg(f_{\varepsilon}%
^{s-2}\sum_{i,j=1}^{m}u_{ij}^{2}+f_{\varepsilon}^{s-2}\left\langle \nabla
u,\nabla\Delta u\right\rangle \\
&  & +\frac{p-4}{2}f_{\varepsilon}^{s-4}\left\langle \nabla u,\nabla
f_{\varepsilon}^{2}\right\rangle \Delta u-bf_{\varepsilon}^{s-2}\left\vert
\nabla f_{\varepsilon}\right\vert ^{2}\bigg)dv\\
& \geq & \varepsilon\int_{M}\eta^{2}f_{\varepsilon}^{s-2}\left\langle \nabla
u,\nabla\Delta u\right\rangle \\
& \geq & -\int_{M}\eta^{2}\left(  \varepsilon f_{\varepsilon}^{s-2}f\right)
\left\vert \nabla\Delta u\right\vert \\
& \rightarrow & 0\text{ whenever }s>-1\text{ and }\varepsilon\rightarrow0.
\end{array}
\]

In particular, if $s>-1$ and $p=4,$ by applying (\ref{s00}), (\ref{s0}),
(\ref{s3}) and $b\leq0$, then%
\[%
\begin{array}
[c]{lll}%
\varepsilon B_{0} & = & \varepsilon\int_{M}\eta^{2}\bigg(f_{\varepsilon}%
^{s-2}\sum_{i,j=1}^{m}u_{ij}^{2}+f_{\varepsilon}^{s-2}\left\langle \nabla
u,\nabla\Delta u\right\rangle -bf_{\varepsilon}^{s-2}\left\vert \nabla
f_{\varepsilon}\right\vert ^{2}\bigg)dv\\
& \geq & -\int_{M}\eta^{2}\left(  \varepsilon f_{\varepsilon}^{s-2}f\right)
\left\vert \nabla\Delta u\right\vert dv\\
& \rightarrow & 0\text{ whenever }s>-1\text{ and }\varepsilon\rightarrow0.
\end{array}
\]
\endproof

\bigskip

\begin{remark}
In Lemma \ref{S}, if $p=4$ and $q>1,$ then $\varepsilon B_{0}\geq0$ as
$\varepsilon\rightarrow0.$
\end{remark}

\begin{proof}
[Proof of Theorem \ref{T2}]Since we assume $q-1+\kappa+b>0,$ the curvature
condition (\ref{Rs}) means that there exists a constant $0<\delta<1$ such
that
\begin{equation}%
\begin{array}
[c]{lll}%
Ric_{M}(\nabla u, \nabla u) = \sum_{i,j=1}^{m}R_{ij}u_{i}u_{j} & \geq &
-\frac{4\left(  q-1+\kappa+b\right)  }{q^{2} }\delta\rho f^{2}.
\end{array}
\label{5.21}%
\end{equation}
To apply Lemmas \ref{S-1} and \ref{S}, we need the following conditions:
\[
(\ast)\left\{
\begin{array}
[c]{l}%
q>2\text{ and }q-1+\kappa+b>0,\\
\text{ or}\quad1<q\leq2\text{ and }q-1+\kappa+b>0, \text{ where }b\leq
-\frac{\left(  p-4\right)  ^{2}m}{4}.
\end{array}
\right.
\]
We first assume $p\neq2.$

\noindent For $p>2\, ,$ the expression $1<q\leq2$ implies that $q < p\, .$
Hence $b=(p-2)(q-p)\, .$ Then
\[
q-1+\kappa+b>0\quad\Longleftrightarrow\quad q>p-1-\frac{\kappa}{p-1}
\]
(cf. Remark \ref{R: 5.10}), and
\[
b\leq-\frac{\left(  p-4\right)  ^{2}m}{4} \quad\Longleftrightarrow\quad q \le
p-\frac{\left(  p-4\right)  ^{2}m}{4\left(  p-2\right)  }\, .
\]
That is, for $p>2\, ,$ $(\ast)$ can be rewritten as
\[
(\ast_{1}) \qquad\left\{
\begin{array}
[c]{l}%
\max\left\{  2,1-\kappa-b\right\}  <q\\
\text{or }\quad\max\left\{  1,p-1-\frac{\kappa}{p-1}\right\}  <q\leq
\min\left\{  2,p-\frac{\left(  p-4\right)  ^{2}m}{4\left(  p-2\right)
}\right\}  \text{ }%
\end{array}
\right.
\]
For $p=4\, ,$ $b\leq-\frac{\left(  p-4\right)  ^{2}m}{4} = 0$ holds and
$(\ast)$ can be simplified as follows:
\[
(\ast_{2})\qquad\max\left\{  1,1-\kappa-b\right\}  < q\, .
\]

\noindent For $1<p<2,$ the expression $1<q\leq2$ implies that $p < q\, .$ Or
$q \le p\, (< 2)\,$ would lead to $0=b\le-\frac{\left(  p-4\right)  ^{2}m}{4}
< 0\, ,$ a contradiction. Hence $b=(p-2)(q-p)\, .$ Then $q-1+\kappa
+b>0\quad\big (\Longleftrightarrow\quad q>p-1-\frac{\kappa}{p-1}%
\big )\quad\operatorname{holds} \, .$ However, $b\leq-\frac{\left(
p-4\right)  ^{2}m}{4} \quad\big (\Longleftrightarrow\quad q \le p-\frac
{\left(  p-4\right)  ^{2}m}{4\left(  p-2\right)  }\big )\quad\operatorname{is}
\quad\operatorname{invalid} \, .$ What remains is the following:

\noindent For $1<p<2,$ the expression $2<q$ implies that $b=(p-2)(q-p)\, .$
Then $q-1+\kappa+b>0\quad\big (\Longleftrightarrow\quad q>p-1-\frac{\kappa
}{p-1}\big )\quad\operatorname{holds}\, . $

Thus, for $1<p<2,$ $(\ast)$ can be rewritten as
\[
(\ast_{3}) \qquad2 < q.
\]

\noindent Similarly, for $p=2\, ,$ we have $b = 0$ and $\kappa= \frac{1}%
{m-1}\, .$ It follows that $q-1+\kappa+b>0$ holds if and only if%

\[
(\ast_{4})\qquad\frac{m-2}{m-1} < q\, .
\]

In view of $(\ast_{1}), (\ast_{2}), (\ast_{3}), (\ast_{4})\, ,$ Lemmas
\ref{S-1} and \ref{S}, we obtain (\ref{vs3-1}). As $\epsilon\to0\, ,$
(\ref{vs3-1}), via (\ref{5.21}) tends to
\begin{equation}%
\begin{array}
[c]{lll}%
\int_{B\left(  R\right)  }A_{3}f^{q}dv & \leq & \frac{100\cdot B_{1}}{R^{2}%
}\int_{B\left(  2R\right)  \backslash B\left(  R\right)  }f^{q}dv,
\end{array}
\label{vs3}%
\end{equation}
where%
\[%
\begin{array}
[c]{lll}%
A_{3} & = & \left(  \frac{4\left(  1-\varepsilon_{2}\right)  \left(
q-1+\kappa+b-\varepsilon_{1}\right)  }{q^{2}}-\frac{4\left(  q-1+\kappa
+b\right)  \delta}{q^{2}}\right)  \rho.
\end{array}
\]
\newline Hence one has $A_{3}>0$ whenever we select $\varepsilon_{1}$ and
$\varepsilon_{2}$ small enough. Suppose $f\in L^{q}\left(  M\right)  $, then
the right hand side of (\ref{vs3}) tends to zero as $R\rightarrow\infty$, and
then we conclude that $f(x) = 0$ for all $x\in M$ and for some $0<\delta<1,$
i.e. $u\left(  x\right)  $ is a constant on $M$ for some $0<\delta<1$.

In particular, if $1<p<2,$ since constant function is not a strongly
$p$-harmonic function, then such $u$ does not exist.
\end{proof}

\begin{remark}
\label{R: 5.10} If $p>2$ and $p\geq q,$ then
\[%
\begin{array}
[c]{lll}%
q-1+\kappa+b & = & q-1+\kappa+\left(  p-2\right)  \left(  q-p\right) \\
& = & \left(  p-1\right)  q-\left(  p-1\right)  ^{2}+\kappa>0,
\end{array}
\]
whenever $q>p-1-\frac{\kappa}{p-1}.$
\end{remark}

\begin{remark}
If we replace the finite $q$-energy by $\int_{B\left(  2R\right)  \backslash
B\left(  R\right)  }\left\vert \nabla u\right\vert ^{q}dv=o\left(
R^{2}\right)  $ as $R\rightarrow\infty,$ then Theorem \ref{T2} is still valid.
\end{remark}

\begin{remark}
Since $\left(  P_{\lambda_{q}}\right)  $ implies $\left(  P_{\lambda_{p}%
}\right)  $ for all $p>q$ (cf. \cite{HKM}). If $M$ satisfies $\left(
P_{\lambda_{2}}\right)  ,$ by using Lemma \ref{Volume estimate}, then
$2$-hyperbolic end is equality to $p$-hyperbolic end since this end has
infinite volume. Hence we may use the method of Theorem 2.1 of \cite{LW1} to
refine the conditions of Theorem \ref{T2} whenever $M$ satisfies $\left(
P_{\lambda_{2}}\right)  .$ But we omit it in this paper.
\end{remark}

\begin{corollary}
Let $M^{m}$ be a complete noncompact Riemannian manifold satisfying $\left(
P_{\rho}\right)  $ and (\ref{Rs}), where%
\[%
\begin{array}
[c]{lll}%
\tau & < & \frac{4\left(  \left(  p-1\right)  q-\left(  p-1\right)
^{2}+\kappa\right)  }{q^{2}},
\end{array}
\]
$\kappa=\min\{\frac{\left(  p-1\right)  ^{2}}{m-1},1\},$ $p>2,$ $p\geq q.$ Let
$u\in C^{3}\left(  M^{m}\right)  $ be a strongly $p$-harmonic function, with
finite $q$-energy $E_{q}\left(  u\right)  .$ Then $u$ is a constant if $p$ and
$q$ satisfy one of the following: \newline(1) $p=4,$ $q>\frac{9-\kappa}{3},$
\newline(2) $p\neq4,$ and either%
\[%
\begin{array}
[c]{l}%
\max\left\{  1,p-1-\frac{\kappa}{p-1}\right\}  <q\leq\min\left\{
2,p-\frac{\left(  p-4\right)  ^{2}m}{4\left(  p-2\right)  }\right\}
\end{array}
\]
or
\[%
\begin{array}
[c]{l}%
\max\{2, p-1-\frac{\kappa}{p-1}\} < q.
\end{array}
\]
In particular, if $p=q,$ then every strongly $p$-harmonic function $u$ with
finite $p$-energy is constant.
\end{corollary}

\begin{corollary}
Let $M^{m}$ be a complete noncompact Riemannian manifold satisfying $\left(
P_{\rho}\right)  $ and (\ref{Rs}), where%
\[%
\begin{array}
[c]{lll}%
\tau & < & \frac{4\left(  p-1+\kappa\right)  }{p^{2}},
\end{array}
\]
$\kappa=\min\{\frac{\left(  p-1\right)  ^{2}}{m-1},1\}.$ If $u\in C^{3}\left(
M^{m}\right)  $ is a strongly $p$-harmonic function for $p\geq2,$ with
$E_{p}\left(  u\right)  <\infty,$ then $u$ is a constant.
\end{corollary}

\begin{remark}
According to the following Lemma \ref{ws}, we can replace ``Let $u\in
C^{3}\left(  M^{m}\right)  $ be a strongly $p$-harmonic function for $1 < p <
\infty\, .$" in Theorem \ref{T2} by ``Let $u\in C^{2}\left(  M^{m}\right)  $
be a weakly $p$-harmonic function for $p\in\left\{  2\right\}  \cup
\lbrack4,\infty),$ and $u\in C^{3}\left(  M^{m}\right)  $ be a strongly $p
$-harmonic function for $p\in\left(  1,2\right)  \cup\left(  2,4\right)  \,
.$" Theorem \ref{T2} remains to be true.
\end{remark}

\begin{lemma}
\label{ws}If $u\in C^{2}\left(  M\right)  \,($resp. $u\in C^{0}\left(
M\right)  \,)$ is a weakly $p$-harmonic function for $p\in\lbrack4,\infty
)\,($resp. $p=2\,)\,,$ then $u$ is a strongly $p$-harmonic function.
\end{lemma}

%

\proof
By assumption, $u$ satisfies
\[%
\begin{array}
[c]{l}%
\int_{M}\left\langle f^{p-2}\nabla u,\nabla\eta\right\rangle \,dv=0
\end{array}
\]
for every $\eta\in C_{0}^{\infty}\left(  M\right)  ,$ where $f=\left\vert
\nabla u\right\vert .$ Since $u\in C^{2}\left(  M\right)  \,,$ and either
$p=2\,,$ or $p\geq4,$ we have $f^{p-2}\in C^{1}\left(  M\right)  .$ Hence
$f^{p-2}\nabla u\in C^{1}\left(  M\right)  ,$ and the divergence theorem
implies%
\[%
\begin{array}
[c]{l}%
0=\int_{M}\left\langle f^{p-2}\nabla u,\nabla\eta\right\rangle \,dv=-\int
_{M}\text{div}\left(  f^{p-2}\nabla u\right)  \,\eta\,dv
\end{array}
\]
for every $\eta\in C_{0}^{\infty}\left(  M\right)  .$ This completes the
proof. \
\endproof

\subsection{Application to $p$-harmonic morphism\label{morphism}}

A $C^{2}$ map $u:M\rightarrow N$ is called a $p$-harmonic morphism if for any
$p$-harmonic function $f$ defined on an open set $V$ of $N$, the composition
$f\circ u$ is $p$-harmonic on $u^{-1}(V)$. Examples of $p$-harmonic morphisms
include the Hopf fibrations. E. Loubeau and J. M. Burel (\cite{BL}) and E.
Loubeau(\cite{Lo}) prove that a $C^{2}$ map $u:M\rightarrow N$ is a
$p$-harmonic morphism with $p\in(1,\infty)$ if and only if $u$ is a
$p$-harmonic and horizontally weak conformal map. We recall a $C^{2}$ map $u:
M \to N$ is \emph{horizontally weak conformal} if for any $x$ such that $du(x)
\ne0$, the restriction of $du(x)$ to the orthogonal complement of Ker $du(x)$
is conformal and surjective.

\begin{theorem}
\label{T:13} Let $M^{m}$ be a complete noncompact Riemannian manifold,
satisfying $\left(  P_{\rho}\right)  $ and (\ref{Rs}), where
\[
\tau<\frac{4\left(  q-1+\kappa+b\right)  }{q^{2}},\quad\kappa=\min
\{\frac{\left(  p-1\right)  ^{2}}{m-1},1\},\quad\operatorname{and}\quad
b=\min\{0,(p-2)(q-p)\}.
\]
Let $u\in C^{3}\left(  M^{m}, \mathbb{R}^{k}\right)  $ is a $p$-harmonic
morphism $\,u:M^{m}\rightarrow\mathbb{R}^{k},$ $k>0$ of finite $q$-energy
$E_{q}\left(  u\right)  <\infty.$ \newline(I). Then $u$ is constant under one
of the following: \newline(1) $p=2$ and $q>\frac{m-2}{m-1},$ \newline(2)
$p=4,$ $q>1$ and $q-1+\kappa+b>0,$ \newline(3) $p>2,$ $p\neq4,$ and either%
\[%
\begin{array}
[c]{l}%
\max\left\{  1,p-1-\frac{\kappa}{p-1}\right\}  <q\leq\min\left\{
2,p-\frac{\left(  p-4\right)  ^{2}m}{4\left(  p-2\right)  }\right\}
\end{array}
\]
or%
\[%
\begin{array}
[c]{l}%
\max\left\{  2,1-\kappa-b\right\}  <q.
\end{array}
\]
\newline(II). Then $u$ does not exit under \newline(4) $1<p<2,$ $q>2.$
\end{theorem}

\begin{lemma}
\cite{WLW}Let $M,N$ and $K$ be manifolds of dimension $m,$ $n,$ and $k$
respectively, and $u:M\rightarrow N$, and $w:N\rightarrow K$ be $C^{2}$. If
$u$ is horizontally weak conformal, then $|d(w\circ u)|^{p-2}=(\frac{1}%
{n})^{\frac{p-2}{2}}|dw|^{p-2}|du|^{p-2}.$\label{L:14}
\end{lemma}

\begin{proof}
[Proof of Theorem \ref{T:13}.]Let $u^{i}=\pi_{i}\circ u\,,$ where $\pi
_{i}:\mathbb{R}^{k}\rightarrow\mathbb{R}$ is the $i$-th projection. Then the
linear function $\pi_{i}$ is a $p$-harmonic function (cf. 2.2 in \cite{W} ).
Hence $u^{i}\,,$ a composition of a $p$-harmonic morphism and a $p$-harmonic
function is $p$-harmonic. Since $u$ is horizontally weak conformal, it follows
from Lemma \ref{L:14} that $E_{p}(u)<\infty$ implies $E_{p}(u^{i})<\infty\,.$
Now apply $u^{i}$ to Theorem \ref{T2}, the assertion follows.
\endproof

These results are in contrast to the following:

\begin{theorem}
\noindent\cite{WLW} If $u:M^{m}\rightarrow\mathbb{R}^{k},$ $k>0,$ is a
$p$-harmonic morphism, and if there exists $i$, such that $u^{i}=\pi_{i}\circ
u\,$ is $p$-finite, i.e.%
\[%
\begin{array}
[c]{l}%
\liminf_{r\rightarrow\infty}\frac{1}{r^{p}}\int_{B\left(  r\right)
}\left\vert u^{i}\right\vert ^{q}dv<\infty
\end{array}
\]
where $B\left(  r\right)  $ is a geodesic ball of radius $r$, for some
$q>p-1.$ Then $u$ must be constant.
\end{theorem}

As further applications, one obtains

\begin{theorem}
\label{T:14} Let $M^{m}$ be a complete noncompact Riemannian manifold,
satisfying $\left(  P_{\rho}\right)  $ and (\ref{Rs}), where
\[
\tau<\frac{4\left(  q-1+\kappa+b\right)  }{q^{2}},\quad\kappa=\min
\{\frac{\left(  p-1\right)  ^{2}}{m-1},1\}\quad\operatorname{and}\quad
b=\min\{0,(p-2)(q-p)\}\, .
\]
Let $u\in C^{3}\left( M^{m}, \mathbb{R}^{k}\right)  $ be a $p$-harmonic
morphism $\,u:M^{m}\rightarrow\mathbb{R}^{k},$ $k>0,$ and $f:u\left(
M^{m}\right)  \subset%
\mathbb{R}
^{k}\rightarrow%
\mathbb{R}
$ be a nonconstant $p$-harmonic function. Assume $f\circ u$ has finite
$q$-energy $E_{q}\left(  f\circ u\right)  < \infty.$ \newline(I). Then $u$ is
constant under one of the following: \newline(1) $p=2$ and $q>\frac{m-2}%
{m-1},$ \newline(2) $p=4,$ $q>1$ and $q-1+\kappa+b>0,$ \newline(3) $p>2,$
$p\neq4,$ and either%
\[%
\begin{array}
[c]{l}%
\max\left\{  1,p-1-\frac{\kappa}{p-1}\right\}  <q\leq\min\left\{
2,p-\frac{\left(  p-4\right)  ^{2}m}{4\left(  p-2\right)  }\right\}
\end{array}
\]
or%
\[%
\begin{array}
[c]{l}%
\max\left\{  2,1-\kappa-b\right\}  <q.
\end{array}
\]
\newline(II). Then $u$ does not exit under \newline(4) $1<p<2,$ $q>2.$
\end{theorem}

\begin{lemma}
\label{open}A nonconstant $p$-harmonic morphism $u:M^{m}\rightarrow
\mathbb{R}^{k}$ is an open map.
\end{lemma}

\begin{proof}
[Proof of Theorem \ref{T:14}.]Since $u$ is a $p$-harmonic morphism, then
$f\circ u$ is a $p$-harmonic function on $M^{m}.$ According to Theorem
\ref{T2}, then $f\circ u$ is a constant $c$. On the other hand, due to Lemma
\ref{open}, $u$ and $f$ are open maps whenever they are not constant. Now we
assume that $u$ is not constant, then the image of $u$ is an open set
$u\left(  M\right)  \subset%
\mathbb{R}
^{k}$. Hence $f\circ u\left(  M^{m}\right)  $ is an open set. This gives a
contradiction to $f\circ u\left(  M^{m}\right)  =c.$ Then we conclude that $u$
is a constant.
\endproof

\begin{theorem}
\noindent$($Picard Theorem for $p$-harmonic morphisms$)$. Let $M^{m}$ be as in
Theorem \ref{T:14}. Suppose that $u\in C^{3}\left( M^{m}, \mathbb{R}%
^{k}\backslash\{y_{0}\}\right)  $ is a $p$-harmonic morphism $\,u:M^{m}%
\rightarrow\mathbb{R}^{k}\backslash\{y_{0}\},$ and the function $x\mapsto
|u(x)-y_{0}|^{\frac{p-k}{p-1}}$ has finite $q$-energy where $p\neq k$, for $p$
and $q$ satisfying one of the following: $(1)$, $(2) $, and $(3)$ as in
Theorem \ref{T:14}. Then $u$ is constant. For $p$ and $q$ satisfying $(4)$ as
in Theorem \ref{T:14}, then $u$ does not exist.\textit{\ }\smallskip
\end{theorem}

%

\proof
Since $y \mapsto|y|^{\frac{p-k}{p-1}}$ is a $p$-harmonic function from
$\mathbb{R}^{k}\backslash\{0\}\, $ to $\mathbb{R}\, ,$ the composite map
$|u(x)-y_{0}|^{\frac{p-k}{p-1}}:M\rightarrow\mathbb{R}$ is a $p$-harmonic
function with finite $q$-energy. By Theorem \ref{T:14}, in which $p\neq k$, we
obtain the conclusion.
\endproof

\subsection{Application to Conformal Maps\label{Maps}}

Our previous result can be applied to weakly conformal maps between equal
dimensional manifolds based on the following:\bigskip

\noindent\textbf{Theorem A} (\cite{OW}) $u:M\rightarrow N$ is an $m$-harmonic
morphism, if and only if $u$ is weakly conformal, where $m=\dim M=\dim N\,.$
\bigskip

\noindent For instance, stereographic projections $u:\mathbb{R}^{m}\rightarrow
S^{m}$ are $m$-harmonic maps and $m$-harmonic morphisms, for all $m \ge1\, .$

\begin{theorem}
Let $M^{m}$ be a complete noncompact $m$-manifold satisfying $\left(  P_{\rho
}\right)  $ and (\ref{Rs}), where $\tau<\frac{4\left(  q+b\right)  }{q^{2}}$
and $b=\min\{0,(m-2)(q-m)\}.$ If $u:M^{m}\rightarrow\mathbb{R}^{m}$ is a
weakly conformal map of finite $q$-energy $E_{q}\left(  u\right)  <\infty.$
Then $u$ is a constant if $m$ and $q$ satisfy one of the following:
\newline(1) $m=2$ and $q>0,$ \newline(2) $m=4,$ $q>1$ and $q+b>0,$ \newline(3)
$m>2,$ $m\neq4,$ and either $\frac{m\left(  m-2\right)  }{m-1}<q\leq
\min\left\{  2,m-\frac{\left(  m-4\right)  ^{2}m}{4\left(  m-2\right)
}\right\}  $ or\textrm{\ }$q>\max\{2,-b\}.$ \newline
\end{theorem}

%

\proof
By Theorem A (\cite{OW}), $u$ is an $m$-harmonic morphism. Now the result
follows immediately from Theorem \ref{T:13} in which $p=m.$ Since $\log|x|$ is
an $m$-harmonic function, $\log|u(x)-y_{0}|:M\rightarrow\mathbb{R}$ is an
$m$-harmonic function with finite $q$-energy. By Theorem \ref{T:14}, in which
$p=m$, we obtain the conclusion.
\endproof

\bigskip

\section{Appendix}

\subsection{The existence of the approximate solution}

In this subsection, we study an approximate solution $u_{\epsilon}$ of the $p
$-Laplace equation or a solution $u_{\epsilon}$ of a perturbed $p$-Laplace
equation
\begin{equation}%
\begin{array}
[c]{lllll}%
\Delta_{p,\epsilon}u_{\epsilon} & = & \text{\textrm{div}}\left(  \left(
\left\vert \nabla u_{\epsilon}\right\vert ^{2}+\epsilon\right)  ^{\frac
{p-2}{2}}\nabla u_{\epsilon}\right)  & = & 0
\end{array}
\label{a1.2}%
\end{equation}
on a domain $\Omega\subset M$ with boundary condition $u_{\epsilon}=u$ on
$\partial\Omega.$ That is, $u_{\epsilon}$ is the Euler-Lagrange equation of
the $\left(  p,\epsilon\right)  $-energy $E_{p,\epsilon}$ functional given by
\begin{equation}%
\begin{array}
[c]{l}%
E_{p,\epsilon}(\Psi)=\int_{\Omega}\left(  |\nabla\Psi|^{2}+\epsilon\right)
^{\frac{p}{2}}\,dv
\end{array}
\, \label{a1.3}%
\end{equation}
with $\Psi\in W^{1,p}\left(  \Omega\right)  ,$ and $\Psi=u$ on $\partial
\Omega.$

\begin{proposition}
[ The existence of $u_{\epsilon}$]\label{ex}Let $u$ be a $W^{1,p}$ function on
the closure $\bar{\Omega}$ of a domain $\Omega\subset M\, .$Then there is a
solution $u_{\epsilon}\in W^{1,p}\left(  \Omega\right)  $ of the
Euler-Lagrange equation of the $\left(  p,\epsilon\right)  $-energy
$E_{p,\epsilon}$ with $u_{\epsilon}=u$ on the boundary of $\Omega$ in the
trace sense.
\end{proposition}

\proof

Let $H$ be the set of functions $v\in W^{1,p}\left(  \Omega\right)  $ such
that $v=u$ on the boundary of $\Omega$ in the trace sense, and $I=\inf
\{E_{p,\epsilon}(v):v\in H\}.$ Then by assumption, $u\in H$, $H$ is nonempty,
and $I$ exists. Furthermore $I\leq E_{p,\epsilon}(u).$

Take a minimizing sequence $\{v_{i}\}_{i=1}^{\infty}$ such that $E_{p,\epsilon
}(v_{i})$ tends to $I$ as $i$ tends to $\infty$.

Then $\{v_{i}\}_{i=1}^{\infty}$ is a bounded sequence in $W^{1,p}\left(
\Omega\right)  $. Hence there exists a subsequence, say $\{u_{i}%
\}_{i=1}^{\infty}\,,$ converges weakly to $u_{\epsilon}$ in $W^{1,p}\left(
\Omega\right)  $, strongly in $L^{p}\left(  \Omega\right)  $, and pointwise
almost everywhere. We infer $u_{\epsilon}$ is in $H$ since $H$ is closed. Thus
$I\leq E_{p,\epsilon}\left(  u_{\epsilon}\right)  .$ \smallskip

To prove $I\geq E_{p,\epsilon}\left(  u_{\epsilon}\right)  \,,$ it suffices to
prove the lower semi-continuity of $E_{p,\epsilon}$ (two methods).

\textbf{Method 1}:

Since Banach-saks Theorem (see, e.g. \cite{Y} p. 120, \cite{RS} p. 80) asserts
there exists some subsequence, say it again $v_{i}$ for simplicity, such that
the average%
\[%
\begin{array}
[c]{l}%
w_{n}=\frac{v_{1}+v_{2}+\cdots+v_{n}}{n}%
\end{array}
\]
converges strongly to $u_{\epsilon}$ in $W^{1,p}\left(  \Omega\right)  $.
Combining this property and Lemma \ref{con}, we have $E_{p,\epsilon}\left(
w_{n}\right)  \rightarrow E_{p,\epsilon}\left(  u_{\epsilon}\right)  $ as
$n\rightarrow\infty.$

Moreover, according to the convexity of $E_{p,\epsilon},$ one has%
\[%
\begin{array}
[c]{l}%
E_{p,\epsilon}\left(  w_{n}\right)  \leq\frac{\sum_{i=1}^{n}E_{p,\epsilon
}\left(  v_{i}\right)  }{n}.
\end{array}
\]
This implies $E_{p,\epsilon}\left(  u_{\epsilon}\right)  \leq I$ as
$n\rightarrow\infty.$

So we obtain lower semi-continuity of $E_{p,\epsilon}.$

\textbf{Method 2}:

If $\dim M>2\,,$ we denote $T_{x} \Omega$ the tangent space to $\Omega\subset
M$ at $x$. Let $\nu_{i}(x) \in T_{x} \Omega$ be a unit vector perpendicular to
$\nabla u_{i}(x)\, , \nabla u_{\epsilon}(x) \in T_{x} \Omega\, ,$ for a.e.
$x\in\Omega\,.$ If $\dim M=2\, ,$ we isometrically embed $M$ into $N= M
\times\mathbb{R}$ with the standard product metric $\langle\, , \, \rangle
_{N}$ and choose $\nu_{i}(x)\, $ to be a unit vector in $\mathbb{R}\,.$

In either case, we set $b(x)=\nabla u_{i}(x)+\sqrt{\epsilon}{\nu}_{i}(x)$ and
$a(x)=\nabla u_{\epsilon}(x)+\sqrt{\epsilon}\nu_{i}(x)$. Then on $\Omega\,,$
$|b|=\sqrt{|\nabla u_{i}|^{2}+\epsilon}$ and $|a|=\sqrt{|\nabla u_{\epsilon
}|^{2}+\epsilon}.$

If $m=2$, applying the inequality
\[
|b|^{p}\geq|a|^{p}+p\langle|a|^{p-2}a,b-a\rangle_{N}%
\]
and integrating it over $\Omega$, we have via $\nu_{i}(x)\bot\nabla
u_{\epsilon}\, ,$ and $\nu_{i}(x)\bot\nabla u_{i}\, ,$ for a.e. $x\in\Omega\,
,$
\[%
\begin{array}
[c]{lll}%
E_{p,\epsilon}(u_{i}) & \geq & E_{p,\epsilon}(u_{\epsilon})+\int_{\Omega
}\langle(|\nabla u_{\epsilon}|^{2}+\epsilon)^{\frac{p-2}{2}}(\nabla
u_{\epsilon}+\sqrt{\epsilon}\nu_{i}),\nabla u_{i}-\nabla u_{\epsilon}%
\rangle_{N}dv\\
& = & E_{p,\epsilon}(u_{\epsilon})+\int_{\Omega}\langle(|\nabla u_{\epsilon
}|^{2}+\epsilon)^{\frac{p-2}{2}}\nabla u_{\epsilon},\nabla u_{i}-\nabla
u_{\epsilon}\rangle_{M}dv
\end{array}
\]
We note that in the last term, $(|\nabla u_{\epsilon}|^{2}+\epsilon
)^{\frac{p-2}{2}}\nabla u_{\epsilon}$ is in $L^{\frac{p}{p-1}}(\Omega
)\,,\nabla u_{i}-\nabla u_{\epsilon}$ is in $L^{p}(\Omega)\,.$ Thus,
$\langle(|\nabla u_{\epsilon}|^{2}+\epsilon)^{\frac{p-2}{2}}\nabla
u_{\epsilon},\nabla u_{i}-\nabla u_{\epsilon}\rangle_{M}$ is in $L^{1}%
(\Omega)\, .$ Since $\nabla u_{i}$ converges weakly to $\nabla u_{\epsilon}$
in $L^{p}$, the last term tends to $0$ as $i$ tends to $\infty$. It follows
that $E_{p,\epsilon}(u_{\epsilon})\leq\lim\inf_{i\rightarrow\infty
}E_{p,\epsilon}(u_{i})=I$. Similarly, if $\dim M > 2\,,$ we obtain directly
\[
E_{p,\epsilon}(u_{i}) \geq E_{p,\epsilon}(u_{\epsilon})+\int_{\Omega}%
\langle(|\nabla u_{\epsilon}|^{2}+\epsilon)^{\frac{p-2}{2}}\nabla u_{\epsilon
},\nabla u_{i}-\nabla u_{\epsilon}\rangle_{M}dv.
\]
Proceed in the same way, the assertion follows. \smallskip

\endproof

\bigskip

\bigskip

\begin{lemma}
\label{con}If $v_{i}$ converges strongly to $v_{0}$ in $W^{1,p},$ then
$E_{p,\epsilon}\left(  v_{i}\right)  $ converges to $E_{p,\epsilon}\left(
v_{0}\right)  .$
\end{lemma}

\proof Step 1: Since $v_{i}$ converges strongly to $v_{0}$ in $W^{1,p}$ i.e.
$\int_{\Omega}\left\vert \nabla v_{i}-\nabla v_{0}\right\vert ^{p}%
dv\rightarrow0$ as $i\rightarrow\infty.$ Then
\[%
\begin{array}
[c]{l}%
\int_{\left\vert \nabla v_{i}\right\vert \geq\left\vert \nabla v_{0}%
\right\vert }\left\vert \nabla v_{i}-\nabla v_{0}\right\vert ^{p}%
dv\rightarrow0\text{ and}\int_{\left\vert \nabla v_{i}\right\vert <\left\vert
\nabla v_{0}\right\vert }\left\vert \nabla v_{i}-\nabla v_{0}\right\vert
^{p}dv\rightarrow0
\end{array}
\]
as $i\rightarrow\infty.$ By using Minkowski's inequality, these also imply%
\[%
\begin{array}
[c]{l}%
\int_{\left\vert \nabla v_{i}\right\vert \geq\left\vert \nabla v_{0}%
\right\vert }\left(  \left\vert \nabla v_{i}\right\vert ^{p}-\left\vert \nabla
v_{0}\right\vert ^{p}\right)  dv\rightarrow0\text{ and }\int_{\left\vert
\nabla v_{i}\right\vert <\left\vert \nabla v_{0}\right\vert }\left(
\left\vert \nabla v_{0}\right\vert ^{p}-\left\vert \nabla v_{i}\right\vert
^{p}\right)  dv\rightarrow0
\end{array}
\]
as $i\rightarrow\infty.$ That is, $\int_{\Omega}\left\vert \left\vert \nabla
v_{0}\right\vert ^{p}-\left\vert \nabla v_{i}\right\vert ^{p}\right\vert
dv\rightarrow0$ as $i\rightarrow\infty.$

Step 2: If we show that, for any positive constant $\delta>0,$
\begin{equation}%
\begin{array}
[c]{l}%
\left\vert \left(  \left\vert \nabla v_{i}\right\vert ^{2}+\epsilon\right)
^{\frac{p}{2}}-\left(  \left\vert \nabla v_{0}\right\vert ^{2}+\epsilon
\right)  ^{\frac{p}{2}}\right\vert \leq a\left\vert \left\vert \nabla
v_{i}\right\vert ^{p}-\left\vert \nabla v_{0}\right\vert ^{p}\right\vert
+\delta
\end{array}
\label{e0-1}%
\end{equation}
where $a$ is a positive constant independent of $i,$ $v_{i}$ and $v_{0}.$ Then
we have, by step 1,%
\[%
\begin{array}
[c]{lll}%
\left\vert E_{p,\epsilon}\left(  v_{i}\right)  -E_{p,\epsilon}\left(
v_{0}\right)  \right\vert  & \leq & \int_{\Omega}\left\vert \left(  \left\vert
\nabla v_{i}\right\vert ^{2}+\epsilon\right)  ^{\frac{p}{2}}-\left(
\left\vert \nabla v_{0}\right\vert ^{2}+\epsilon\right)  ^{\frac{p}{2}%
}\right\vert dv\\
& \leq & a\int_{\Omega}\left\vert \left\vert \nabla v_{i}\right\vert
^{p}-\left\vert \nabla v_{0}\right\vert ^{p}\right\vert dv+\delta\left\vert
\Omega\right\vert \\
& \rightarrow & \delta\left\vert \Omega\right\vert \text{ as }i\rightarrow
\infty.
\end{array}
\]
This implies $E_{p,\epsilon}\left(  v_{i}\right)  \rightarrow E_{p,\epsilon
}\left(  v_{0}\right)  .$

To show (\ref{e0-1}), we only claim that, $X,Y\in%
\mathbb{R}
^{n}$ with $\left\vert X\right\vert \geq\left\vert Y\right\vert ,$
\begin{equation}%
\begin{array}
[c]{l}%
\left(  \left\vert X\right\vert ^{2}+\epsilon\right)  ^{\frac{p}{2}}-\left(
\left\vert Y\right\vert ^{2}+\epsilon\right)  ^{\frac{p}{2}}\leq a\left(
\left\vert X\right\vert ^{p}-\left\vert Y\right\vert ^{p}\right)  +\delta.
\end{array}
\label{e00}%
\end{equation}

Let $f(t)=\left(  \left\vert X\right\vert ^{2}+t\right)  ^{\frac{p}{2}%
}-\left(  \left\vert Y\right\vert ^{2}+t\right)  ^{\frac{p}{2}},$ $t\geq0.$
Then we have $f(0)=\left\vert X\right\vert ^{p}-\left\vert Y\right\vert ^{p}$
and $f(\epsilon)=\left(  \left\vert X\right\vert ^{2}+\epsilon\right)
^{\frac{p}{2}}-\left(  \left\vert Y\right\vert ^{2}+\epsilon\right)
^{\frac{p}{2}}.$

Since
\[%
\begin{array}
[c]{l}%
f^{\prime}(t)=\frac{p}{2}\left(  \left(  \left\vert X\right\vert
^{2}+t\right)  ^{\frac{p-2}{2}}-\left(  \left\vert Y\right\vert ^{2}+t\right)
^{\frac{p-2}{2}}\right)  ,
\end{array}
\]
then $f(t)$ is a decreasing function for $1\leq p\leq2.$ Hence we have
$f(\epsilon)\leq f(0)$ whenever $1\leq p\leq2.$

If $2<p\leq4,$ then, for $s>0,$%
\begin{equation}%
\begin{array}
[c]{lll}%
f(s)-f(0) & = & \int_{0}^{s}f^{\prime}(t)dt\\
& = & \frac{p}{2}\int_{0}^{s}\left(  \left\vert X\right\vert ^{2}+t\right)
^{\frac{p-2}{2}}-\left(  \left\vert Y\right\vert ^{2}+t\right)  ^{\frac
{p-2}{2}}dt,\\
& \leq & \frac{ps}{2}\left(  \left\vert X\right\vert ^{p-2}-\left\vert
Y\right\vert ^{p-2}\right)  ,
\end{array}
\label{e1}%
\end{equation}
since $1<p-2\leq2.$

For any $\delta_{1}>0,$
\[%
\begin{array}
[c]{l}%
\left\vert X\right\vert ^{p-2}-\left\vert Y\right\vert ^{p-2}\leq\left\{
\begin{array}
[c]{ll}%
\left\vert X\right\vert ^{p-2} & \text{if }\left\vert X\right\vert +\left\vert
Y\right\vert <\delta_{1},\\
\frac{\left(  \left\vert X\right\vert +\left\vert Y\right\vert \right)
^{2}\left(  \left\vert X\right\vert ^{p-2}-\left\vert Y\right\vert
^{p-2}\right)  }{\delta_{1}^{2}} & \text{if }\left\vert X\right\vert
+\left\vert Y\right\vert \geq\delta_{1}.
\end{array}
\right.
\end{array}
\]
Since%
\[
\left\{
\begin{array}
[c]{ll}%
\left\vert X\right\vert ^{p-2}\leq\delta_{1}^{p-2} & \text{if }\left\vert
X\right\vert +\left\vert Y\right\vert <\delta_{1},\\
\frac{\left(  \left\vert X\right\vert +\left\vert Y\right\vert \right)
^{2}\left(  \left\vert X\right\vert ^{p-2}-\left\vert Y\right\vert
^{p-2}\right)  }{\delta_{1}^{2}}\leq\frac{2}{\delta_{1}^{2}}\left(  \left\vert
X\right\vert ^{p}-\left\vert Y\right\vert ^{p}\right)  & \text{if }\left\vert
X\right\vert +\left\vert Y\right\vert \geq\delta_{1}.
\end{array}
\right.
\]
So we have
\begin{equation}%
\begin{array}
[c]{l}%
\left\vert X\right\vert ^{p-2}-\left\vert Y\right\vert ^{p-2}\leq\frac
{2}{\delta_{1}^{2}}\left(  \left\vert X\right\vert ^{p}-\left\vert
Y\right\vert ^{p}+\delta_{1}^{p}\right)  ,
\end{array}
\label{e3}%
\end{equation}
and then (\ref{e1}) can be rewritten as%
\begin{equation}%
\begin{array}
[c]{l}%
\left(  \left\vert X\right\vert ^{2}+s\right)  ^{\frac{p}{2}}-\left(
\left\vert Y\right\vert ^{2}+s\right)  ^{\frac{p}{2}}\leq\left(  1+\frac
{ps}{\delta_{1}^{2}}\right)  \left(  \left\vert X\right\vert ^{p}-\left\vert
Y\right\vert ^{p}\right)  +\left(  \frac{ps}{\delta_{1}^{2}}\right)
\delta_{1}^{p}.
\end{array}
\label{e3-1}%
\end{equation}
Hence we have%
\[%
\begin{array}
[c]{l}%
\left(  \left\vert X\right\vert ^{2}+\epsilon\right)  ^{\frac{p}{2}}-\left(
\left\vert Y\right\vert ^{2}+\epsilon\right)  ^{\frac{p}{2}}\leq a\left(
\left\vert X\right\vert ^{p}-\left\vert Y\right\vert ^{p}\right)  +\delta,
\end{array}
\]
where $a=1+\frac{p\epsilon}{\delta_{1}^{2}}$ and $\delta=\left(
\frac{p\epsilon}{\delta_{1}^{2}}\right)  \delta_{1}^{p}.$

If $4<p\leq6,$ then one has $2<p-2\leq4,$ so (\ref{e3}) and (\ref{e3-1}) imply%
\[%
\begin{array}
[c]{lll}%
f(s)-f(0) & = & \frac{p}{2}\int_{0}^{s}\left(  \left\vert X\right\vert
^{2}+t\right)  ^{\frac{p-2}{2}}-\left(  \left\vert Y\right\vert ^{2}+t\right)
^{\frac{p-2}{2}}dt\\
& \leq & \frac{p}{2}\int_{0}^{s}\left(  1+\frac{pt}{\delta_{1}^{2}}\right)
\left(  \left\vert X\right\vert ^{p-2}-\left\vert Y\right\vert ^{p-2}\right)
+\left(  \frac{pt}{\delta_{1}^{2}}\right)  \delta_{1}^{p-2}dt\\
& \leq & \frac{p}{2}\left(  s+\frac{ps^{2}}{2\delta_{1}^{2}}\right)  \left(
\frac{2}{\delta_{1}^{2}}\left(  \left\vert X\right\vert ^{p}-\left\vert
Y\right\vert ^{p}+\delta_{1}^{p}\right)  \right)  +\frac{p}{2}\left(
\frac{ps^{2}}{2\delta_{1}^{2}}\right)  \delta_{1}^{p-2}\\
& \leq & \left(  \frac{ps}{\delta_{1}^{2}}+\frac{1}{2}\left(  \frac{ps}%
{\delta_{1}^{2}}\right)  ^{2}\right)  \left(  \left\vert X\right\vert
^{p}-\left\vert Y\right\vert ^{p}\right)  +\left(  \frac{ps}{\delta_{1}^{2}%
}+\left(  \frac{ps}{\delta_{1}^{2}}\right)  ^{2}\right)  \delta_{1}^{p}.
\end{array}
\]
Hence
\[%
\begin{array}
[c]{lll}%
\left(  \left\vert X\right\vert ^{2}+s\right)  ^{\frac{p}{2}}-\left(
\left\vert Y\right\vert ^{2}+s\right)  ^{\frac{p}{2}} & \leq & \left(
1+\frac{ps}{\delta_{1}^{2}}+\frac{1}{2}\left(  \frac{ps}{\delta_{1}^{2}%
}\right)  ^{2}\right)  \left(  \left\vert X\right\vert ^{p}-\left\vert
Y\right\vert ^{p}\right) \\
&  & +\left(  \frac{ps}{\delta_{1}^{2}}+\left(  \frac{ps}{\delta_{1}^{2}%
}\right)  ^{2}\right)  \delta_{1}^{p}.
\end{array}
\]

In particular, we obtain%
\[%
\begin{array}
[c]{lll}%
\left(  \left\vert X\right\vert ^{2}+\epsilon\right)  ^{\frac{p}{2}}-\left(
\left\vert Y\right\vert ^{2}+\epsilon\right)  ^{\frac{p}{2}} & \leq & a\left(
\left\vert X\right\vert ^{p}-\left\vert Y\right\vert ^{p}\right)  +\delta,
\end{array}
\]
where $a=1+\frac{p\epsilon}{\delta_{1}^{2}}+\frac{1}{2}\left(  \frac
{p\epsilon}{\delta_{1}^{2}}\right)  ^{2}$ and $\delta=\left(  \frac{p\epsilon
}{\delta_{1}^{2}}+\left(  \frac{p\epsilon}{\delta_{1}^{2}}\right)
^{2}\right)  \delta_{1}^{p}.$

By mathematical induction, we conclude that, for any $p>2$ satisfying
$2q<p\leq2q+2,$ $q\in%
\mathbb{Z}
^{+},$%
\[%
\begin{array}
[c]{lll}%
\left(  \left\vert X\right\vert ^{2}+\epsilon\right)  ^{\frac{p}{2}}-\left(
\left\vert Y\right\vert ^{2}+\epsilon\right)  ^{\frac{p}{2}} & \leq & \left(
1+\sum_{n=1}^{q}\frac{1}{n!}\left(  \frac{p\epsilon}{\delta_{1}^{2}}\right)
^{n}\right)  \left(  \left\vert X\right\vert ^{p}-\left\vert Y\right\vert
^{p}\right) \\
&  & +\left(  \sum_{n=1}^{q}\left(  \frac{p\epsilon}{\delta_{1}^{2}}\right)
^{n}\right)  \delta_{1}^{p}.
\end{array}
\]
If we select $\delta_{1}$ small enough such that $\left(  \sum_{n=1}%
^{q}\left(  \frac{p\epsilon}{2\delta_{1}^{2}}\right)  ^{n}\right)  \delta
_{1}^{p}=\delta,$ then we have (\ref{e00}) with $a=\left(  1+\sum_{n=1}%
^{q}\left(  \frac{p\epsilon}{2\delta_{1}^{2}}\right)  ^{n}\right)  . $

\endproof

\bigskip

\bigskip

\bigskip

\bigskip

\bigskip

\subsection{$\epsilon$-regularization of $p$-Laplacian}

\begin{proposition}
\label{3.1}Let $u$ be a weak solution of the $p$-Laplace equation (\ref{1.0}).
For every $\epsilon>0,$ let $u_{\epsilon}$ be a solution of the Euler-Lagrange
equation (\ref{a1.2}) with $u-u_{\epsilon}\in W_{0}^{1,p}\left(
\Omega\right)  ,$ where $\Omega$ is a domain in $M\,.$ Then $u_{\epsilon}\in
C_{loc}^{\infty}\left(  \Omega\right)  $ is a strong solution of (\ref{a1.2}),
and $u_{\epsilon}$ converges strongly to $u$ in $W^{1,p}\left(  \Omega\right)
$ as $\epsilon\rightarrow0\,.$
\end{proposition}

\bigskip

\proof

Such solution $u_{\epsilon}$ exists (Proposition \ref{ex}), and $u_{\epsilon
}\in C_{loc}^{\infty}\left(  \Omega\right)  $ by the usual arguments of
boot-strap (see, e.g. \cite{LU} Chapter 4, \cite{S} Theorem 3.3, \cite{HS}
Theorem 14.2, \cite{Ho}). That is, $u_{\epsilon}$ is the strong solution of
the partial differential equation (\ref{1.2}).

Since $u_{\epsilon}$ and $u$ are the minimizers of the energy functions
\[%
\begin{array}
[c]{c}%
\int_{\Omega}\left\vert \left\vert \nabla\phi\right\vert ^{2}+\epsilon
\right\vert ^{p/2}dv\text{ and }\int_{\Omega}\left\vert \nabla\phi\right\vert
^{p}dv,
\end{array}
\]
respectively, over all functions $\phi\in W^{1,p}\left(  \Omega\right)  $ and
$\phi=u$ on $\partial\Omega.$ Then one has
\begin{equation}%
\begin{array}
[c]{c}%
\int_{\Omega}\left\vert \nabla u\right\vert ^{p}dv\leq\int_{\Omega}\left\vert
\nabla u_{\epsilon}\right\vert ^{p}dv
\end{array}
\label{p1}%
\end{equation}
and
\begin{equation}%
\begin{array}
[c]{c}%
\int_{\Omega}\left\vert \left\vert \nabla u_{\epsilon}\right\vert
^{2}+\epsilon\right\vert ^{p/2}dv\text{ }\leq\int_{\Omega}\left\vert
\left\vert \nabla u\right\vert ^{2}+\epsilon\right\vert ^{p/2}dv.
\end{array}
\label{p2}%
\end{equation}

Combining (\ref{p1}) and (\ref{p2}),
\[%
\begin{array}
[c]{c}%
\int_{\Omega}\left\vert \nabla u\right\vert ^{p}dv\leq\int_{\Omega}\left\vert
\nabla u_{\epsilon}\right\vert ^{p}dv\leq\int_{\Omega}\left\vert \left\vert
\nabla u_{\epsilon}\right\vert ^{2}+\epsilon\right\vert ^{p/2}dv\text{ }%
\leq\int_{\Omega}\left\vert \left\vert \nabla u\right\vert ^{2}+\epsilon
\right\vert ^{p/2}dv,
\end{array}
\]
one has $\left\Vert \nabla u_{\epsilon}\right\Vert _{p}\rightarrow\left\Vert
\nabla u\right\Vert _{p}$ as $\epsilon\rightarrow0.$ Moreover, by Lemma
\ref{c1} $\nabla u_{\epsilon}\rightarrow\nabla u$ a.e. on $\Omega$ for $p>1,$
we have $\nabla u_{\epsilon}\rightarrow\nabla u$ in $L^{p}\left(
\Omega\right)  ,$ and then $p$-Poincar\'{e} inequality implies $u_{\epsilon
}\rightarrow u$ in $W^{1,p}\left(  \Omega\right)  .$

\endproof

\begin{lemma}
\label{c1}$\nabla u_{\epsilon}\rightarrow\nabla u$ a.e. on $\Omega$ for $p>1.
$
\end{lemma}

\bigskip

First, we recall the following inequality (cf. \cite{Li} Chapter 10, or
\cite{HPV} Lemma 4)

\begin{proposition}
\label{12}Let $X$ and $Y$ be vector fields on $\Omega.$ Then%
\begin{equation}%
\begin{array}
[c]{l}%
\left\langle X-Y,\left\vert X\right\vert ^{p-2}X-\left\vert Y\right\vert
^{p-2}Y\right\rangle \geq C\Psi\left(  X,Y\right)  ,
\end{array}
\label{c0}%
\end{equation}
where%
\begin{equation}%
\begin{array}
[c]{l}%
\Psi\left(  X,Y\right)  =\left\{
\begin{array}
[c]{ll}%
\left\vert X-Y\right\vert ^{p} & \text{if }p\geq2,\\
\frac{\left(  p-1\right)  \left\vert X-Y\right\vert ^{2}}{\left(  1+\left\vert
X\right\vert ^{2}+\left\vert Y\right\vert ^{2}\right)  ^{\frac{2-p}{2}}} &
\text{if }1<p<2.
\end{array}
\right.
\end{array}
\label{c0-1}%
\end{equation}

\end{proposition}

\bigskip

\proof

Since $u-u_{\epsilon}\in W_{0}^{1,p}\left(  \Omega\right)  ,$ one has
\[%
\begin{array}
[c]{l}%
\int_{\Omega}\left\vert \nabla u\right\vert ^{p-2}\left\langle \nabla
u,\nabla\left(  u-u_{\epsilon}\right)  \right\rangle dv=0
\end{array}
\]
and
\[%
\begin{array}
[c]{l}%
\int_{\Omega}\left\vert \left\vert \nabla u_{\epsilon}\right\vert
^{2}+\epsilon\right\vert ^{\frac{p-2}{2}}\left\langle \nabla u_{\epsilon
},\nabla\left(  u-u_{\epsilon}\right)  \right\rangle dv=0.
\end{array}
\]
Then%
\[%
\begin{array}
[c]{lll}%
0 & = & \int_{\Omega}\left\vert \nabla u\right\vert ^{p-2}\left\langle \nabla
u,\nabla\left(  u-u_{\epsilon}\right)  \right\rangle -\left\vert \left\vert
\nabla u_{\epsilon}\right\vert ^{2}+\epsilon\right\vert ^{\frac{p-2}{2}%
}\left\langle \nabla u_{\epsilon},\nabla\left(  u-u_{\epsilon}\right)
\right\rangle dv\\
& = & \int_{\Omega}\left\vert \nabla u\right\vert ^{p}-\left\vert \nabla
u\right\vert ^{p-2}\left\langle \nabla u,\nabla u_{\epsilon}\right\rangle \\
&  & -\left\vert \left\vert \nabla u_{\epsilon}\right\vert ^{2}+\epsilon
\right\vert ^{\frac{p-2}{2}}\left\langle \nabla u_{\epsilon},\nabla
u\right\rangle +\left\vert \left\vert \nabla u_{\epsilon}\right\vert
^{2}+\epsilon\right\vert ^{\frac{p-2}{2}}\left\vert \nabla u_{\epsilon
}\right\vert ^{2}dv.
\end{array}
\]
This equality can be rewritten as $LHS1=RHS,$ where%
\[%
\begin{array}
[c]{l}%
LHS1=\int_{\Omega}\left\vert \nabla u\right\vert ^{p}-\left\vert \nabla
u\right\vert ^{p-2}\left\langle \nabla u,\nabla u_{\epsilon}\right\rangle
-\left\vert \nabla u_{\epsilon}\right\vert ^{p-2}\left\langle \nabla
u_{\epsilon},\nabla u\right\rangle +\left\vert \left\vert \nabla u_{\epsilon
}\right\vert ^{2}+\epsilon\right\vert ^{\frac{p}{2}}dv
\end{array}
\]
and%
\[%
\begin{array}
[c]{l}%
RHS=\int_{\Omega}\left(  \left\vert \left\vert \nabla u_{\epsilon}\right\vert
^{2}+\epsilon\right\vert ^{\frac{p-2}{2}}-\left\vert \nabla u_{\epsilon
}\right\vert ^{p-2}\right)  \left\langle \nabla u_{\epsilon},\nabla
u\right\rangle +\epsilon\left\vert \left\vert \nabla u_{\epsilon}\right\vert
^{2}+\epsilon\right\vert ^{\frac{p-2}{2}}dv.
\end{array}
\]

It is easy to see that $LHS1\geq LHS2$ where%
\[%
\begin{array}
[c]{l}%
LHS2=\int_{\Omega}\left\vert \nabla u\right\vert ^{p}-\left\vert \nabla
u\right\vert ^{p-2}\left\langle \nabla u,\nabla u_{\epsilon}\right\rangle
-\left\vert \nabla u_{\epsilon}\right\vert ^{p-2}\left\langle \nabla
u_{\epsilon},\nabla u\right\rangle +\left\vert \nabla u_{\epsilon}\right\vert
^{p}dv.
\end{array}
\]

So, we select $X=$ $\nabla u$ and $Y=\nabla u_{\epsilon}$, then Proposition
\ref{12} implies%
\[%
\begin{array}
[c]{l}%
LHS2\geq C\int_{\Omega}\Psi\left(  \nabla u,\nabla u_{\epsilon}\right)
dv\geq0
\end{array}
\]
where%
\[%
\begin{array}
[c]{l}%
\Psi\left(  \nabla u,\nabla u_{\epsilon}\right)  =\left\{
\begin{array}
[c]{ll}%
\left\vert \nabla u-\nabla u_{\epsilon}\right\vert ^{p} & \text{if }p\geq2,\\
\frac{\left(  p-1\right)  \left\vert \nabla u-\nabla u_{\epsilon}\right\vert
^{2}}{\left(  1+\left\vert \nabla u\right\vert ^{2}+\left\vert \nabla
u_{\epsilon}\right\vert ^{2}\right)  ^{\frac{2-p}{2}}} & \text{if }1<p<2.
\end{array}
\right.
\end{array}
\]

If we can show that%
\[%
\begin{array}
[c]{l}%
RHS\rightarrow0\text{ as }\epsilon\rightarrow0,
\end{array}
\]

Then we have%
\[%
\begin{array}
[c]{l}%
\int_{\Omega}\Psi\left(  \nabla u,\nabla u_{\epsilon}\right)  dv\rightarrow
0\text{ as }\epsilon\rightarrow0.
\end{array}
\]
Therefore $\nabla u_{\epsilon}\rightarrow\nabla u$ a.e. on $\Omega.$

Now we claim that
\[%
\begin{array}
[c]{l}%
RHS=RHS1+RHS2\rightarrow0
\end{array}
\]
as $\epsilon\rightarrow0,$ where%
\[%
\begin{array}
[c]{l}%
RHS1=\int_{\Omega}\epsilon\left\vert \left\vert \nabla u_{\epsilon}\right\vert
^{2}+\epsilon\right\vert ^{\frac{p-2}{2}}dv
\end{array}
\]
and%
\[%
\begin{array}
[c]{l}%
RHS2=\int_{\Omega}\left(  \left\vert \left\vert \nabla u_{\epsilon}\right\vert
^{2}+\epsilon\right\vert ^{\frac{p-2}{2}}-\left\vert \nabla u_{\epsilon
}\right\vert ^{p-2}\right)  \left\langle \nabla u_{\epsilon},\nabla
u\right\rangle dv.
\end{array}
\]

It is easy to see that, if $\left\vert \nabla u_{\epsilon}\right\vert ^{2}%
\geq1,$%
\[%
\begin{array}
[c]{lll}%
\int_{\Omega}\epsilon\left\vert \left\vert \nabla u_{\epsilon}\right\vert
^{2}+\epsilon\right\vert ^{\frac{p-2}{2}}dv & \leq & \int_{\Omega}%
\epsilon\left\vert \nabla u_{\epsilon}\right\vert ^{2}\left\vert \left\vert
\nabla u_{\epsilon}\right\vert ^{2}+\epsilon\right\vert ^{\frac{p-2}{2}}dv\\
& \leq & \epsilon\int_{\Omega}\left\vert \left\vert \nabla u_{\epsilon
}\right\vert ^{2}+\epsilon\right\vert ^{\frac{p}{2}}dv,
\end{array}
\]
and if $\left\vert \nabla u_{\epsilon}\right\vert ^{2}<1,$%
\[%
\begin{array}
[c]{l}%
\int_{\Omega}\epsilon\left\vert \left\vert \nabla u_{\epsilon}\right\vert
^{2}+\epsilon\right\vert ^{\frac{p-2}{2}}dv\leq\left\{
\begin{array}
[c]{ll}%
\epsilon\left(  1+\epsilon\right)  ^{\frac{p-2}{2}}\cdot vol\left(
\Omega\right)  & \text{if }p\geq2\\
\epsilon^{\frac{p}{2}}\cdot vol\left(  \Omega\right)  & \text{if }p<2.
\end{array}
\right.
\end{array}
\]

So we have $RHS1\rightarrow0$ as $\epsilon\rightarrow0.$

Now we focus on the term $RHS2$,%
\[%
\begin{array}
[c]{lll}%
RHS2 & = & \int_{\Omega}\left(  \left\vert \left\vert \nabla u_{\epsilon
}\right\vert ^{2}+\epsilon\right\vert ^{\frac{p-2}{2}}-\left\vert \nabla
u_{\epsilon}\right\vert ^{p-2}\right)  \left\langle \nabla u_{\epsilon},\nabla
u\right\rangle dv\\
& \leq & \int_{\Omega}\left\vert \left\vert \left\vert \nabla u_{\epsilon
}\right\vert ^{2}+\epsilon\right\vert ^{\frac{p-2}{2}}-\left\vert \nabla
u_{\epsilon}\right\vert ^{p-2}\right\vert \left\vert \nabla u_{\epsilon
}\right\vert \left\vert \nabla u\right\vert dv.
\end{array}
\]

In the case $p\geq2,$ one may rewrite it as%
\[%
\begin{array}
[c]{lll}%
RHS2 & \leq & \int_{\Omega}\left(  \left\vert \left\vert \nabla u_{\epsilon
}\right\vert ^{2}+\epsilon\right\vert ^{\frac{p-2}{2}}-\left\vert \nabla
u_{\epsilon}\right\vert ^{p-2}\right)  \left\vert \nabla u_{\epsilon
}\right\vert \left\vert \nabla u\right\vert dv\\
& \leq & \int_{\Omega}\left(  \left\vert \left\vert \nabla u_{\epsilon
}\right\vert ^{2}+\epsilon\right\vert ^{\frac{p-1}{2}}-\left\vert \nabla
u_{\epsilon}\right\vert ^{p-1}\right)  \left\vert \nabla u\right\vert dv.
\end{array}
\]

If $p\geq3,$ using mean value theorem, we have the inequality%
\[%
\begin{array}
[c]{l}%
\left(  x+\epsilon\right)  ^{q}-x^{q}=q\epsilon\left(  x+\epsilon_{1}\right)
^{q-1}\leq q\epsilon\left(  x+\epsilon\right)  ^{q-1}%
\end{array}
\]
here $q=\frac{p-1}{2}\geq1,$ $x\geq0$ and $\epsilon_{1}\in\left(
0,\epsilon\right)  .$ Hence%
\[%
\begin{array}
[c]{lll}%
RHS2 & \leq & \frac{\left(  p-1\right)  \epsilon}{2}\int_{\Omega}\left\vert
\left\vert \nabla u_{\epsilon}\right\vert ^{2}+\epsilon\right\vert
^{\frac{p-3}{2}}\left\vert \nabla u\right\vert dv\\
& \leq & \left\{
\begin{array}
[c]{ll}%
\frac{\left(  p-1\right)  \epsilon}{2}\int_{\Omega}\left\vert \left\vert
\nabla u_{\epsilon}\right\vert ^{2}+\epsilon\right\vert ^{\frac{p-1}{2}%
}\left\vert \nabla u\right\vert dv\text{ \ \ } & \text{if }\left\vert \nabla
u_{\epsilon}\right\vert ^{2}>1\\
\frac{\left(  p-1\right)  \epsilon}{2}\left(  1+\epsilon\right)  ^{\frac
{p-3}{2}}\int_{\Omega}\left\vert \nabla u\right\vert dv & \text{if }\left\vert
\nabla u_{\epsilon}\right\vert ^{2}\leq1
\end{array}
\right. \\
& \leq & \left\{
\begin{array}
[c]{ll}%
\frac{\left(  p-1\right)  \epsilon}{2}\left(  \int_{\Omega}\left\vert
\left\vert \nabla u_{\epsilon}\right\vert ^{2}+\epsilon\right\vert ^{\frac
{p}{2}}\right)  ^{\frac{p-1}{p}}\left(  \int_{\Omega}\left\vert \nabla
u\right\vert ^{p}\right)  ^{\frac{1}{p}}dv & \text{if }\left\vert \nabla
u_{\epsilon}\right\vert ^{2}>1\\
\frac{\left(  p-1\right)  \epsilon}{2}\left(  1+\epsilon\right)  ^{\frac
{p-3}{2}}\left(  vol\left(  \Omega\right)  \right)  ^{\frac{p-1}{p}}\left(
\int_{\Omega}\left\vert \nabla u\right\vert ^{p}\right)  ^{\frac{1}{p}%
}dv\text{ \ } & \text{if }\left\vert \nabla u_{\epsilon}\right\vert ^{2}\leq1
\end{array}
\right. \\
& \rightarrow & 0\text{ as }\epsilon\rightarrow0.
\end{array}
\]

If $2\leq p\leq3,$ using the inequality%
\[%
\begin{array}
[c]{l}%
\left(  x+\epsilon\right)  ^{q}-x^{q}\leq\epsilon^{q}%
\end{array}
\]
here $\frac{1}{2}\leq q=\frac{p-1}{2}\leq1,$ $x\geq0,$ then%
\[%
\begin{array}
[c]{lll}%
RHS2 & \leq & \int_{\Omega}\left(  \left\vert \left\vert \nabla u_{\epsilon
}\right\vert ^{2}+\epsilon\right\vert ^{\frac{p-1}{2}}-\left\vert \nabla
u_{\epsilon}\right\vert ^{p-1}\right)  \left\vert \nabla u\right\vert \\
& \leq & \epsilon^{\frac{p-1}{2}}\int_{\Omega}\left\vert \nabla u\right\vert
\\
& \leq & \epsilon^{\frac{p-1}{2}}\left(  vol\left(  \Omega\right)  \right)
^{\frac{p-1}{p}}\left(  \int_{\Omega}\left\vert \nabla u\right\vert
^{p}\right)  ^{\frac{1}{p}}\\
& \rightarrow & 0\text{ as }\epsilon\rightarrow0.
\end{array}
\]

In the case $1<p<2,$ one may rewrite $RHS2$ as%
\[%
\begin{array}
[c]{l}%
RHS2\leq\int_{\Omega}\left(  \left\vert \nabla u_{\epsilon}\right\vert
^{p-2}-\left\vert \left\vert \nabla u_{\epsilon}\right\vert ^{2}%
+\epsilon\right\vert ^{\frac{p-2}{2}}\right)  \left\vert \nabla u_{\epsilon
}\right\vert \left\vert \nabla u\right\vert .
\end{array}
\]
Since $0<\frac{2-p}{2}<1,$ then we have%
\[%
\begin{array}
[c]{lll}%
RHS2 & = & \int_{\Omega}\frac{\left\vert \left\vert \nabla u_{\epsilon
}\right\vert ^{2}+\epsilon\right\vert ^{\frac{2-p}{2}}-\left\vert \nabla
u_{\epsilon}\right\vert ^{2-p}}{\left\vert \left\vert \nabla u_{\epsilon
}\right\vert ^{2}+\epsilon\right\vert ^{\frac{2-p}{2}}}\left\vert \nabla
u_{\epsilon}\right\vert ^{p-1}\left\vert \nabla u\right\vert dv\\
& \leq & \int_{\Omega}\frac{\epsilon^{\frac{2-p}{2}}}{\left\vert \left\vert
\nabla u_{\epsilon}\right\vert ^{2}+\epsilon\right\vert ^{\frac{3-2p}{2}}%
}\cdot\frac{\left\vert \nabla u_{\epsilon}\right\vert ^{p-1}}{\left\vert
\left\vert \nabla u_{\epsilon}\right\vert ^{2}+\epsilon\right\vert
^{\frac{p-1}{2}}}\left\vert \nabla u\right\vert \\
& \leq & \epsilon^{\frac{p-1}{2}}\int_{\Omega}\left\vert \nabla u\right\vert
dv\\
& \leq & \epsilon^{\frac{p-1}{2}}\left(  vol\left(  \Omega\right)  \right)
^{\frac{p-1}{p}}\left(  \int_{\Omega}\left\vert \nabla u\right\vert
^{p}dv\right)  ^{\frac{1}{p}}\\
& \rightarrow & 0\text{ as }\epsilon\rightarrow0.
\end{array}
\]

Hence we conclude that%
\[%
\begin{array}
[c]{l}%
RHS=RHS1+RHS2\rightarrow0\text{ as }\epsilon\rightarrow0.
\end{array}
\]

\endproof

\bigskip

\subsection{Non-trivial $p$-harmonic function with finite $q$-energy}

In this subsection, we construct an example of non-trivial $p$-harmonic
function $u$ with finite $q$-energy, $q>p-1,$ on a complete noncompact
manifold with weighted Poincar\'{e} inequality $\left(  P_{\rho}\right)  .$

Let $M=%
\mathbb{R}
\times N^{m-1},$ $m\geq3,$ with a metric $ds^{2}=dt^{2}+\eta^{2}\left(
t\right)  g_{N},$ where $\eta\left(  t\right)  :%
\mathbb{R}
\rightarrow\left(  0,\infty\right)  $ is a smooth function with $\eta
^{\prime\prime}>0,$ $\left(  m-2\right)  \left(  \log\eta\right)
^{\prime\prime}+\eta^{-2}Ric_{N}\geq0,$ and $\left(  N,g_{N}\right)  $ is a
compact Riemannian manifold with $vol\left(  N^{m-1}\right)  =1.$

According to \cite{LW3} Proposition 6.1, $M$ satisfies weighted Poincar\'{e}
inequality $\left(  P_{\rho}\right)  $ and $Ric_{M}\geq-\frac{m-1}{m-2}\rho$
with $\rho=\left(  m-2\right)  \eta^{\prime\prime}\eta^{-1}.$

Let $A\left(  t\right)  $ be the volume of $\left\{  t\right\}  \times
N^{m-1},$ then $A\left(  t\right)  =$ $\eta^{m-1}\left(  t\right)  .$

Now we select $\eta\left(  t\right)  $ such that each end of $M$ is
$p$-hyperbolic, and%
\[
A\left(  t\right)  \geq d_{1}\left\vert t\right\vert ^{\frac{p-1}%
{q-p+1-\delta}},\text{ if }\left\vert t\right\vert \geq1,
\]
where $d_{1}>0$ and $0<\delta<q-p+1$ are positive constants.

By using \cite{Tr} Proposition 5.3,
\[%
\begin{array}
[c]{l}%
\emph{Cap}_{p}\left(  \left(  -\infty,a\right)  \times N^{m-1},\left(
b,\infty\right)  \times N^{m-1};M\right)  =\left(  \int_{a}^{b}\left(
\frac{1}{A\left(  t\right)  }\right)  ^{1/\left(  p-1\right)  }dt\right)
^{1-p},
\end{array}
\]
for any $-\infty<a<b<\infty.$ \ If we define $u$ by%
\[%
\begin{array}
[c]{l}%
u\left(  t\right)  =\int_{-\infty}^{t}\left(  \frac{1}{A\left(  s\right)
}\right)  ^{1/\left(  p-1\right)  }ds
\end{array}
\]
then%
\[%
\begin{array}
[c]{l}%
\emph{Cap}_{p}\left(  \left(  -\infty,a\right)  \times N^{n-1},\left(
b,\infty\right)  \times N^{n-1};M\right)  =\left(  u\left(  b\right)
-u\left(  a\right)  \right)  ^{1-p},
\end{array}
\]
$u\left(  t\right)  \rightarrow0$ as $t\rightarrow-\infty$, and $u$ is
uniformly bounded for all $t\in\left(  -\infty,\infty\right)  .$

Moreover, define a function $v$ as follows,%
\[
v\left(  t\right)  =\left\{
\begin{array}
[c]{ll}%
0, & \text{if }t\leq a,\\
\frac{u\left(  t\right)  -u\left(  a\right)  }{u\left(  b\right)  -u\left(
a\right)  },\text{ \ \ } & \text{if }a<t<b,\\
1, & \text{if }t\geq b.
\end{array}
\right.
\]
then%
\[%
\begin{array}
[c]{l}%
\int_{M}\left\vert \nabla v\right\vert ^{p}dv=\int_{a}^{b}\frac{\left(
u^{\prime}\left(  t\right)  \right)  ^{p}}{\left(  u\left(  b\right)
-u\left(  a\right)  \right)  ^{p}}A\left(  t\right)  dt=\left(  u\left(
b\right)  -u\left(  a\right)  \right)  ^{1-p}%
\end{array}
\]
which implies $v$ is extremal of $p$-energy for every $-\infty<a<b<\infty.$
Hence $u\left(  t\right)  $ is $p$-harmonic in $M$ with finite $q$ energy
\[%
\begin{array}
[c]{l}%
\int_{M}\left\vert \nabla u\right\vert ^{q}dv=\int_{-\infty}^{\infty}%
A^{\frac{p-1-q}{p-1}}\left(  t\right)  dt<\infty
\end{array}
\]
for all $q>p-1.$ Moreover, by \cite{LW3} Proposition 6.1, we have
\[%
\begin{array}
[c]{l}%
Ric_{M}\left(  \nabla u,\nabla u\right)  =-\frac{m-1}{m-2}\rho\left\vert
\nabla u\right\vert ^{2}.
\end{array}
\]

\bigskip

\bigskip

\bigskip

\subsection{Volume estimate and $p$-Poincar\'{e} inequality}

In this subsection, we study a complete noncompact manifold $M$ with the
$p$-Poincar\'{e} inequality $\left(  P_{\lambda_{p}}\right)  ,$ $p>1,$ that
is, the inequality
\begin{equation}%
\begin{array}
[c]{l}%
\lambda_{p}\int_{M}|\Psi|^{p}\leq\int_{M}\left\vert \nabla\Psi\right\vert ^{p}%
\end{array}
\label{p poincare}%
\end{equation}
holds for every $\Psi\in W_{0}^{1,p}\left(  M\right)  ,$ where $\lambda_{p} >
0\, .$ In particular, if $p=2,$ this formula is the general Poincar\'{e}
inequality, and $\lambda_{2}$ is the spectrum of $M.$ In \cite{HKM}, they show
that a complete manifold $M$ with positive spectrum $\lambda_{2}>0,$ then it
must have $\lambda_{p}>0$ for all $p\geq2.$ In fact, the following inequality%
\[%
\begin{array}
[c]{l}%
p\left(  \lambda_{p}\right)  ^{1/p}\geq2\left(  \lambda_{2}\right)  ^{1/2}%
\end{array}
\]
holds on $M$ for all $p\geq2.$

\begin{lemma}
\label{decay estimate}Let $M$ be a complete noncompact manifold satisfying
$\left(  P_{\lambda_{p}}\right)  ,$ $p>1$. Suppose $w$ is a positive,
$p$-subharmonic function with a finite $p$-energy on $M.$ If $w$ satisfies
\begin{equation}%
\begin{array}
[c]{l}%
\int_{B\left(  2R\right)  \backslash B\left(  R\right)  }\exp(-\frac{\left(
\lambda_{p}\right)  ^{1/p}r\left(  x\right)  }{p+1})\left\vert w\right\vert
^{p}dv=o\left(  R\right)  ,
\end{array}
\label{growth control}%
\end{equation}
where $R\geq R_{0}+1.$ Then,
\[%
\begin{array}
[c]{l}%
\left(  1-\delta\right)  \left\Vert \exp(\frac{\delta\left(  \lambda
_{p}\right)  ^{1/p}r\left(  x\right)  }{p+1})w\right\Vert _{L_{p}\left(
M\backslash B\left(  R_{0}+1\right)  \right)  }\leq C,
\end{array}
\]
and
\[%
\begin{array}
[c]{l}%
\left(  1-\delta\right)  \left\Vert \exp(\frac{\delta\left(  \lambda
_{p}\right)  ^{1/p}r\left(  x\right)  }{p+1})\nabla w\right\Vert
_{L_{p}\left(  M\backslash B\left(  R_{0}+1\right)  \right)  }\leq C,
\end{array}
\]
for all $0<\delta<1,$ and for some constant $C$ depending on $p$ and
$\lambda_{p}.$
\end{lemma}

\bigskip

\bigskip

\proof Let $\psi$ be a non-negative cut-off function, then we have%
\begin{equation}%
\begin{array}
[c]{lll}%
0 & \geq & \int_{M}\psi^{p}w\left(  -\Delta_{p}w\right) \\
& = & \int_{M}\left\langle \nabla(\psi^{p}w),\left\vert \nabla w\right\vert
^{p-2}\nabla w\right\rangle \\
& = & \int_{M}\psi^{p}\left\vert \nabla w\right\vert ^{p}+pw\left\vert \nabla
w\right\vert ^{p-2}\psi^{p-1}\left\langle \nabla\psi,\nabla w\right\rangle \\
& \geq & \int_{M}\psi^{p}\left\vert \nabla w\right\vert ^{p}-p\int_{M}%
w\psi^{p-1}\left\vert \nabla w\right\vert ^{p-1}\left\vert \nabla
\psi\right\vert .
\end{array}
\label{ed1}%
\end{equation}
By using H\"{o}lder inequality%
\[%
\begin{array}
[c]{l}%
\int_{M}w\psi^{p-1}\left\vert \nabla w\right\vert ^{p-1}\left\vert \nabla
\psi\right\vert \leq\left(  \int_{M}\left\vert \nabla w\right\vert ^{p}%
\psi^{p}\right)  ^{(p-1)/p}\left(  \int_{M}w^{p}\left\vert \nabla
\psi\right\vert ^{p}\right)  ^{1/p},
\end{array}
\]
then (\ref{ed1}) can be rewritten as%
\begin{equation}%
\begin{array}
[c]{l}%
\left\Vert \psi\nabla w\right\Vert _{L_{p}}\leq p\left\Vert \nabla\psi\cdot
w\right\Vert _{L_{p}},
\end{array}
\label{ed2}%
\end{equation}
and this inequality is the Caccioppoli type estimate.

Since Minkowski inequality yields%
\[%
\begin{array}
[c]{l}%
\left\Vert \nabla(\psi w)\right\Vert _{L_{p}}\leq\left\Vert \nabla\psi\cdot
w\right\Vert _{L_{p}}+\left\Vert \psi\nabla w\right\Vert _{L_{p}},
\end{array}
\]
then (\ref{ed2}) implies%
\begin{equation}%
\begin{array}
[c]{l}%
\left\Vert \nabla(\psi w)\right\Vert _{L_{p}}\leq\left(  p+1\right)
\left\Vert \nabla\psi\cdot w\right\Vert _{L_{p}}.
\end{array}
\label{de1}%
\end{equation}
This inequality is not sharp enough whenever $p=2.$ In fact, if $p=2,$ one can
easy to show $\left\Vert \nabla\left(  \psi w\right)  \right\Vert _{L_{2}}%
\leq\left\Vert \nabla\psi w\right\Vert _{L_{2}}$ by the similar method (cf.
\cite{LW1}\cite{LW3}).

By scaling the metric, we may assume $\lambda_{p}=1.$ Combining (\ref{de1})
and (\ref{p poincare}), then
\begin{equation}%
\begin{array}
[c]{l}%
\left\Vert \psi w\right\Vert _{L_{p}}\leq\left(  p+1\right)  \left\Vert
\nabla\psi\cdot w\right\Vert _{L_{p}},
\end{array}
\label{de1-1}%
\end{equation}
where $\psi$ is a cut off function on $M.$

Now we select $\psi=\phi(r(x))\exp(a\left(  r(x)\right)  )$, then%
\begin{equation}%
\begin{array}
[c]{lll}%
\frac{1}{p+1}\left\Vert \psi w\right\Vert _{L_{p}} & \leq & \left\Vert \left(
\nabla\phi+\phi\nabla a\right)  \exp(a(x))w\right\Vert _{L_{p}}\\
& \leq & \left\Vert \left(  \nabla\phi\right)  \exp(a(x))w\right\Vert _{L_{p}%
}+\left\Vert \left(  \nabla a\right)  \phi\exp(a(x))w\right\Vert _{L_{p}}%
\end{array}
\label{de11}%
\end{equation}
where $\phi$ is a non-negative cut-off function defined by $\phi=\phi_{+}%
+\phi_{-}$ where%
\[%
\begin{array}
[c]{ll}%
\phi_{+}(r)=\left\{
\begin{array}
[c]{ll}%
r-R_{0} & \text{for }R_{0}\leq r\leq R_{0}+1,\vspace{3mm}\\
1 & \text{for }r>R_{0}+1,\vspace{3mm}%
\end{array}
\right.  & \phi_{-}(r)=\left\{
\begin{array}
[c]{ll}%
\frac{R-r}{R} & \text{for }R\leq r\leq2R,\vspace{3mm}\\
-1 & \text{for }r>2R,
\end{array}
\right.
\end{array}
\]
and we also choose $a=a_{+}(r(x))+a_{-}(r(x))$ as%
\[%
\begin{array}
[c]{l}%
a_{+}(r)=\left\{
\begin{array}
[c]{ll}%
\frac{\delta r\left(  x\right)  }{p+1} & \text{for }r\leq\frac{K}{1+\delta
},\vspace{3mm}\\
\frac{\delta K}{\left(  1+\delta\right)  \left(  p+1\right)  } & \text{for
}r>\frac{K}{1+\delta},
\end{array}
\right. \\
a_{-}(r)=\left\{
\begin{array}
[c]{ll}%
0 & \text{for }r\leq\frac{K}{1+\delta},\vspace{3mm}\\
\frac{1}{p+1}\left(  \frac{2K}{1+\delta}-r\left(  x\right)  \right)  &
\text{for }r>\frac{K}{1+\delta},
\end{array}
\right.
\end{array}
\]
for some fixed $K>\left(  R_{0}+1\right)  \left(  1+\delta\right)  ,$
$0<\delta<1,$ and $R\geq\frac{K}{1+\delta},$ it's easy to check that%
\[%
\begin{array}
[c]{l}%
\left\vert \nabla\phi\right\vert ^{2}\left(  x\right)  =\left\{
\begin{array}
[c]{ll}%
1\text{ \ \ \ \ } & \text{on }B(R_{0}+1)\backslash B(R_{0}),\vspace{3mm}\\
0 & \text{on }B\left(  R_{0}\right)  ,\text{ }B\left(  R\right)  \backslash
B(R_{0}+1)\text{ and }M\backslash B\left(  2R\right)  ,\vspace{3mm}\\
\frac{1}{R^{2}} & \text{on }B\left(  2R\right)  \backslash B\left(  R\right)
,\vspace{3mm}%
\end{array}
\right.
\end{array}
\]
and%
\[%
\begin{array}
[c]{l}%
\left\vert \nabla a\right\vert ^{2}\left(  x\right)  =\left\{
\begin{array}
[c]{ll}%
\frac{\delta^{2}}{\left(  p+1\right)  ^{2}} & \text{for }r<\frac{K}{1+\delta
},\vspace{3mm}\\
\frac{1}{\left(  p+1\right)  ^{2}}\text{\ \ \ } & \text{for }r>\frac
{K}{1+\delta}.
\end{array}
\right.
\end{array}
\]
Substituting into (\ref{de11}), we obtain%
\[%
\begin{array}
[c]{lll}
&  & \frac{1}{p+1}\left\Vert \phi\exp(a(x))w\right\Vert _{L_{p}\left(
M\right)  }\\
& \leq & \left\Vert \left(  \nabla\phi_{+}\right)  \exp(a(x))w\right\Vert
_{L_{p}\left(  M\right)  }+\left\Vert \left(  \nabla\phi_{-}\right)
\exp(a(x))w\right\Vert _{L_{p}\left(  M\right)  }\\
&  & +\left\Vert \left(  \nabla a_{+}\right)  \phi\exp(a(x))w\right\Vert
_{L_{p}\left(  M\right)  }+\left\Vert \left(  \nabla a_{-}\right)  \phi
\exp(a(x))w\right\Vert _{L_{p}\left(  M\right)  }\\
& \leq & \left\Vert \exp(a(x))w\right\Vert _{L_{p}\left(  B(R_{0}+1)\backslash
B(R_{0})\right)  }+\frac{1}{R}\left\Vert \exp(a(x))w\right\Vert _{L_{p}\left(
B\left(  2R\right)  \backslash B\left(  R\right)  \right)  }\\
&  & +\frac{\delta}{p+1}\left\Vert \phi\exp(a(x))w\right\Vert _{L_{p}\left(
B\left(  \frac{K}{1+\delta}\right)  \right)  }+\frac{1}{p+1}\left\Vert
\phi\exp(a(x))w\right\Vert _{L_{p}\left(  M\backslash B\left(  \frac
{K}{1+\delta}\right)  \right)  },
\end{array}
\]
hence%
\[%
\begin{array}
[c]{lll}
&  & \left(  \frac{1-\delta}{p+1}\right)  \left\Vert \phi\exp
(a(x))w\right\Vert _{L_{p}\left(  B\left(  \frac{K}{1+\delta}\right)
\backslash B\left(  R_{0}+1\right)  \right)  }\\
& \leq & \left\Vert \exp(a(x))w\right\Vert _{L_{p}\left(  B(R_{0}+1)\backslash
B(R_{0})\right)  }+\frac{1}{R}\left\Vert \exp(a(x))w\right\Vert _{L_{p}\left(
B\left(  2R\right)  \backslash B\left(  R\right)  \right)  }.
\end{array}
\]
The definition of $a(x)$ and the growth condition (\ref{growth control}) imply
that the last term on the right hand side tends to $0$ as $R\rightarrow\infty
$. Thus one has the following inequality,%
\begin{equation}%
\begin{array}
[c]{cl}
& \left(  \frac{1-\delta}{p+1}\right)  \left\Vert \exp(a(x))w\right\Vert
_{L_{p}\left(  B\left(  \frac{K}{1+\delta}\right)  \backslash B\left(
R_{0}+1\right)  \right)  }\\
\leq & \left\Vert \exp(a(x))w\right\Vert _{L_{p}\left(  B(R_{0}+1)\backslash
B(R_{0})\right)  }.
\end{array}
\label{ed3-0}%
\end{equation}
Since the right hand side of (\ref{ed3-0}) is independent of $K$ and
$0<\delta<1,$ by letting $K\rightarrow\infty$ we obtain that
\begin{equation}%
\begin{array}
[c]{l}%
\left(  1-\delta\right)  \left\Vert \exp(a(x))w\right\Vert _{L_{p}\left(
M\backslash B\left(  R_{0}+1\right)  \right)  }\leq C_{1},
\end{array}
\label{ed3}%
\end{equation}
for some constant $0<C_{1}=C_{1}\left(  p\right)  <\infty.$

Moreover, by (\ref{ed2}) and similar process as above, we have%
\[%
\begin{array}
[c]{ll}
& \frac{1}{p}\left\Vert \psi\nabla w\right\Vert _{L_{p}\left(  M\right)  }\\
\leq & \left\Vert \nabla\psi\cdot w\right\Vert _{L_{p}\left(  M\right)  }\\
\leq & \left\Vert \exp(a(x))w\right\Vert _{L_{p}\left(  B(R_{0}+1)\backslash
B(R_{0})\right)  }+\frac{1}{R}\left\Vert \exp(a(x))w\right\Vert _{L_{p}\left(
B\left(  2R\right)  \backslash B\left(  R\right)  \right)  }\\
& +\frac{\delta}{p+1}\left\Vert \phi\exp(a(x))w\right\Vert _{L_{p}\left(
B\left(  \frac{K}{1+\delta}\right)  \right)  }+\frac{1}{p+1}\left\Vert
\phi\exp(a(x))w\right\Vert _{L_{p}\left(  B\left(  2R\right)  \backslash
B\left(  \frac{K}{1+\delta}\right)  \right)  }\\
\leq & 2\left\Vert \exp(a(x))w\right\Vert _{L_{p}\left(  B(R_{0}+1)\backslash
B(R_{0})\right)  }+3\left\Vert \phi\exp(a(x))w\right\Vert _{L_{p}\left(
B\left(  2R\right)  \backslash B\left(  R_{0}+1\right)  \right)  }\\
\leq & C_{2}+\frac{3C_{1}}{1-\delta}.
\end{array}
\]
Hence, by letting $R\rightarrow\infty$ and then letting $K\rightarrow\infty,$
we conclude%
\[%
\begin{array}
[c]{l}%
\left(  1-\delta\right)  \left\Vert \exp(\delta r\left(  x\right)  )\nabla
w\right\Vert _{L_{p}\left(  M\backslash B\left(  R_{0}+1\right)  \right)
}\leq C_{3}%
\end{array}
\]
for some constant $0<C_{3}=C_{2}+\frac{3C_{1}}{1-\delta}<\infty.$

Then lemma now follows.

\endproof

\bigskip

\bigskip

\begin{lemma}
\label{decay estimate on E}Let $M$ be a complete noncompact manifold
satisfying $\left(  P_{\lambda_{p}}\right)  ,$ $p>1$. Suppose $E$ is an end of
$M$ respective to a compact set, $w_{i}$ is a positive, $p$-harmonic function
with a finite $p$-energy on $E\left(  R_{i}\right)  $ and $w_{i}=1$ on
$\partial E$ and $w_{i}=0$ on $S\left(  R_{i}\right)  =\partial E\left(
R_{i}\right)  \backslash\partial E.$ If $R_{i}\rightarrow\infty$ and
$w_{i}\rightarrow w$ as $i\rightarrow\infty.$ Then,
\begin{equation}%
\begin{array}
[c]{l}%
\int_{E\backslash E\left(  R\right)  }\left\vert \nabla w\right\vert
^{p}dv\leq C_{3}R^{p}\exp(\frac{-\left(  \lambda_{p}\right)  ^{1/p}\left(
R-1\right)  }{\left(  p+1\right)  }),
\end{array}
\label{DDw-2}%
\end{equation}
and
\begin{equation}%
\begin{array}
[c]{l}%
\int_{E\left(  kR\right)  \backslash E\left(  R\right)  }\left\vert
w\right\vert ^{p}dv\leq C_{1}R^{p}\exp(\frac{-\left(  \lambda_{p}\right)
^{1/p}\left(  R-1\right)  }{p+1}),
\end{array}
\label{Dw-2}%
\end{equation}
for some constant $C$ depending on $p.$
\end{lemma}

\bigskip

\proof As in the proof of Lemma \ref{decay estimate}. If $\phi$ is a
non-negative cut-off function defined by%
\[%
\begin{array}
[c]{l}%
\phi(r(x))=\left\{
\begin{array}
[c]{ll}%
r(x)-R_{0}\text{ \ \ \ \ \ \ \ \ } & \text{on }E(R_{0}+1)\backslash
E(R_{0}),\vspace{3mm}\\
1 & \text{on }E\backslash E(R_{0}+1),\vspace{3mm}%
\end{array}
\right.
\end{array}
\]
$\newline$and we choose $a=\frac{\delta r\left(  x\right)  }{p+1}$ for
$0<\delta<1.$ It's easy to check that%
\[%
\begin{array}
[c]{lll}%
\left\vert \nabla\phi\right\vert ^{2}\left(  x\right)  =\left\{
\begin{array}
[c]{ll}%
1\text{\ \ \ } & \text{on }E(R_{0}+1)\backslash E(R_{0}),\vspace{3mm}\\
0 & \text{on }E\backslash E\left(  R_{0}+1\right)  ,
\end{array}
\right.  & \text{and} & \left\vert \nabla a\right\vert ^{2}\left(  x\right)
=\frac{\delta^{2}}{\left(  p+1\right)  ^{2}}.
\end{array}
\]
By the formula (\ref{de11}), we obtain%
\[%
\begin{array}
[c]{ll}
& \frac{1}{p+1}\left\Vert \phi\exp(a(x))w\right\Vert _{L_{p}}\\
\leq & \left\Vert \left(  \nabla\phi\right)  \exp(a(x))w\right\Vert _{L_{p}%
}+\left\Vert \left(  \nabla a\right)  \phi\exp(a(x))w\right\Vert _{L_{p}}\\
\leq & \left\Vert \exp(a(x))w\right\Vert _{L_{p}\left(  E\left(
R_{0}+1\right)  \backslash E\left(  R_{0}\right)  \right)  }+\frac{\delta
}{p+1}\left\Vert \phi\exp(a(x))w\right\Vert _{L_{p}\left(  E\right)  }%
\end{array}
\]
hence%
\[%
\begin{array}
[c]{l}%
\left(  \frac{1-\delta}{p+1}\right)  \left\Vert \phi\exp(a(x))w\right\Vert
_{L_{p}\left(  E\backslash E\left(  R_{0}+1\right)  \right)  }\leq\left\Vert
\exp(a(x))w\right\Vert _{L_{p}\left(  E(R_{0}+1)\backslash E(R_{0})\right)  }.
\end{array}
\]
Then we obtain that
\begin{equation}%
\begin{array}
[c]{l}%
\left(  1-\delta\right)  \left\Vert \exp(\delta r)w\right\Vert _{L_{p}\left(
E\backslash E\left(  R_{0}+1\right)  \right)  }\leq C_{1},
\end{array}
\label{DDw1}%
\end{equation}
for some constant $0<C_{1}=C_{1}\left(  p\right)  <\infty.$

Moreover, since%
\[%
\begin{array}
[c]{lll}%
\frac{1}{p}\left\Vert \psi\nabla w\right\Vert _{L_{p}} & \leq & \left\Vert
\nabla\psi w\right\Vert _{L_{p}}\\
& \leq & \left\Vert \exp(a(x))w\right\Vert _{L_{p}\left(  E(R_{0}+1)\backslash
E(R_{0})\right)  }+\delta\left\Vert \phi\exp(a(x))w\right\Vert _{L_{p}\left(
E\right)  }\\
& \leq & C_{2}+\frac{\delta C_{1}}{1-\delta}.
\end{array}
\]
Hence, we conclude%
\begin{equation}%
\begin{array}
[c]{l}%
\left(  1-\delta\right)  \left\Vert \exp(\delta r\left(  x\right)  )\nabla
w\right\Vert _{L_{p}\left(  E\backslash E\left(  R_{0}+1\right)  \right)
}\leq C_{3}%
\end{array}
\label{DDw2}%
\end{equation}
for some constant $0<C_{3}=C_{3}\left(  p\right)  <\infty.$

If we select $\delta=\left(  1-\frac{1}{R}\right)  $ and $R_{0}>1,$ then
\ref{DDw1}) gives%
\[%
\begin{array}
[c]{lll}%
C_{3} & \geq & \frac{1}{R^{p}}\int_{E\backslash E\left(  R_{0}+1\right)  }%
\exp\left(  (1-\frac{1}{R})\frac{\left(  \lambda_{p}\right)  ^{1/p}r}%
{p+1}\right)  \left\vert \nabla w\right\vert ^{p}dv\\
& \geq & \frac{1}{R^{p}}\int_{E\left(  kR\right)  \backslash E\left(
R_{0}+1\right)  }\exp\left(  (1-\frac{1}{R})\frac{\left(  \lambda_{p}\right)
^{1/p}r}{p+1}\right)  \left\vert \nabla w\right\vert ^{p}dv.
\end{array}
\]
Hence%
\[%
\begin{array}
[c]{l}%
\int_{E\left(  kR\right)  }\exp(\frac{\left(  \lambda_{p}\right)
^{1/p}\left(  R-1\right)  r}{\left(  p+1\right)  R})\left\vert \nabla
w\right\vert ^{p}dv\leq C_{3}R^{p},
\end{array}
\]
and then we have%
\[%
\begin{array}
[c]{l}%
\int_{E\left(  kR\right)  \backslash E\left(  R\right)  }\left\vert \nabla
w\right\vert ^{p}dv\leq C_{3}R^{p}\exp(\frac{-\left(  \lambda_{p}\right)
^{1/p}\left(  R-1\right)  }{\left(  p+1\right)  }),
\end{array}
\]
for all constant $k>1.$

Similarly, \ref{DDw2} implies%
\[%
\begin{array}
[c]{l}%
\int_{E\left(  kR\right)  }\exp(\frac{\left(  \lambda_{p}\right)
^{1/p}\left(  R-1\right)  r}{\left(  p+1\right)  R})\left\vert w\right\vert
^{p}dv\leq C_{1}R^{p}%
\end{array}
\]
and%
\[%
\begin{array}
[c]{l}%
\int_{E\left(  kR\right)  \backslash E\left(  R\right)  }\left\vert
w\right\vert ^{p}dv\leq C_{1}R^{p}\exp(\frac{-\left(  \lambda_{p}\right)
^{1/p}\left(  R-1\right)  }{p+1}),
\end{array}
\]
for any constant $k>1.$

\endproof

\bigskip

\begin{lemma}
\label{Volume estimate}Let $M$ be a complete noncompact manifold satisfying
$\left(  P_{\lambda_{p}}\right)  ,$ $p>1$. If $E$ is a $p$-hyperbolic end of
$M^{n}$, then
\[%
\begin{array}
[c]{l}%
V\left(  E(R+1)\right)  -V\left(  E\left(  R\right)  \right)  \geq
CR^{-p\left(  p-1\right)  }\exp(\frac{\left(  p-1\right)  \left(  \lambda
_{p}\right)  ^{1/p}\left(  R-1\right)  }{p+1}).
\end{array}
\]
for some constant $C>0,$ and for $R$ sufficiently large. If $E$ is
$p$-parabolic, then%
\[%
\begin{array}
[c]{l}%
V\left(  E\right)  <\infty
\end{array}
\]
and%
\[%
\begin{array}
[c]{l}%
V\left(  E\right)  -V\left(  E\left(  R\right)  \right)  \leq CR^{p}\exp
(\frac{-\left(  \lambda_{p}\right)  ^{1/p}\left(  R-1\right)  }{p+1})
\end{array}
\]
for some constant $C>0,$ for any $0<\delta<1,$ and for $R$ sufficiently large.
\end{lemma}

\bigskip

\proof If $E$ is $p$-parabolic, we select the barrier function $w=1$ on $E,$
then (\ref{Dw-2}) implies%
\[%
\begin{array}
[c]{l}%
\int_{E\backslash E\left(  R\right)  }dv\leq CR^{p}\exp(\frac{-\left(
\lambda_{p}\right)  ^{1/p}\left(  R-1\right)  }{p+1})
\end{array}
\]
for all $R$ large enough and for any $\delta$ satisfying $0<\delta<1.$ This
implies $V\left(  E\right)  <\infty.$

If $E$ is $p$-hyperbolic. Let $w$ be the barrier function on $E$, and
$S\left(  R\right)  =\partial E\left(  R\right)  \backslash\partial E,$ then%
\begin{equation}%
\begin{array}
[c]{lll}%
C & = & \int_{\partial E}\left\vert \nabla w\right\vert ^{p-2}\frac{\partial
w}{\partial\nu}dA\\
& \leq & \int_{S\left(  r\right)  }\left\vert \nabla w\right\vert ^{p-1}dA\\
& \leq & \left(  \int_{S\left(  r\right)  }\left\vert \nabla w\right\vert
^{p}dA\right)  ^{\left(  p-1\right)  /p}\left(  \int_{S\left(  r\right)
}dA\right)  ^{1/p}.
\end{array}
\label{13}%
\end{equation}
Then (\ref{13}) imply%
\[%
\begin{array}
[c]{lll}%
\int_{R}^{R+1}\left(  \int_{S\left(  r\right)  }dA\right)  ^{-1/\left(
p-1\right)  }dr & \leq & C\int_{R}^{R+1}\int_{S\left(  r\right)  }\left\vert
\nabla w\right\vert ^{p}dAdr\\
& = & C\int_{E\left(  R+1\right)  \backslash E\left(  R\right)  }\left\vert
\nabla w\right\vert ^{p}dv.
\end{array}
\]
By using Schwarz inequality,%
\[%
\begin{array}
[c]{lll}%
1 & = & \int_{R}^{R+1}\left(  \int_{S\left(  r\right)  }dA\right)  ^{-\frac
{1}{p}}\left(  \int_{S\left(  r\right)  }dA\right)  ^{\frac{1}{p}}dr\\
& \leq & \left(  \int_{R}^{R+1}\left(  \int_{S\left(  r\right)  }dA\right)
^{-\frac{1}{p-1}}dr\right)  ^{\frac{p-1}{p}}\cdot\left(  \int_{R}^{R+1}%
\int_{S\left(  r\right)  }dAdr\right)  ^{\frac{1}{p}}\\
& \leq & C\left(  \int_{E\left(  R+1\right)  \backslash E\left(  R\right)
}\left\vert \nabla w\right\vert ^{p}dv\right)  ^{\frac{p-1}{p}}\cdot\left(
\int_{R}^{R+1}\int_{S\left(  R\right)  }dAdr\right)  ^{\frac{1}{p}}.
\end{array}
\]
Then co-area formula and (\ref{DDw-2}) give%
\[%
\begin{array}
[c]{l}%
\int_{E\left(  R+1\right)  \backslash E\left(  R\right)  }dv\geq CR^{-p\left(
p-1\right)  }\exp(\frac{\left(  p-1\right)  \left(  \lambda_{p}\right)
^{1/p}\left(  R-1\right)  }{p+1}).
\end{array}
\]

\endproof

\bigskip

\bigskip

Since $\left(  P_{\lambda_{p}}\right)  $ implies the volume of $M$ is
infinity, then Lemma \ref{Volume estimate} implies the following property.

\begin{theorem}
\label{hy}If $M$ is a complete noncompact manifold satisfying $\left(
P_{\lambda_{p}}\right)  ,$ then $M$ must be $p$-hyperbolic.
\end{theorem}

\begin{remark}
One can also prove the above theorem by contradiction. That is, if $M$ were
$p$-parabolic, then $\lambda_{p}\, $ would be zero, a contradiction by a
different approach (cf. e.g. \cite{WLW} proof of Theorem 6.1).
\end{remark}

\textbf{Acknowledgments. }S.C. Chang and J.T. Chen were partially supported by
NSC, and S.W. Wei was partially supported by NSF (DMS-1447008), and the OU
Arts and Sciences Travel Assistance Program Fund. Part of this paper has been
written while the second author was visiting University of Oklahoma in summer
2011. He is grateful to Professor S.W. Wei and to the Department of
Mathematics for the kind hospitality. The authors would like to express their
thanks to the referee for valuable comments.

\bigskip
\end{proof}
\end{proof}

\end{document}